\documentclass[11pt,reqno]{amsart}

\usepackage{mathrsfs}
\usepackage{amsfonts}
\usepackage{graphicx,epsfig,amsmath,amsfonts}
\usepackage{color}
\usepackage{latexsym,amssymb}
\usepackage{multirow}
\usepackage{bm}
\usepackage{float}
\usepackage{subfigure}
\usepackage{placeins}

\usepackage [latin1] {inputenc}

\usepackage{bm}

\makeatletter
\renewcommand{\@thesubfigure}{\hskip\subfiglabelskip}
\makeatother
\usepackage{subeqnarray}
\usepackage{cases}
\usepackage{color}
\usepackage{caption}
\newtheorem{theorem}{\hskip\parindent\bf Theorem}[section]

\newtheorem{remark}{\hskip\parindent\bf Remark}[section]
\newtheorem{example}{\hskip\parindent\bf Example}[section]
\newtheorem{algorithm}{\hskip\parindent\bf Algorithm}[section]

\def\bc{\begin{center}}
\def\ec{\end{center}}
\numberwithin{equation}{section}{}

\newcommand{\ba}{\begin{array}}\newcommand{\ea}{\end{array}}
\newcommand{\be}{\begin{eqnarray}}\newcommand{\ee}{\end{eqnarray}}
\newcommand{\bex}{\begin{eqnarray*}}
\newcommand{\eex}{\end{eqnarray*}}
\begin{document}
\graphicspath{{./}}

\title[Linear adaptive variable time step schemes]
{Linear adaptive variable time step schemes for the anisotropic phase field dendritic crystal growth model}
\author{Yayu Guo$^{1}$, Ziqi Ma$^{2,*}$}

\thanks{\hskip -12pt
${}^{1}$School of Science, Changzhou Institute of Technology, 213000 Changzhou, China.
 \\
 ${}^{2}$AVIC Leihua Electronic Technology Research Institute, 214063 Wuxi, China. \\
${}^{*}$ Corresponding author.\\
Emails:
yayuguo@163.com (Yayu Guo), ziqima05@163.com (Ziqi Ma).\\}


\date {\today}
\maketitle

\begin{abstract}
In this paper, we propose a class of efficient, linear and decoupled
adaptive variable time step schemes
for the anisotropic phase field dendritic crystal growth model.
The main difficulty is twofold: the gradient dependent anisotropic coefficient and the strongly nonlinear terms make unconditional energy stability difficult to achieve, and on nonuniform meshes the discrete energy dissipation law, which is straightforward to establish on uniform meshes, has to be reestablished under varying step sizes. To overcome these difficulties, two linear stabilizers and a time dependent auxiliary variable with relaxation are introduced and incorporated into a second-order generalized BDF-type operator, with the spatial derivatives discretized by the Fourier spectral method.
The proposed schemes are  energy stable for arbitrary positive adjacent step ratios, so that no restriction on the step ratio is imposed, and only two linear equations together with
a few algebraic equations, which are solved independently in Fourier space, need to be solved at each time step. Moreover, the
schemes are naturally compatible with energy based time adaptive strategies. The associated adaptive strategy reduces the number of time steps by $24.6\%$--$86.6\%$ in comparison with uniform time stepping, and numerical experiments on fourfold and sixfold dendritic growth in two dimensions and on three-dimensional growth confirm the accuracy, efficiency and stability of the proposed schemes.
\end{abstract}

\section{Introduction}
Dendrites are the most common growth patterns during the solidification of metals, alloys and other substances. Inspired by the snowflake structure, a fairly complete theory of dendrite growth has been established at the microscopic scale. For metallic materials, this microscopic process greatly affects the material properties, and therefore dendritic crystal growth has become an important research topic \cite{chalmers1970principles,glicksman1989fundamentals}. Experimental investigations of dendrite growth are difficult and costly, while numerical simulation provides a powerful alternative. Among various numerical approaches, the phase field method \cite{fix1982phase,caginalp1986analysis,kim1999universal} has the most mature theoretical basis and produces the most accurate simulation results: a phase field variable is introduced to distinguish the solid and liquid states in space and time, so that the explicit tracking of the complex moving interface is avoided.

To observe dendritic crystal growth, Karma and Rappel developed a quantitative phase field formulation for the solidification of a pure melt \cite{karma1998quantitative}.
Other related simulations of the dendritic crystal growth model can be found in \cite{kobayashi1993modeling,karma1999phase,ramirez2004phase,ferreira2011numerical}. The dendritic crystal growth model is composed of an Allen Cahn equation with a gradient dependent anisotropy coefficient coupled with a heat transfer equation. It preserves an energy dissipation law and meanwhile is strongly nonlinear, so designing efficient numerical schemes for this model is the main problem we want to solve.

In recent years, many numerical schemes have been constructed to preserve the energy dissipation property of gradient flow models; see, e.g., \cite{zhu1999coarsening,xu2006stability,shen2010numerical,feng2013stabilized,yang2016linear,yang2017linearly}. For the nonlinear terms in the model, there are several linearization techniques, such as:
 the Lagrange multiplier approach \cite{cheng2020global},
 the invariant energy quadratization (IEQ) method \cite{yang2017numerical,zhao2017numerical}, the scalar auxiliary variable (SAV) method \cite{shen2018scalar}, and some improved auxiliary variable methods (such as the generalized SAV method \cite{shen2018scalar,huang2022new}, the relaxed SAV method \cite{jiang2022improving}, and the generalized SAV method with relaxation \cite{zhang2022generalized}); another influential variant can be found in \cite{hou2019variant}. There are also many works devoted to the phase field dendritic crystal growth model; see, e.g., \cite{li2011fast,yang2016linear,wang2022accurate,li2023second,wang2025class,zhou2026high,li2026imex,yao2026semi}. Although these schemes are efficient for solving the model, it is still desirable to design more cost effective, linear, decoupled and higher order schemes.

The difficulties of solving the dendritic crystal growth model lie in the strong nonlinearity, the gradient dependent anisotropic coefficient, and the coupling between the phase field and the temperature field.
In order to overcome these difficulties, the following methods have been proposed:
Zhao et al. \cite{zhao2017numerical} considered the dendrite growth model in the isothermal case, and proposed a linear, decoupled numerical scheme by employing the invariant energy quadratization (IEQ) approach, in which the auxiliary variable depends on space. Then, to balance the anisotropic coefficient and the nonlinear term, an IEQ method with two stabilization terms was proposed in \cite{yang2019efficient}, but the scheme is coupled, which means one needs to solve a system that couples the phase field and the temperature field at each time step, and this undoubtedly increases the computational complexity.
To remedy the drawbacks of the above schemes, two aspects have been considered: one is to completely decouple the phase field and the temperature field, and the other is to improve the form of the auxiliary variables so that they only depend on time.
For the first aspect, a splitting method based on the temperature field was proposed in \cite{zhang2019novel}, which is feasible in practical calculations, but the unconditional stability of its second-order version is not available.
For the second aspect, the scalar auxiliary variable approach \cite{yang2020efficient,yang2021fully} has been adopted.
Recently, an auxiliary variable satisfying the "zero-energy-contribution" property was proposed by Yang \cite{yang2021fully} and applied to simulate this model; the resulting scheme is decoupled, linear, and unconditionally energy stable. In addition, a more efficient numerical method which only needs to solve three linear equations was proposed in \cite{li2022new}. Recently, Wang and Xu \cite{wang2025class} proposed a class of efficient high order time stepping methods for the anisotropic phase field dendritic crystal growth model.
To the best of our knowledge, \cite{yang2021fully} needs to solve four linear equations at each time step, while \cite{li2022new} needs to solve three linear equations.
Most recently, Guo \cite{guo2024efficient} further improved the efficiency by combining the relaxed generalized scalar auxiliary variable (R-GSAV) approach with the Fourier spectral method, so that only two linear equations need to be solved at each time step, which is the closest work to the present manuscript. Nevertheless, the proposed scheme is still implemented on uniform temporal meshes.

Besides these efforts, two recent studies on this model are also closely related to the present work. Zhou and Guo \cite{zhou2026high} proposed high order efficient and energy stable Runge--Kutta scalar auxiliary variable/local discontinuous Galerkin (RK--SAV/LDG) schemes for this model, and Li et al. \cite{li2026imex} developed efficient and decoupled second- and third-order implicit--explicit Runge--Kutta schemes with original-form energy stability for the same model. These works confirm that constructing efficient and energy stable schemes with higher temporal accuracy for this model is still an active research topic, while the present work focuses on nonuniform temporal meshes and energy based time adaptivity.

In addition, it is worth noting that the above mentioned schemes for the dendritic crystal growth model are implemented on uniform temporal meshes. Since the solidification process usually involves multiple time scales, i.e., the initial rapid growth stage with fast energy decay and the later slowly evolving stage, a uniform small time step leads to a huge waste of computational cost for long time simulations, while a uniform large time step may lose accuracy at the early stage. Variable time stepping techniques have been shown to be effective for gradient flow models; see \cite{hou2023implicit,zhang2025energy,shen2026stability,liao2022adaptive}. Therefore, it is desirable to combine the R-GSAV approach with variable time steps for the dendritic crystal growth model, which, to the best of our knowledge, has not been studied in the literature.

The aim of this paper is to construct  energy stable, efficient and high precision schemes for the anisotropic phase field dendritic crystal growth model. We construct a temporal discretization framework on nonuniform temporal meshes for this model. The main difficulty lies in that the discrete energy dissipation law, which is straightforward to establish on uniform temporal meshes, has to be re-established under varying step sizes. To overcome this difficulty, we combine the R-GSAV approach \cite{zhang2022generalized} with the second-order generalized BDF operator \cite{hou2023implicit}, prove the discrete energy dissipation law for arbitrary positive adjacent step ratios, and retain the efficiency of only two linear solves at each time step.
The main contributions are as follows:
\begin{itemize}
  \item We construct the second-order generalized SAV scheme with relaxation on nonuniform temporal meshes for solving the dendritic crystal growth model. The proposed scheme is energy stable, and the phase field and the temperature field are decoupled, so that one only needs to solve two linear equations and some algebraic equations at each time step, which are pointwise solves in Fourier space and hence cost the same on uniform and nonuniform meshes;
  \item We give a rigorous proof of the energy stability of the proposed scheme  for arbitrary positive adjacent step ratios, whereas the variable-step BDF2 analysis in \cite{hou2023implicit} requires an additional bound $\gamma_{n+1}\leq4.8645$;
  \item Moreover, the scheme allows variable time steps and is naturally compatible with energy based time adaptive strategies, and the modified energy obtained by the proposed scheme is equal to the original energy almost every time in our numerical results, i.e., the scheme exhibits original energy stability in all the tests performed. The associated adaptive strategy reduces the number of time steps by $24.6\%$--$86.6\%$ in the two and three dimensional simulations.
\end{itemize}

The rest of this paper is organized as follows. In Section 2, we briefly introduce the phase field dendritic crystal growth model and its energy dissipation law in the continuous case. In Section 3, we construct the R-GSAV scheme based on the nonuniform BDF2 operator for the model, and prove the corresponding energy stability. Numerical tests verifying the effectiveness and accuracy of the proposed scheme, including an energy based time adaptive strategy, are provided in Section 4. Some concluding remarks are given in Section 5.

\section{The model}
In this section, we briefly introduce the phase field dendrite crystal model. Consider the following energy functional as \cite{karma1998quantitative}
\begin{equation}\label{energy}
E_{tot}(\phi, U)=\int_{\Omega}\left(\frac{1}{2}m^2(\nabla \phi)| \nabla \phi|^2
+\frac{\lambda}{2 \varepsilon K} U^2+ F(\phi)\right) d x,
\end{equation}
where $\phi$ stands for the phase field, where $\phi=1$ for the solid and $\phi=-1$ for the fluid, $U$ denotes the temperature field, $m(\nabla \phi)$ is the anisotropic coefficient, $\varepsilon$, $\lambda$ and $K$ are positive constants. $F(\phi)=\frac{1}{4 \varepsilon^2}(\phi^2-1)^2$ is a double-well type Ginzburg-Landau potential.
In 2D, the anisotropic coefficient $m(\nabla\phi)$ is denoted as follows
\begin{equation}\label{1.1}
m(\nabla\phi)=1+\varepsilon_4\cos(\beta\Theta),
\end{equation}
where $\beta$ is the number of folds of anisotropy, $\varepsilon_4$ is the anisotropy strength, $\Theta=\arctan(\frac{\phi_y}{\phi_x})$.
 For $\beta=4$, the anisotropy coefficient can be rewritten as
\begin{equation}\label{1.2}
m(\nabla\phi):=\left\{
             \begin{array}{ll}
(1-3\varepsilon_4)(1+\frac{4\varepsilon_4}{1-3\varepsilon_4}\frac{\phi_x^4+\phi_y^4}{|\nabla\phi|^4}),   & \text{in 2D},\\
(1-3\varepsilon_4)(1+\frac{4\varepsilon_4}{1-3\varepsilon_4}\frac{\phi_x^4+\phi_y^4+\phi_z^4}{|\nabla\phi|^4}),    & \text{in 3D}.
             \end{array}
           \right.
\end{equation}
Denote by $\boldsymbol{H}(\phi)$ the variational derivative of $m(\nabla\phi)$, which can be expressed in the following form
\begin{equation}\label{1.3}
 \boldsymbol{H}(\phi):=\left\{
             \begin{array}{ll}
4\varepsilon_4\frac{4}{|\nabla\phi|^6}(\phi_x(\phi_x^2\phi_y^2-\phi_y^4),\phi_y(\phi_x^2\phi_y^2-\phi_x^4)),   & \text{in 2D},\\
4\varepsilon_4\frac{4}{|\nabla\phi|^6}(\phi_x(\phi_x^2\phi_y^2+\phi_x^2\phi_z^2-\phi_y^4-\phi_z^4),\phi_y(\phi_y^2\phi_z^2+\phi_x^2\phi_y^2-\phi_x^4-\phi_z^4),
 \\ \phi_z(\phi_x^2\phi_z^2+\phi_y^2\phi_z^2-\phi_x^4-\phi_y^4)),    & \text{in 3D}.
             \end{array}
           \right.
\end{equation}
The governing equations of the phase field dendrite crystal model, derived by the variational approach, are given as the following formulation:
\begin{subequations}\label{1.4}
\begin{align}
\sigma\phi_t=
&-\frac{\delta E_{tot}}{\delta\phi}-\frac{\lambda}{\varepsilon} g(\phi)U \nonumber\\
=&\nabla \cdot(m^2(\nabla \phi) \nabla \phi+|\nabla \phi|^2 m(\nabla \phi) \boldsymbol{H}(\phi))-f(\phi)-\frac{\lambda}{\varepsilon} g(\phi)U, \label{1.4a}\\
U_t=&D\Delta U+Kg(\phi)\phi_t, \label{1.4b}
\end{align}
\end{subequations}
where $\sigma$ is the mobility constant, $D$ is the diffusion rate of the temperature, the function $g(\phi)=(1-\phi^2)^2$, $f(\phi)=F'(\phi)$.
To avoid the complexity of integrating on the boundary, we adopt periodic or Neumann boundary conditions for the system \eqref{1.4}.
Taking the $L^2$ inner product of the system \eqref{1.4a} and \eqref{1.4b} with $\phi_t$ and $\frac{\lambda}{\varepsilon K}U$ respectively, then the energy dissipation laws can be obtained as follows
\begin{equation}\label{1.5}
\frac{dE_{tot}(\phi, U)}{dt}=-\|\sqrt{\sigma}\phi_t\|^2-\frac{\lambda D}{\varepsilon K}\|\nabla U\|^2\leq0.
\end{equation}
where $\|\cdot\|$ denotes the corresponding $L^2$ norm over $\Omega$. From \eqref{1.5}, one can see that the original energy decays in time.

\section{Numerical schemes}
In this section, we first introduce the auxiliary variable approach.
Denote the variable as follows
\begin{equation}\label{reforenergy}
\widetilde{E}(\phi, U)=\int_{\Omega}\left(\frac{1}{2}(m^2(\nabla \phi)-S_1) |\nabla \phi|^2+\frac{\lambda}{2 \varepsilon K} U^2+ (F(\phi)-\frac{S_2}{2\varepsilon^2}\phi^2)\right) d x,
\end{equation}
from which we obtain
\begin{equation}\label{2.1}
E_{tot}(\phi, U)=\widetilde{E}(\phi, U)+\int_{\Omega}\left(\frac{S_1}{2} |\nabla \phi|^2+ \frac{S_2}{2\varepsilon^2}\phi^2\right) d x.
\end{equation}
Here we add two positive constants $S_1, S_2$ to balance the anisotropic coefficient $m(\nabla\phi)$ and the nonlinear term $F(\phi)$ respectively, inspired by \cite{yang2019efficient}.
 For convenience, let $E(\phi, U)=E_{tot}(\phi, U)+1$, then we introduce the variable $q(t)$, defined by
\begin{equation}\label{2.2}
 q(t)=E(\phi, U),
 \end{equation}
 we rewrite \eqref{1.4} with the energy law \eqref{1.5} as the following expanded system
\begin{subequations}\label{2.3}
\begin{align}
& \sigma\phi_t- S_1\Delta\phi+\frac{S_2}{\varepsilon^2}\phi=-\frac{\delta \widetilde{E}}{\delta\phi}-\frac{\lambda}{\varepsilon} g(\phi)U, \\
& U_t-D\Delta U=Kg(\phi)\phi_t,\\
& \frac{d q(t)}{d t}=-\frac{q(t)}{E(\phi, U)} \mathcal{H}(\phi, U) .
\end{align}
\end{subequations}
where
\begin{equation}\label{2.4H}
\begin{split}
&\frac{\delta \widetilde{E}}{\delta\phi}=-\nabla \cdot(m^2(\nabla \phi) \nabla \phi+|\nabla \phi|^2 m(\nabla \phi) \boldsymbol{H}(\phi))+ S_1\Delta\phi+f(\phi)-\frac{S_2}{\varepsilon^2}\phi,\\
&\mathcal{H}(\phi, U)=\int_{\Omega}(\frac{1}{\sigma}| \frac{\delta E}{\delta \phi}+\frac{\lambda}{\varepsilon} g(\phi)U|^2+\frac{\lambda D}{\varepsilon K}|\nabla U|^2)dx \geq 0.
\end{split}
\end{equation}

\subsection{Second-order scheme on a nonuniform temporal mesh.}
We now turn to the construction of efficient time stepping schemes for \eqref{2.3}
that approximate the solution $\phi,U$ on a nonuniform temporal mesh. Let
$\tau_n := t_n - t_{n-1}$, $n = 1, 2, \ldots, M$, be the time step sizes,
denote the adjacent time step ratios as $\gamma_{n+1} := \tau_{n+1} / \tau_n$ , and
$\tau := \max\limits_{1\leq n\leq M} \tau_n$ be the maximum time step size.

Following the approach in \cite{hou2023implicit}, denote $F_2^{n+1}$ the nonuniform operator:
\begin{equation}\label{eq:bdf2}
F_2^{n+\delta}\bar{v}=\frac{1}{\tau_{n+1}}\Bigl[\frac{1+2\delta\gamma_{n+1}}{1+\gamma_{n+1}}\bar{v}^{n+1}
-(1+(2\delta-1)\gamma_{n+1})v^n
+\frac{(2\delta-1)\gamma_{n+1}^2}{1+\gamma_{n+1}}v^{n-1}\Bigr].
\end{equation}
In short, setting $\delta=1$ in \eqref{eq:bdf2} yields the classical second-order BDF scheme, while setting $\delta=\frac{1}{2}$ yields the Crank--Nicolson scheme.

Now we propose the  second-order generalized BDF-type scheme for the system \eqref{2.3} on a
nonuniform temporal mesh:

{\bf{Step 1}:} Compute $\bar{\phi}^{n+1},\bar{U}^{n+1},\overline{q}^{n+1}$ by:
\begin{subequations}\label{2.4}
\begin{align}
& \sigma F_2^{n+\delta}\bar{\phi}- S_1\Delta\bar{\phi}^{n+\delta}+\frac{S_2}{\varepsilon^2}\bar{\phi}^{n+\delta}
=-\frac{\delta \widetilde{E}}{\delta\phi}(\phi^{\ast,n+\delta})-\frac{\lambda}{\varepsilon }g(\phi^{\ast,n+\delta})U^{\ast,n+\delta}, \label{2.4a}\\
& F_2^{n+\delta} \bar{U}
-D\Delta\bar{U}^{n+\delta}=Kg(\phi^{\ast,n+\delta})
F_2^{n+\delta}\bar{\phi}, \label{2.4b}\\
& \frac{\overline{q}^{n+1}-q^n}{\tau_{n+1}}=-\frac{\overline{q}^{n+1}}{E(\bar{\phi}^{n+1},\bar{U}^{n+1})} \mathcal{H}(\bar{\phi}^{n+1},\bar{U}^{n+1}), \label{2.4c}
\end{align}
\end{subequations}
where we define $\bar{\phi}^{n+\delta}=\delta\bar{\phi}^{n+1}+(1-\delta)\phi^n$ and
$\phi^{*,n+\delta}=\phi^n+\delta\gamma_{n+1}(\phi^n-\phi^{n-1})$.

{\bf{Step 2}:} Compute the coefficients and the solution at step $n+1$,
\begin{subequations}\label{2.5}
\begin{align}
& \xi^{n+1}=\frac{\overline{q}^{n+1}}{E(\bar{\phi}^{n+1},\bar{U}^{n+1})}, \quad  \eta^{n+1}=1-(1-\xi^{n+1})^{3}, \label{2.4d}\\
& \phi^{n+1}=\eta^{n+1} \bar{\phi}^{n+1}, \quad
U^{n+1}=\eta^{n+1} \bar{U}^{n+1}, \label{2.4e}
\end{align}
\end{subequations}
then update the SAV $q^{n+1}$ via the following relaxation:
\begin{equation}\label{2.10}
q^{n+1}=\zeta^{n+1}\overline{q}^{n+1}+(1-\zeta^{n+1}) E(\phi^{n+1},U^{n+1}), \quad \zeta^{n+1} \in [0,1],
\end{equation}
where
\begin{equation}\label{2.6}
\frac{q^{n+1}-\overline{q}^{n+1}}{\tau_{n+1}}
\leq\frac{\overline{q}^{n+1}}{E(\bar{\phi}^{n+1},\bar{U}^{n+1})} \mathcal{H}(\bar{\phi}^{n+1},\bar{U}^{n+1}).
\end{equation}
We now establish the existence of $\zeta^{n+1}\in[0,1]$ in \eqref{2.10}.
Observe that this is equivalent to finding $\zeta^{n+1}\in[0,1]$ such that
\begin{equation}\label{opt3}
\left(\bar{q}^{n+1}-E(\phi^{n+1},U^{n+1})\right)\zeta^{n+1}
\leq R^{n+1},
\end{equation}
where
\[
R^{n+1}:=\bar{q}^{n+1}-E(\phi^{n+1},U^{n+1})
+\tau_{n+1}\,\frac{\bar{q}^{n+1}}{E(\bar{\phi}^{n+1},\bar{U}^{n+1})}
\mathcal{H}(\bar{\phi}^{n+1},\bar{U}^{n+1}).
\]

Since such a $\zeta_0^{n+1}$ is not unique, we select the value closest to
zero, which by \eqref{2.10} makes $q^{n+1}$ as close as possible to
$E(\phi^{n+1},U^{n+1})$. The choice is summarized as follows.

\begin{algorithm}\label{alg1}
(Determination of the parameter $\zeta^{n+1}\in[0,1]$)
\begin{enumerate}
\item If $\bar{q}^{n+1}\geq E(\phi^{n+1},U^{n+1})$,
      set $\zeta_0^{n+1}=0$.
\item If $\bar{q}^{n+1}<E(\phi^{n+1},U^{n+1})$ and $R^{n+1}\geq 0$, we set $\zeta_0^{n+1}=0$.
\item If $\bar{q}^{n+1}<E(\phi^{n+1},U^{n+1})$ and $R^{n+1}< 0$, we set
$
\zeta^{n+1}=1-\frac{\tau_{n+1} \bar{q}^{n+1} \mathcal{H}(\bar{\phi}^{n+1},\bar{U}^{n+1})}
       {E(\bar{\phi}^{n+1},\bar{U}^{n+1})(E(\phi^{n+1},U^{n+1})-\bar{q}^{n+1})}.
$
In this case, one readily verifies that $\zeta^{n+1}\in[0,1]$.
Indeed, since $q^n\ge 0$ and $E(\bar{\phi}^{n+1},\bar{U}^{n+1})>0$,
it follows from \eqref{2.4c} that
\begin{equation}\label{2.12}
\bar{q}^{n+1}=
\frac{q^n}{1+\tau_{n+1}\,
\dfrac{\mathcal{H}(\bar{\phi}^{n+1},\bar{U}^{n+1})}
{E(\bar{\phi}^{n+1},\bar{U}^{n+1})}}\ge 0,
\end{equation}
which, together with $\mathcal{H}(\bar{\phi}^{n+1},\bar{U}^{n+1})\ge 0$,
yields $0\le \zeta^{n+1}\le 1$.

\end{enumerate}
\end{algorithm}

\begin{remark}At the starting step $n=0$, the  second-order generalized BDF-type operator
is not available, we instead compute $\phi^{n}, U^{n}, q^n$ by the backward
Euler scheme  introduced in \cite{guo2024efficient}
\begin{subequations}\label{2.7}
\begin{align}
& \sigma\frac{\bar{\phi}^{n+1}-\phi^n}{\tau_{n+1}}- S_1\Delta\bar{\phi}^{n+1}+\frac{S_2}{\varepsilon^2}\bar{\phi}^{n+1}
=-\frac{\delta \widetilde{E}}{\delta\phi}(\phi^n)-\frac{\lambda}{\varepsilon }g(\phi^n)U^n, \label{2.7a}\\
& \frac{\bar{U}^{n+1}-U^n}{\tau_{n+1}}
-D\Delta\bar{U}^{n+1}=Kg(\phi^n)\frac{ \bar{\phi}^{n+1}-\phi^n}{\tau_{n+1}}, \label{2.7b}\\
& \frac{\overline{q}^{n+1}-q^n}{\tau_{n+1}}=-\frac{\overline{q}^{n+1}}{E(\bar{\phi}^{n+1},\bar{U}^{n+1})} \mathcal{H}(\bar{\phi}^{n+1},\bar{U}^{n+1}), \label{2.7c}
\end{align}
\end{subequations}
followed by the SAV relaxation in Step 2. Moreover, as established in \cite{guo2024efficient}, the scheme is unconditionally energy stable and obeys a discrete energy dissipation law, i.e.
\begin{equation}\label{2.14}
q^{n+1} \leq q^n.
\end{equation}
\end{remark}

\begin{theorem}\label{thm1}
 For $\frac{1}{2}\leq\delta\leq1$, the energy dissipation law of the scheme \eqref{2.4}-\eqref{2.6} is satisfied as follows
\begin{equation}\label{2.9}
q^{n+1}-q^n \leq 0 .
\end{equation}
In addition, there exists the constant $c>0$ such that
\begin{equation}\label{2.11}
\int_{\Omega}\Big(\frac{(1-\varepsilon_4)^2}{2} |\nabla \phi^{n}|^2+\frac{\lambda}{2 \varepsilon K} (U^{n})^2\Big)dx  \leq c, \quad \forall n .
\end{equation}
\end{theorem}

\begin{proof}
 Combining \eqref{2.4c} and \eqref{2.6}, we obtain the desired result \eqref{2.9}.

Denote $M:=q^0=E(\phi^{0}, U^{0})$, then we can infer from \eqref{2.9} and \eqref{2.12} that for every time $n$, $\overline{q}^{n+1} \leq M$ holds. For convenience, we divide the energy functional into two parts
\begin{equation}\label{energysplit}
E_{tot}(\phi, U)=\int_{\Omega}(\frac{1}{2}m^2(\nabla \phi) |\nabla \phi|^2+\frac{\lambda}{2 \varepsilon K} U^2)dx +\int_{\Omega}F(\phi) d x:=E_0(\phi,U)+E_1(\phi).
\end{equation}
 Applying \eqref{2.4d}, we derive that
\begin{equation}\label{2.13}
|\xi^{n+1}|=\frac{\overline{q}^{n+1}}{E(\bar{\phi}^{n+1},\bar{U}^{n+1})} \leq \frac{M}{E_{tot}(\bar{\phi}^{n+1}, \bar{U}^{n+1})+1} .
\end{equation}
According to the second equation in \eqref{2.4d}, we have $\eta^{n+1}=1-(1-\xi^{n+1})^{3}=\xi^{n+1}(3-3\xi^{n+1}+(\xi^{n+1})^2)$. Following \eqref{2.13}, we know there exists $c>0$ such that
$$
|\eta^{n+1}|=|\xi^{n+1}(3-3\xi^{n+1}+(\xi^{n+1})^2)| \leq \frac{c}{E_0(\bar{\phi}^{n+1}, \bar{U}^{n+1})+1}.
$$
 Using \eqref{2.4e}, i.e. $\phi^{n+1}=\eta^{n+1} \bar{\phi}^{n+1}$, $U^{n+1}=\eta^{n+1} \bar{U}^{n+1}$, we obtain
$$
\begin{aligned}
&\int_{\Omega}\Big(\frac{(1-\varepsilon_4)^2}{2} |\nabla \phi^{n+1}|^2+\frac{\lambda}{2 \varepsilon K} (U^{n+1})^2\Big)dx\\
 =
 &(\eta^{n+1})^2\int_{\Omega}\Big(\frac{1}{2}m^2(\nabla \bar{\phi}^{n+1}) |\nabla \bar{\phi}^{n+1}|^2+\frac{\lambda}{2 \varepsilon K} (\bar{U}^{n+1})^2\Big)dx \\
  \leq
&\Big(\frac{c}{E_0(\bar{\phi}^{n+1}, \bar{U}^{n+1})+1}\Big)^2 E_0(\bar{\phi}^{n+1}, \bar{U}^{n+1})\leq c,
\end{aligned}
$$
where
$m^2(\nabla \phi^{n+1}) \geq (1-\varepsilon_4)^2$. This completes the proof.
\end{proof}

\begin{remark}
In the analysis of variable-step BDF2 schemes, a mild ratio bound such as $\gamma_{n+1}\leq4.8645$ is frequently adopted, e.g.\ in \cite{hou2023implicit}. Such a bound stems from the analytical techniques employed for the variable step BDF2 operator rather than from the R-GSAV relaxation mechanism: as shown above, \eqref{2.9} follows purely from \eqref{2.4c} and \eqref{2.6} and therefore holds for arbitrary positive adjacent step ratios. Indeed, the numerical tests in Section 4 remain stable and accurate even for step ratios exceeding $4.8645$ (e.g.\ $\max\{\gamma_{n}\}=7.254$), so we retain this bound only in the adaptive strategy \eqref{eq:strategy} as a safeguard for temporal accuracy.
\end{remark}

\begin{remark}
From the point of view of energy approximation, we hope that $q^{n+1}=E(\phi^{n+1},U^{n+1})$. However, it is hard to prove energy stability when $q^{n+1}=$ $E(\phi^{n+1},U^{n+1})$ at each time step. Following Theorem \ref{thm1}, we can find that, if we take $\zeta^{n+1}=0$ in Cases 1  and 2, such that $q^{n+1}=E(\phi^{n+1},U^{n+1})$. Meanwhile, using \eqref{2.9} and \eqref{2.10}, we can deduce that
$$
E(\phi^{n+1},U^{n+1})=q^{n+1} \leq q^n \leq E(\phi^{n},U^{n}) \quad(\text {Cases } 1,2) .
$$
It is only in the third case that we cannot prove the dissipative nature of the original energy.
Here we need to point out in particular that in our experiments, the dissipation of the original energy is maintained at every time level, i.e. $\zeta^{n+1}=0$.
\end{remark}

\subsection{The spatial discretization}

Consider the two dimensional periodic domain $\Omega=(0,L)^2$ with a uniform spatial mesh of $N\times N$ grid points. Define the discrete Fourier basis functions
\[
\psi_{k,l}(x,y)=e^{\mathrm{i}(k x+l y)},\qquad k=\frac{2\pi k}{L},\quad l=\frac{2\pi l}{L},
\]
where the wave number indices range over
\[
(k,l)\in\mathcal{I}_N:=\Big\{-\frac{N}{2},-\frac{N}{2}+1,\dots,\frac{N}{2}-1\Big\}^2.
\]
The trial space for the pseudo-spectral approximation is
\[
X_N:=span\big\{\psi_{k,l}(x,y):\;(k,l)\in\mathcal{I}_N\big\}.
\]

For any periodic function $v\in L^2(\Omega)$, its truncated Fourier series and the corresponding discrete Fourier approximation are defined by
\[
v_N(x,y)=\sum_{(k,l)\in\mathcal{I}_N}\widehat{v}_{k,l}\,\psi_{k,l}(x,y),
\qquad
\widehat{v}_{k,l}=\frac{1}{N^2}\sum_{j_1,j_2=0}^{N-1}v(x_{j_1},y_{j_2})\,e^{-\mathrm{i}(k x_{j_1}+l y_{j_2})},
\]
where $\{(x_{j_1},y_{j_2})\}_{j_1,j_2=0}^{N-1}$ denotes the uniform collocation points. In practice, the discrete Fourier coefficients $\widehat{v}_{k,l}$ are computed by the Fast Fourier Transform (FFT), and the physical values are recovered by the Inverse FFT (IFFT).

For convenience, we denote the nonuniform BDF2 operator
\begin{equation}\label{oper}
F_2^{n+1}\widehat{v}:=a_{n+1}\bar{v}^{\,n+1}+b_{n+1}v^{\,n}+c_{n+1}v^{\,n-1},
\end{equation}
where
\[
a_{n+1}=\frac{1+2\delta\gamma_{n+1}}{\tau_{n+1}(1+\gamma_{n+1})},\qquad
b_{n+1}=-\frac{1+(2\delta-1)\gamma_{n+1}}{\tau_{n+1}},\qquad
c_{n+1}=\frac{(2\delta-1)\gamma_{n+1}^2}{\tau_{n+1}(1+\gamma_{n+1})}.
\]

{\bf{Step 1: Linear Prediction $(\bar{\phi}_{N}^{n+1},\bar{U}_{N}^{n+1},\bar{q}_{N}^{n+1})$}}

The linear equation for $\bar{\phi}^{n+1}\in X_N$ reads, in Fourier space,
\begin{equation}\label{eq:phi_bar}
\Big(\sigma a_{n+1}+S_1\delta(k^2+l^2)+\frac{S_2\delta}{\varepsilon^2}\Big)\widehat{\bar{\phi}}^{\,n+1}_{k,l}
=\widehat{\mathcal{G}}^{\,\phi}_{k,l}
-\sigma\big(b_{n+1}\widehat{\phi}^{\,n}_{k,l}+c_{n+1}\widehat{\phi}^{\,n-1}_{k,l}\big),
\qquad\forall\,(k,l)\in\mathcal{I}_N,
\end{equation}
where the right hand side is assembled pseudo spectrally:
\begin{equation}\label{eq:rhs_phi}
\widehat{\mathcal{G}}^{\,\phi}
=-\frac{\delta \widetilde{E}}{\delta\phi}(\phi_{N}^{\ast,n+\delta})-\frac{\lambda}{\varepsilon }g(\phi_{N}^{\ast,n+\delta})U_{N}^{\ast,n+\delta}
+S_1(1-\delta)\Delta\phi_{N}^n-\frac{S_2}{\varepsilon^2}(1-\delta)\phi_{N}^n.
\end{equation}
 Since the left hand side of \eqref{eq:phi_bar} is diagonal in $(k,l)$, it is inverted pointwise:
\begin{equation}\label{solve_phi}
\widehat{\bar{\phi}}^{\,n+1}_{k,l}
=\frac{\widehat{\mathcal{G}}^{\,\phi}_{k,l}
-\sigma\big(b_{n+1}\widehat{\phi}^{\,n}_{k,l}+c_{n+1}\widehat{\phi}^{\,n-1}_{k,l}\big)}
{\sigma a_{n+1}+S_1\delta|k|^2+S_2\delta/\varepsilon^2}.
\end{equation}
Once $\bar{\phi}^{n+1}$ is obtained, compute the discrete time derivative in physical space using the definition of \eqref{oper}.
The right hand side for the temperature equation is assembled pseudo spectrally:
\[
\mathcal{G}^U=Kg(\phi_{N}^{\ast,n+\delta})
F_2^{n+\delta}\bar{\phi_{N}}
+D(1-\delta)\Delta U_{N}^n.
\]
Then $\bar{U}^{n+1}\in X_N$ is determined in Fourier space by
\begin{equation}\label{eq:U_bar}
\big(a_{n+1}+D\delta(k^2+l^2)\big)\widehat{\bar{U}}^{\,n+1}_{k,l}
=\widehat{\mathcal{G}}^{\,U}_{k,l}
-\big(b_{n+1}\widehat{U}^{\,n}_{k,l}+c_{n+1}\widehat{U}^{\,n-1}_{k,l}\big),
\qquad\forall\,(k,l)\in\mathcal{I}_N,
\end{equation}
which is again inverted pointwise:
\begin{equation}\label{solve_U}
\widehat{\bar{U}}^{\,n+1}_{k,l}
=\frac{\widehat{\mathcal{G}}^{\,U}_{k,l}
-\big(b_{n+1}\widehat{U}^{\,n}_{k,l}+c_{n+1}\widehat{U}^{\,n-1}_{k,l}\big)}
{a_{n+1}+D\delta(k^2+l^2)}.
\end{equation}
Define the discrete original energy of the intermediate solution:
\begin{equation}\label{eq:E_bar}
\bar{E}=\frac{L^2}{N^2}\sum_{j_1,j_2=0}^{N-1}\Bigg[
\frac{1}{2}m^2(\nabla\bar{\phi}_{N}^{n+1})\big|\nabla\bar{\phi}_{N}^{n+1}\big|^2
+\frac{\lambda}{2\varepsilon K}\big(\bar{U}_{N}^{n+1}\big)^2
+\frac{1}{4\varepsilon^2}\Big(\big(\bar{\phi}_{N}^{n+1}\big)^2-1\Big)^2
\Bigg]_{(x_{j_1},y_{j_2})}+\,1,
\end{equation}
and the discrete dissipation rate:
\begin{equation}\label{eq:H_bar}
\bar{\mathcal{H}}
=\frac{L^2}{N^2}\sum_{j_1,j_2=0}^{N-1}\Bigg[
\sigma\big|F_2^{n+1}\bar{\phi}_{N}\big|^2
+\frac{\lambda D}{\varepsilon K}\Big(\big|\nabla\bar{U}_{N}^{n+1}\big|^2\Big)
\Bigg]_{(x_{j_1},y_{j_2})}.
\end{equation}
Then the intermediate auxiliary variable is updated explicitly:
\begin{equation}\label{eq:q_bar}
\bar{q}_{N}^{n+1}=\frac{q_{N}^n}{1+\tau_{n+1}\,\dfrac{\bar{\mathcal{H}}}{\bar{E}}}.
\end{equation}
{\bf{Step 2: Relaxation Correction $(\phi_{N}^{n+1},U_{N}^{n+1},q_{N}^{n+1})$}}
\begin{equation}\label{eq:xi_eta}
\xi^{n+1}=\frac{\bar{q}_{N}^{n+1}}{\bar{E}},\qquad
\eta^{n+1}=1-(1-\xi_{N}^{n+1})^3.
\end{equation}
The corrected solution is obtained by uniform scaling in physical space:
\begin{equation}\label{eq:correction_fields}
\widehat{\phi}^{\,n+1}_{k,l}=\eta^{n+1}\,\widehat{\bar{\phi}}^{\,n+1}_{k,l},\qquad
\widehat{U}^{\,n+1}_{k,l}=\eta^{n+1}\,\widehat{\bar{U}}^{\,n+1}_{k,l},
\qquad\forall\,(k,l)\in\mathcal{I}_N.
\end{equation}
Define
\[
R^{n+1}=\bar{q}_{N}^{n+1}-E(\phi_{N}^{n+1},U_{N}^{n+1})+\tau_{n+1}\,\frac{\bar{q}_{N}^{n+1}\,\bar{\mathcal{H}}}{\bar{E}}.
\]
The relaxation parameter $\zeta^{n+1}\in[0,1]$ is chosen according to Algorithm \ref{alg1}:
\begin{equation}\label{eq:zeta}
\zeta^{n+1}=
\begin{cases}
0, & \text{if }\bar{q}_{N}^{n+1}\ge E(\phi_{N}^{n+1},U_{N}^{n+1}),\\[4pt]
0, & \text{if }\bar{q}_{N}^{n+1}<E(\phi_{N}^{n+1},U_{N}^{n+1})  \text{ and }R^{n+1}\ge 0,\\[8pt]
\displaystyle 1-\frac{\tau_{n+1}\,\bar{q}_{N}^{n+1}\,\bar{\mathcal{H}}}{\bar{E}\,\big(E^{n+1}-\bar{q}^{n+1}\big)},
& \text{if }\bar{q}_{N}^{n+1}<E(\phi_{N}^{n+1},U_{N}^{n+1})\text{ and }R^{n+1}<0.
\end{cases}
\end{equation}
Finally, update the auxiliary variable:
\begin{equation}\label{eq:q_update}
q_{N}^{n+1}=\zeta^{n+1}\,\bar{q}_{N}^{n+1}+(1-\zeta^{n+1})\,E(\phi_{N}^{n+1},U_{N}^{n+1}).
\end{equation}

\begin{theorem}\label{thm:fd}
 Without any restriction on the step ratios, the scheme satisfies the energy dissipation law
\[
q_{N}^{n+1}- q_{N}^{n}\le 0,\qquad\forall\,n\ge 1.
\]
Moreover, in Cases 1 and 2 of \eqref{eq:zeta} where $\zeta^{n+1}=0$, we have $q_{N}^{n+1}=E(\phi_{N}^{n+1},U_{N}^{n+1})$, and consequently the original energy is dissipated:
\[
E(\phi_{N}^{n+1},U_{N}^{n+1})=q_{N}^{n+1}\le q_{N}^{n}\le E(\phi_{N}^{n},U_{N}^{n}).
\]
\end{theorem}

\section{Numerical experiments}
In this section, some numerical tests are implemented to verify the feasibility of the proposed scheme.
Let us start by verifying the convergence of the proposed schemes \eqref{2.4} by constructing appropriate exact solutions. Then, we simulate the dendrite growth with fourfold anisotropy and sixfold anisotropy under different parameters $K$. Finally, we study the 3D growth of the dendrite with fourfold anisotropy. In all examples, the domain is chosen as $[0,2\pi]^d$, where $d=2,3$. For simplicity, we adopt periodic boundary conditions and the Fourier spectral method for space discretization.

\subsection{The test of the convergence order}
In this subsection, we test the convergence of the proposed schemes. We use the exact solution hereafter, and add the source term to the original equation. Here, we use the $256\times256$ Fourier modes for the spatial discretization so that the error of the  spatial discretization is negligible compared to the temporal discretization.

To examine the temporal convergence order of the proposed scheme \eqref{2.4} under nonuniform temporal partitions, we construct a family of random time grids by perturbing the uniform mesh $\{t_n=n\tau\}_{n=0}^{M}$ with a $40\%$ random disturbance, i.e.
\[
\widehat{t}_{n}=t_{n}+0.4\,\tau\cdot\mathrm{rand}(\cdot),\qquad n=0,1,\ldots,M,
\]
where $\mathrm{rand}(\cdot)$ is a uniformly distributed random number in $[-1,1]$.
The actual time-step sizes are then obtained by $\tau_{n+1}=\widehat{t}_{n+1}-\widehat{t}_{n}$ after fixing the end points $\widehat{t}_{0}=0$ and $\widehat{t}_{M}=T$.
Such a construction produces highly irregular temporal meshes with adjacent step ratios $\gamma_{n+1}=\tau_{n+1}/\tau_{n}$ varying dramatically across the interval $[0,T]$.

 \begin{example}\label{exp.1}
\rm{We solve the system \eqref{1.4} with the exact solution}
\[
\phi(x,y,t)=\sin(x)\cos(y)\cos(t),\qquad U(x,y,t)=\cos(x)\sin(y)\cos(t),
\]
\rm{on the domain $\Omega=[0,2\pi]^{2}$ up to $T=1$.}
\rm{The physical parameters are taken as}
 $\sigma=4\times10^{3}$, $\varepsilon=1$, $\lambda=1$, $D=2.25\times10^{-2}$, $K=0.01$, $\varepsilon_{4}=0.05$ and $\beta=4$.
\end{example}
Tables 1--3 list the $L^{2}$-errors of the phase field $\phi$ and the temperature $U$ at $T=1$ for three different choices of the  parameter $\delta$, namely $\delta=1/2$, $3/4$, and $1$, on the above random nonuniform meshes.
The observed adjacent step ratios $\max\{\gamma_{n}\}$ are also reported.
It can be seen from the table that the proposed scheme achieves approximately second-order accuracy in time for both $\phi$ and $U$ for all tested $\delta$.
 It is worth noting that some step ratios exceed the commonly used analysis bound $4.8645$ (e.g. $\max\{\gamma_{n}\}=7.254$ for $M=160$), yet the scheme remains stable and retains the expected convergence rate, which is consistent with the unconditional energy dissipation laws established in Theorem \ref{thm1} and Theorem \ref{thm:fd}.
Moreover, the least squares fitted convergence orders are all above $2$ for both $\phi$ and $U$ for all three choices of $\delta$, confirming the robustness of the adaptive BDF2-RGSAV solver on irregular temporal meshes.

\begin{table}[htbp]\label{tab1}
\centering
\caption{Convergence rates for the phase function $\phi^n$ and the temperature $U^n$ in \eqref{2.4} on nonuniform temporal meshes with $\delta=\frac{1}{2}$.}
\label{tab:convergence1}
\begin{tabular}{ccc|cc|cc}
  \hline
$M$ & $\tau$ & $\max\{\gamma_n\}$ & $\phi_{error}$ & $\phi_{order}$ & $U_{error}$ & $U_{order}$ \\
  \hline
    20 &   8.2391e-02 &      2.649 &     6.9570e-04 & --    &     6.9695e-04 & --    \\
    40 &   3.9791e-02 &      3.104 &     1.3869e-04 &  2.216 &     1.3839e-04 &  2.221 \\
    80 &   2.0986e-02 &      4.150 &     3.0873e-05 &  2.348 &     3.0938e-05 &  2.342 \\
   160 &   1.0672e-02 &      7.254 &     7.9409e-06 &  2.008 &     7.9641e-06 &  2.007 \\
   320 &   5.5151e-03 &      5.836 &     1.7555e-06 &  2.286 &     1.7627e-06 &  2.285 \\
   640 &   2.7229e-03 &      6.690 &     4.4966e-07 &  1.930 &     4.2182e-07 &  2.026 \\
  \hline
\end{tabular}
\end{table}

\begin{table}[htbp]\label{tab2}
\centering
\caption{Convergence rates for the phase function $\phi^n$ and the temperature $U^n$ in \eqref{2.4} on nonuniform temporal meshes with $\delta=\frac{3}{4}$.}
\label{tab:convergence2}
\begin{tabular}{ccc|cc|cc}
  \hline
$M$ & $\tau$ & $\max\{\gamma_n\}$ & $\phi_{error}$ & $\phi_{order}$ & $U_{error}$ & $U_{order}$ \\
  \hline
    20 &   8.2391e-02 &      2.649 &     7.7251e-04 & --    &     7.7338e-04 & --    \\
    40 &   3.9791e-02 &      3.104 &     1.5849e-04 &  2.176 &     1.5816e-04 &  2.181 \\
    80 &   2.0986e-02 &      4.150 &     3.7381e-05 &  2.258 &     3.7407e-05 &  2.253 \\
   160 &   1.0672e-02 &      7.254 &     9.5723e-06 &  2.015 &     9.5849e-06 &  2.014 \\
   320 &   5.5151e-03 &      5.836 &     2.1443e-06 &  2.266 &     2.1505e-06 &  2.264 \\
   640 &   2.7229e-03 &      6.690 &     5.5034e-07 &  1.927 &     5.2786e-07 &  1.990 \\
  \hline
\end{tabular}
\end{table}

\begin{table}[htbp]\label{tab3}
\centering
\caption{Convergence rates for the phase function $\phi^n$ and the temperature $U^n$ in \eqref{2.4} on nonuniform temporal meshes with $\delta=1$.}
\label{tab:convergence3}
\begin{tabular}{ccc|cc|cc}
  \hline
$M$ & $\tau$ & $\max\{\gamma_n\}$ & $\phi_{error}$ & $\phi_{order}$ & $U_{error}$ & $U_{order}$ \\
  \hline
   20 &   8.2391e-02 &      2.649 &     8.3934e-04 & --    &     8.3955e-04 & --    \\
    40 &   3.9791e-02 &      3.104 &     1.8023e-04 &  2.114 &     1.7969e-04 &  2.118 \\
    80 &   2.0986e-02 &      4.150 &     4.3721e-05 &  2.214 &     4.3684e-05 &  2.211 \\
   160 &   1.0672e-02 &      7.254 &     1.1202e-05 &  2.014 &     1.1198e-05 &  2.013 \\
   320 &   5.5151e-03 &      5.836 &     2.5466e-06 &  2.244 &     2.5496e-06 &  2.242 \\
   640 &   2.7229e-03 &      6.690 &     6.5315e-07 &  1.928 &     6.3338e-07 &  1.973 \\
  \hline
\end{tabular}
\end{table}

\FloatBarrier
\subsection{Time adaptive strategy}

In this part, the efficiency of the fully discrete scheme \eqref{2.4} with a time adaptive strategy is investigated through a dendritic
solidification process. As we know, the solidification process usually goes
through different stages: the circular seed grows rapidly at the beginning due
to the large undercooling, accompanied by rapid latent heat release and fast
energy decay, and then the evolution becomes rather slow afterwards until the
system reaches a steady state. Thus it would be very efficient to apply some
time adaptive strategy in the simulation based on our proposed nonuniform BDF2
scheme \eqref{2.4}.

 \begin{example}\label{exp.2}
\rm{  We consider the initial condition as a circle shape below}
\begin{equation}
\begin{cases}
\phi(x,y,t=0)=\tanh\Big(\dfrac{1.5-\sqrt{(x-\pi)^2+(y-\pi)^2}}{0.1}\Big),\\[2mm]
U(x,y,t=0)=-0.55.
\end{cases}
\end{equation}
\rm{The computational domain is set to be $\Omega=(0,2\pi)^2$, and the parameters
are taken as}
\[
\sigma=1\mathrm{e}3,\quad \varepsilon=0.1,\quad \lambda=1,\quad
D=2.25\mathrm{e}{-1},\quad K=1,\quad \varepsilon_4=0.05,\quad \beta=4.
\]
\end{example}
The time adaptive strategy proposed in \cite{hou2023implicit} is adopted in the
simulation with our fully discrete scheme \eqref{2.4}:
\begin{equation}\label{eq:strategy}
\begin{aligned}
&\tau_{n+1}=\min\Big\{\max\Big\{\tau_{\min},
\frac{\tau_{\max}}{\sqrt{1+\mu|E'(t)|^2}}\Big\},\ 4.8645\cdot\tau_n\Big\},
\end{aligned}
\end{equation}
The mild ratio bound $4.8645$ is not required for the energy dissipation laws in Theorem \ref{thm1} and Theorem \ref{thm:fd}; it is retained in the strategy \eqref{eq:strategy} only as a safeguard for temporal accuracy.
where $\tau_{\min}$, $\tau_{\max}$ are predetermined minimum and maximum time
step sizes, and  $\mu$ is a positive constant. This type of time adaptive
strategy is based on the energy variation. It will automatically select small
time steps when the energy decays rapidly and choose larger time steps
otherwise.  In this simulation, we set $\tau_{\min}=0.01$, $\mu=10^3$ and
the maximum step sizes $\tau_{\max}=1$.
Furthermore, $256\times256$ Fourier modes are used for the spatial
discretization.

In Figure 1, we plot the numerical results up to $T=100$ with three different
types of temporal meshes, i.e., the fixed large time step mesh with $\tau=1$, the uniform small time step mesh with $\tau=0.01$ as the reference solution, and the adaptive temporal mesh with $\tau_{\max}=1$.
It is observed in Figure 1 (a) that the computed energy is dissipative in time for all the tested temporal meshes. Note that $\zeta^{n+1}=0$ holds at every time level in this simulation, hence the modified energy is exactly the original energy, i.e., $q^{n+1}=E(\phi^{n+1},U^{n+1})$. However, after the energy goes through a large variation stage, the fixed large time step $\tau=1$ yields an inaccurate energy evolution, while the adaptive time strategies with $\tau_{\max}=1$  give energy curves that
agree very well with the one using the fixed small time step $\tau=0.01$.
Figure 1 (b) shows the evolution of $|\xi^n-1|$, which remains small
(about $10^{-7}$--$10^{-4}$) during the whole simulation, indicating that the
auxiliary variable $\xi^n$ of the GSAV approach is always an accurate
approximation of $1$ even on the nonuniform temporal meshes. The corresponding
time step sizes are plotted in Figure 1 (c): small time steps near $\tau_{\min}$
are automatically selected at the early rapid solidification stage with fast
energy decay, while the step size grows towards $\tau_{\max}$ once the
evolution slows down, which shows the efficiency of our proposed scheme \eqref{2.4}
with the adaptive strategy \eqref{eq:strategy} during this long time simulation.
Figure \ref{fig:adaptiveDynamics} presents the snapshots of the phase field $\phi$ at $t=0,20,60$ and $100$ computed by the scheme \eqref{2.4} with the adaptive strategy \eqref{eq:strategy} ($\tau_{\max}=1$).

 Although the scalar coefficients $a_{n+1}$, $b_{n+1}$ and $c_{n+1}$ in \eqref{oper} vary with the time step $\tau_{n+1}$ and the adjacent step ratio $\gamma_{n+1}$, the spatial operators of the two linear systems retain the same structure and remain diagonal in Fourier space, so that each time step only requires pointwise FFT-based solves whose cost is independent of the temporal mesh. Consequently, the computational cost per time step is nearly identical for all the tested meshes, and hence the ratio of the step counts is approximately the ratio of the CPU times. As listed in Table \ref{tab:steps42}, the adaptive simulation reaches almost the same accuracy as the reference solution with $272$ steps and CPU cost of about $3.2\%$ of that of the reference simulation, while the fixed large step $\tau=1$, though the cheapest one, produces a relative energy error two orders of magnitude larger than that of the adaptive one, which clearly demonstrates the efficiency of the proposed adaptive scheme.

\begin{figure}[htbp]
\subfigure[\scriptsize (a) Energy evolution with different time step.]{
\begin{minipage}{0.45 \linewidth}
\centering
\includegraphics[width=\textwidth]{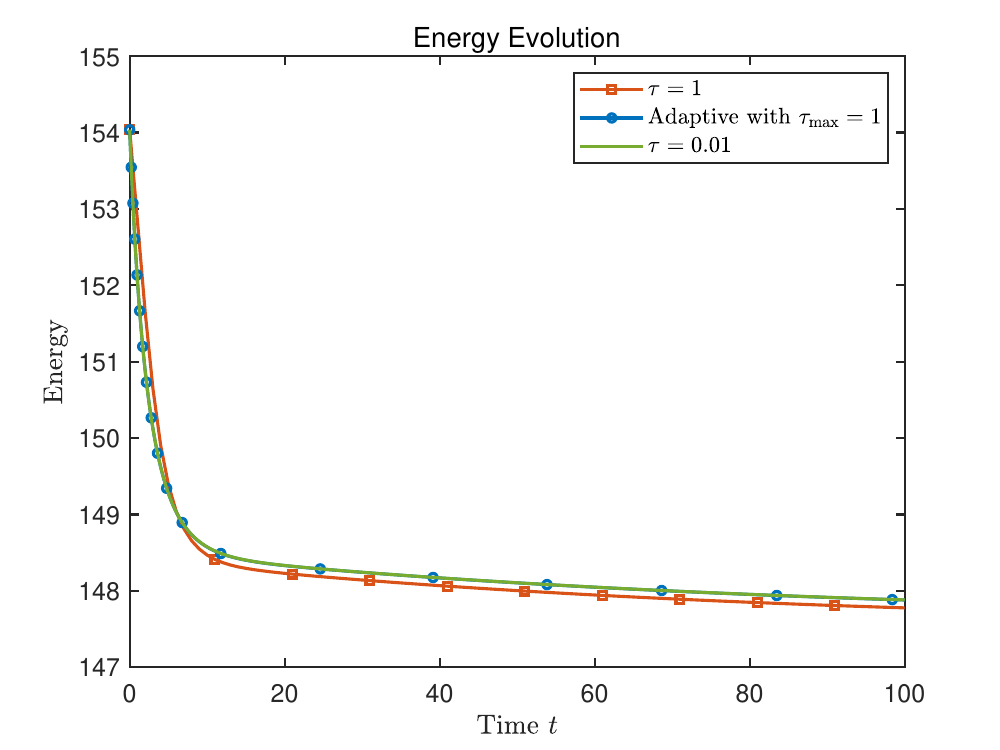}
\end{minipage}
}
\subfigure[\scriptsize (b) $|\xi^n-1|$]{
\begin{minipage}{0.45\linewidth}
\centering
\includegraphics[width=\textwidth]{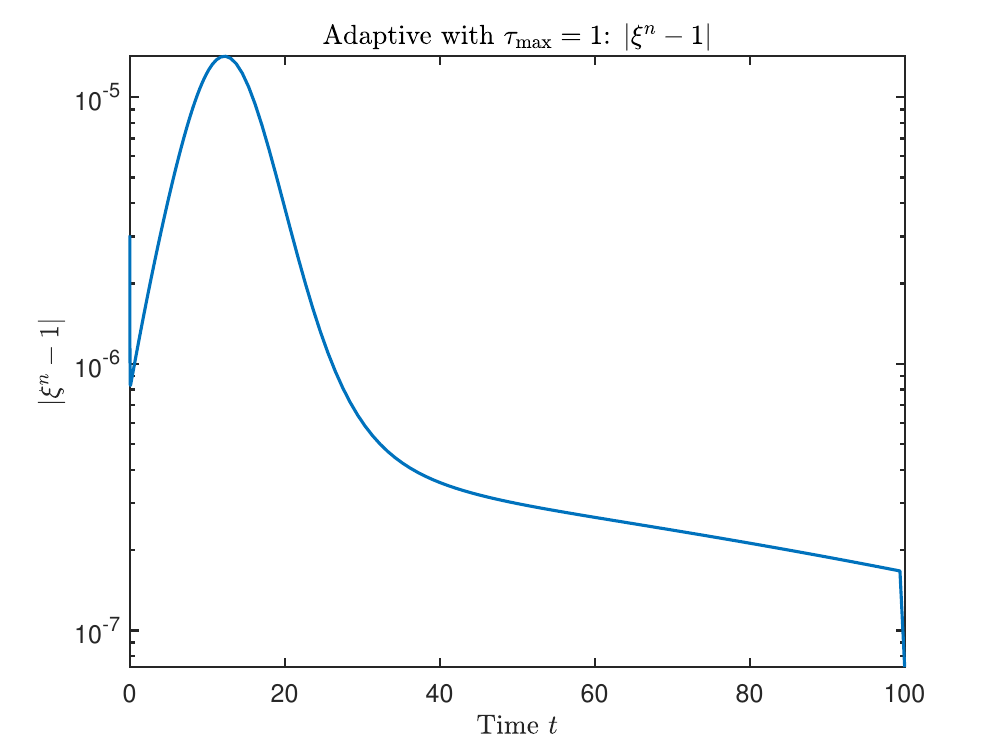}
\end{minipage}
}

\subfigure[\scriptsize (c) adaptive time step size with $\tau_{\max}=1$]{
\begin{minipage}{0.45 \linewidth}
\centering
\includegraphics[width=\textwidth]{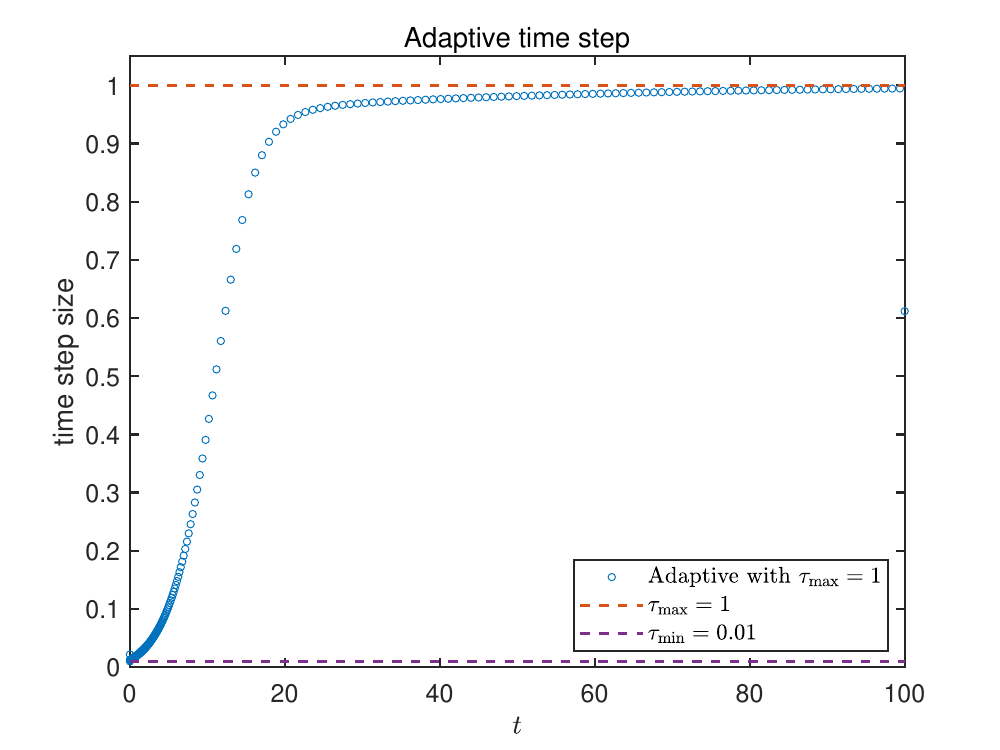}
\end{minipage}
}

\caption{(Example \ref{exp.2}) (a) The evolution of energy with different time steps. (b) The evolution of the deviation $|\xi^{n}-1|$ for the adaptive time stepping simulation. (c) The evolution of the adaptive time step size with $\tau_{\max}=1$.}
\label{energy_stable}
\end{figure}

\begin{figure}[htbp]

	\vspace{-0.35cm}
	\subfigtopskip=2pt
	\subfigbottomskip=2pt
	\subfigcapskip=-5pt %

	\centering	

    \subfigure{
		\begin{minipage}[t]{0.185\linewidth}
			\centering
			\includegraphics[width=1\linewidth]{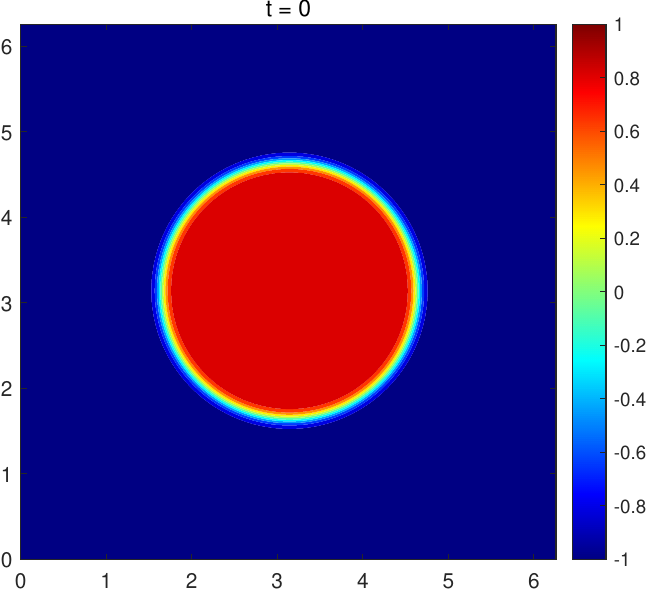}
		\end{minipage}
		\begin{minipage}[t]{0.185\linewidth}
			\centering
			\includegraphics[width=1\linewidth]{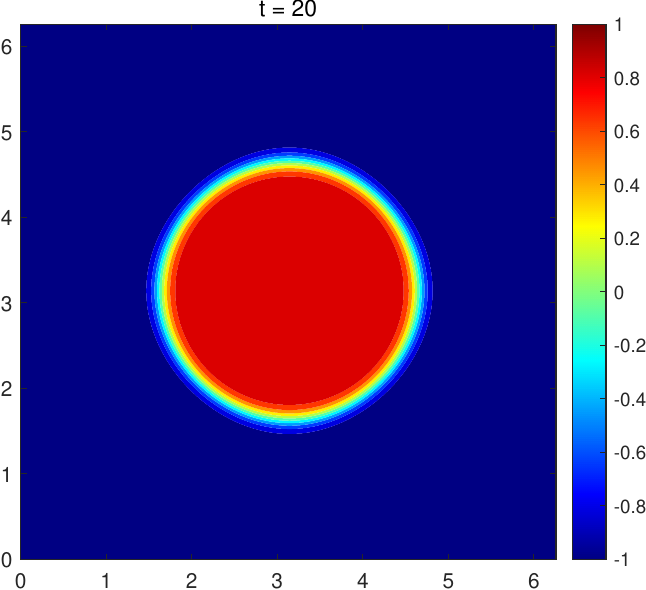}
		\end{minipage}
		\begin{minipage}[t]{0.185\linewidth}
			\centering
			\includegraphics[width=1\linewidth]{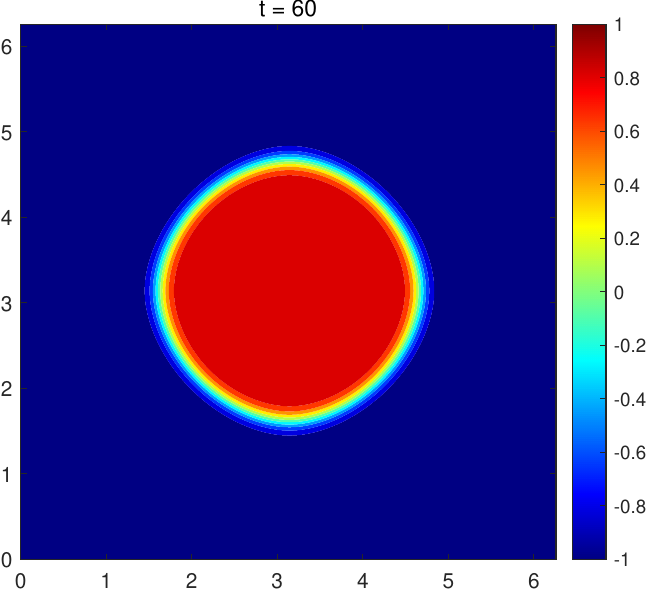}
		\end{minipage}
		\begin{minipage}[t]{0.185\linewidth}
			\centering
			\includegraphics[width=1\linewidth]{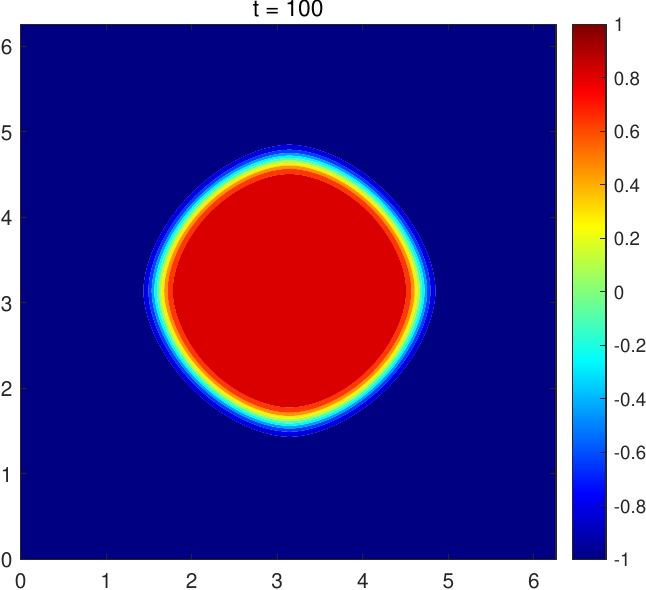}
		\end{minipage}
	}
	
	\caption{(Example \ref{exp.2}) Dynamics driven by the dendritic crystal growth model with adaptive time steps.}
	\label{fig:adaptiveDynamics}
\end{figure}

\begin{table}[htbp]
\centering
\caption{(Example \ref{exp.2}) Comparison of the numbers of time steps, the CPU times and the relative errors of the computed energy at $T=100$ for different temporal meshes, where the reference solution is computed with the uniform small time step $\tau=0.01$.}
\label{tab:steps42}
\begin{tabular}{l|ccc}
  \hline
 temporal mesh & number of time steps & CPU time (s) & relative energy error\\
  \hline
 uniform $\tau=1$ & 100 & 6.75 & $7.05\times10^{-4}$ \\
 uniform $\tau=0.01$ (reference) & 10000 & 645.63 & -- \\
 adaptive strategy \eqref{eq:strategy}  & 272 & 20.69 & $6.76\times10^{-6}$ \\
  \hline
\end{tabular}
\end{table}
\subsection{Fourfold dendrite crystal growth in 2D}
In this subsection, we investigate the fourfold dendrite crystal growth and the effect of the latent heat parameter $K$ on the crystal shape. Different from the uniform time stepping simulations in \cite{yang2020efficient,guo2024efficient}, the temporal discretization here is carried out by the proposed variable step R-GSAV BDF2 scheme \eqref{2.4} (with $\delta=1$) combined with the time adaptive strategy \eqref{eq:strategy}, in which the step sizes are automatically adjusted according to the energy variation  with $\tau_{\min}=10^{-2}$, $\tau_{\max}=1$, $\mu=10^{-2}$ and the initial step size $\tau_0=10^{-2}$. The Fourier spectral method with $512\times512$ Fourier modes is employed for the spatial discretization, and the physical parameters are taken the same as those in \cite{yang2020efficient}.
 \begin{example}\label{exp.4}
\rm{The initial phase field is taken as }
\begin{equation}\label{e4}
   \phi(x,y,t=0)=\tanh\Big(\frac{0.02-\sqrt{(x-\pi)^2+(y-\pi)^2}}{0.072}\Big),
\end{equation}
\rm{and the initial temperature field is set as follows}
\begin{equation*}
U(x,y,t=0)=\left\{
             \begin{array}{ll}
             0,  &\phi>0,\\
             -0.55, &otherwise.
                           \end{array}
           \right.
\end{equation*}
We choose the parameters as
$$
\begin{aligned}
\sigma=4.4e3, \varepsilon=1.12e-2, \lambda=380, D=2.25e-4, \varepsilon_4=0.05,\beta=4, S_1=4, S_2=4.
\end{aligned}
$$
\end{example}

\begin{figure}[htbp]

	\vspace{-0.35cm}
	\subfigtopskip=2pt
	\subfigbottomskip=2pt
	\subfigcapskip=-5pt

	\centering	

    \subfigure[\scriptsize (a) Snapshots of the phase field $\phi$.]{

		\begin{minipage}[t]{0.185\linewidth}
			\centering
			\includegraphics[width=1\linewidth]{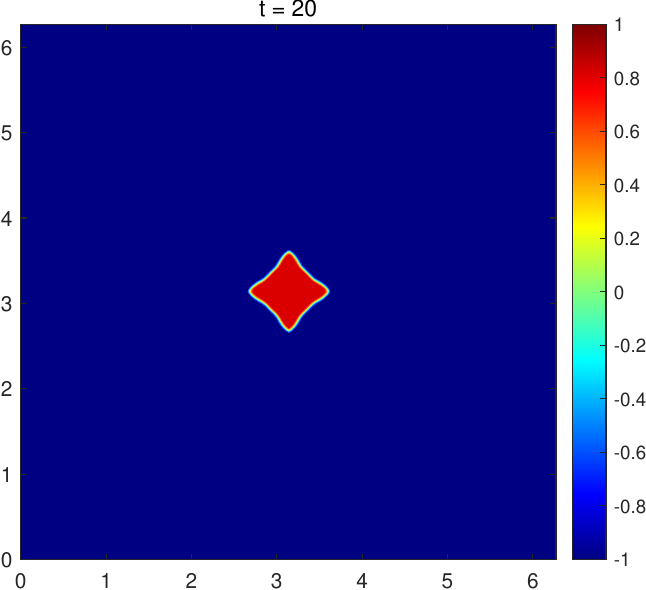}
		\end{minipage}
		\begin{minipage}[t]{0.185\linewidth}
			\centering
			\includegraphics[width=1\linewidth]{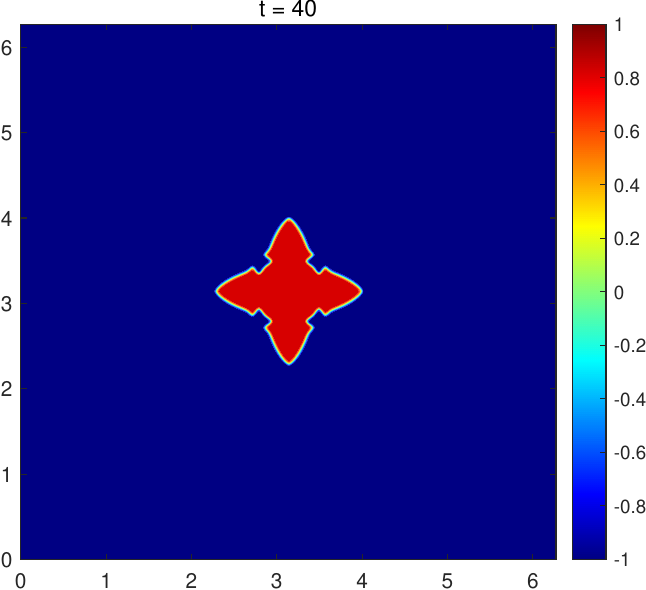}
		\end{minipage}
		\begin{minipage}[t]{0.185\linewidth}
			\centering
			\includegraphics[width=1\linewidth]{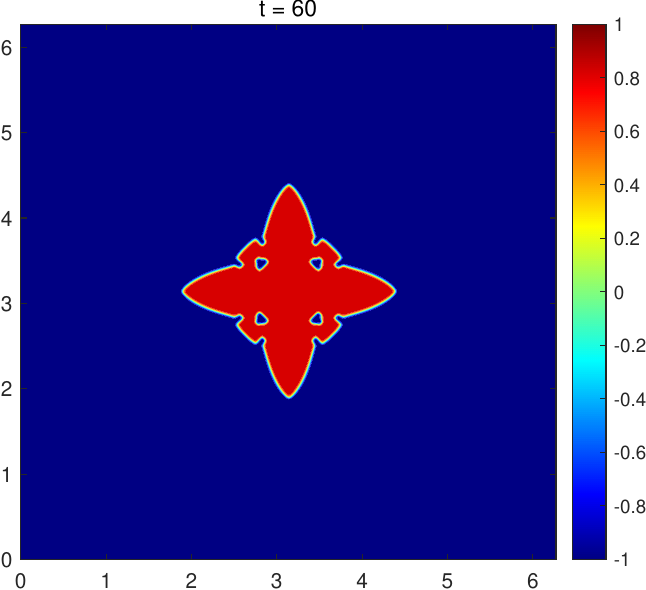}
		\end{minipage}
		\begin{minipage}[t]{0.185\linewidth}
			\centering
			\includegraphics[width=1\linewidth]{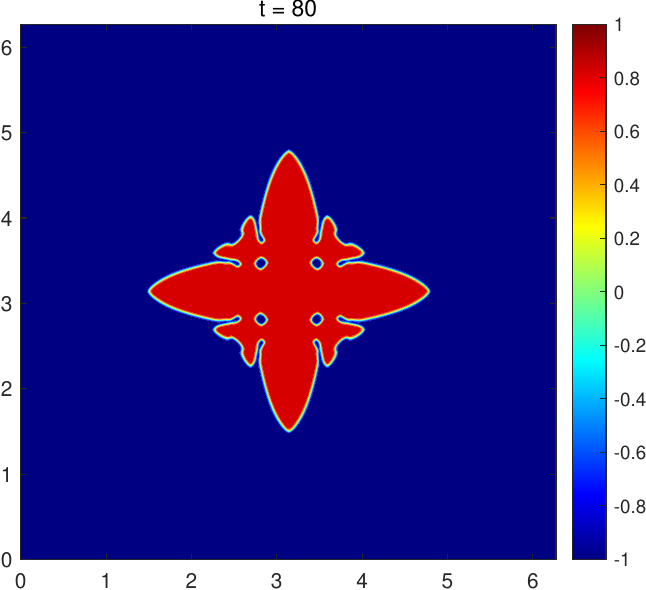}
		\end{minipage}
		\begin{minipage}[t]{0.185\linewidth}
			\centering
			\includegraphics[width=1\linewidth]{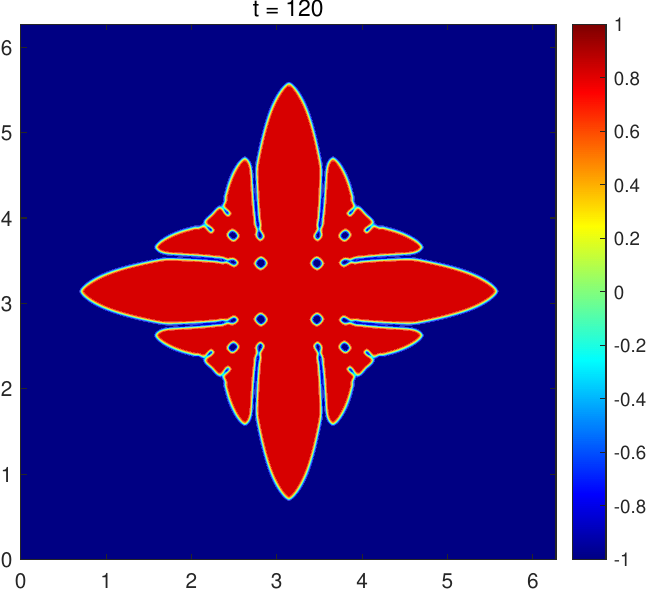}
		\end{minipage}
	}

    \subfigure[\scriptsize (b) Snapshots of the temperature field $U$.]{

		\begin{minipage}[t]{0.185\linewidth}
			\centering
			\includegraphics[width=1\linewidth]{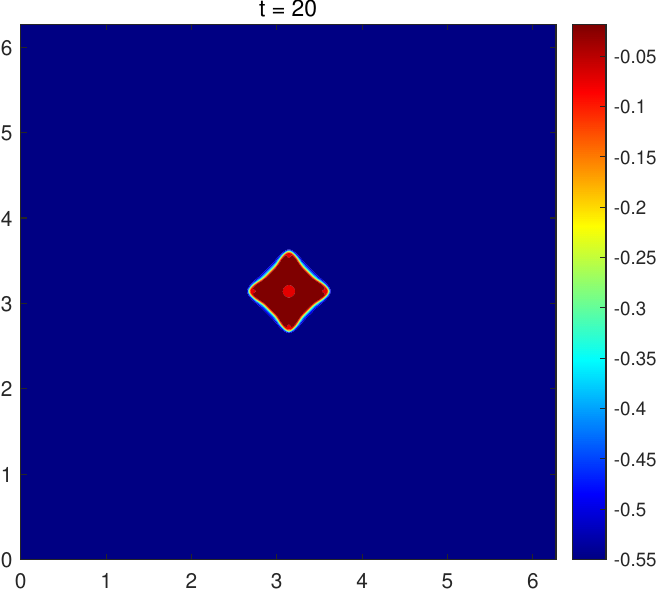}
		\end{minipage}
		\begin{minipage}[t]{0.185\linewidth}
			\centering
			\includegraphics[width=1\linewidth]{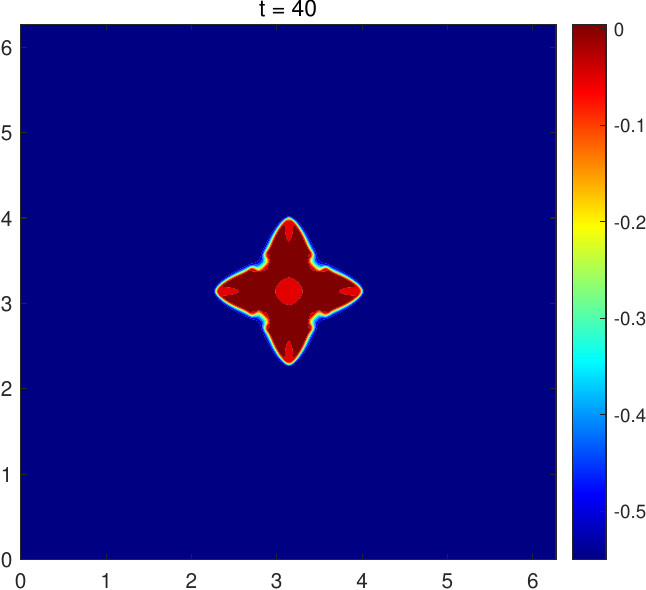}
		\end{minipage}
		\begin{minipage}[t]{0.185\linewidth}
			\centering
			\includegraphics[width=1\linewidth]{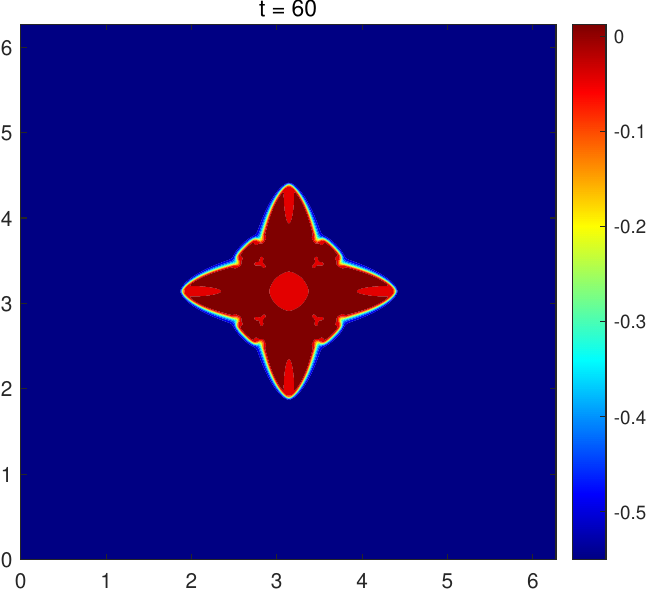}
		\end{minipage}
		\begin{minipage}[t]{0.185\linewidth}
			\centering
			\includegraphics[width=1\linewidth]{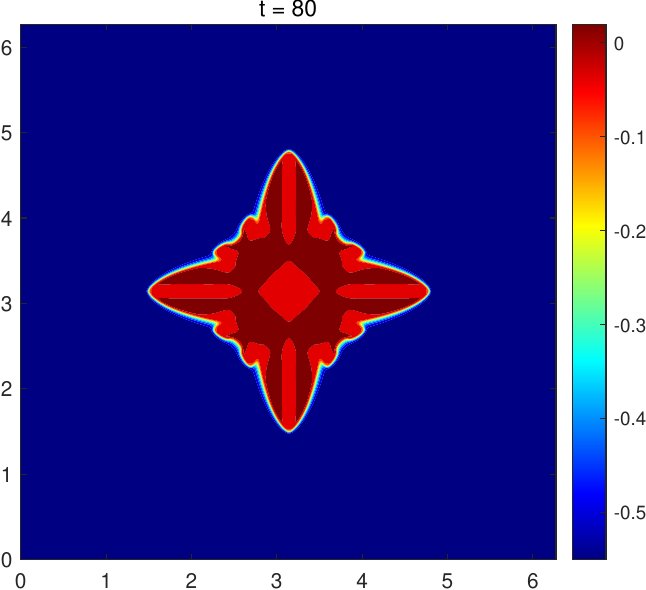}
		\end{minipage}
		\begin{minipage}[t]{0.185\linewidth}
			\centering
			\includegraphics[width=1\linewidth]{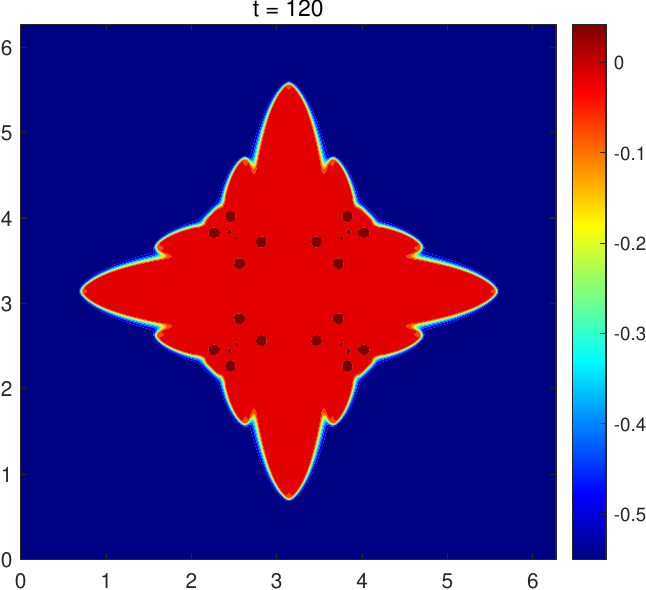}
		\end{minipage}
	}
	
	\caption{(Example \ref{exp.4}) Dynamics driven by the dendritic crystal growth model with fourfold anisotropy $\beta=4$ and the parameter $K=0.6$ using the adaptive variable time step R-GSAV BDF2 scheme.}
	\label{fig:fourfoldK06}
	\label{fig:fourfoldK06snap}
\end{figure}

\begin{figure}[htbp]

	\vspace{-0.35cm}
	\subfigtopskip=2pt
	\subfigbottomskip=2pt
	\subfigcapskip=-5pt

	\centering	

    \subfigure[\scriptsize (a) Snapshots of the phase field $\phi$.]{

		\begin{minipage}[t]{0.185\linewidth}
			\centering
			\includegraphics[width=1\linewidth]{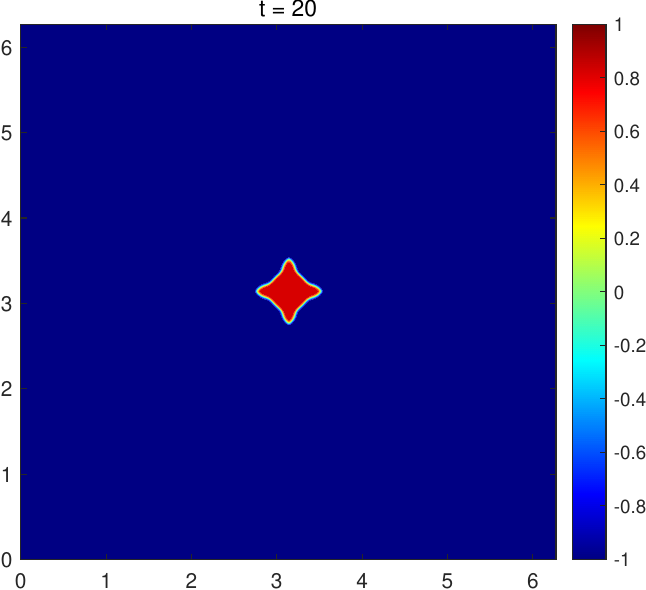}
		\end{minipage}
		\begin{minipage}[t]{0.185\linewidth}
			\centering
			\includegraphics[width=1\linewidth]{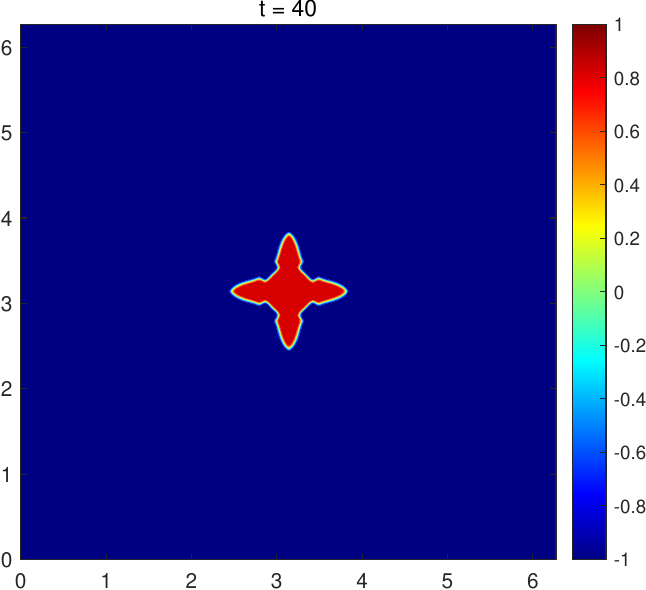}
		\end{minipage}
		\begin{minipage}[t]{0.185\linewidth}
			\centering
			\includegraphics[width=1\linewidth]{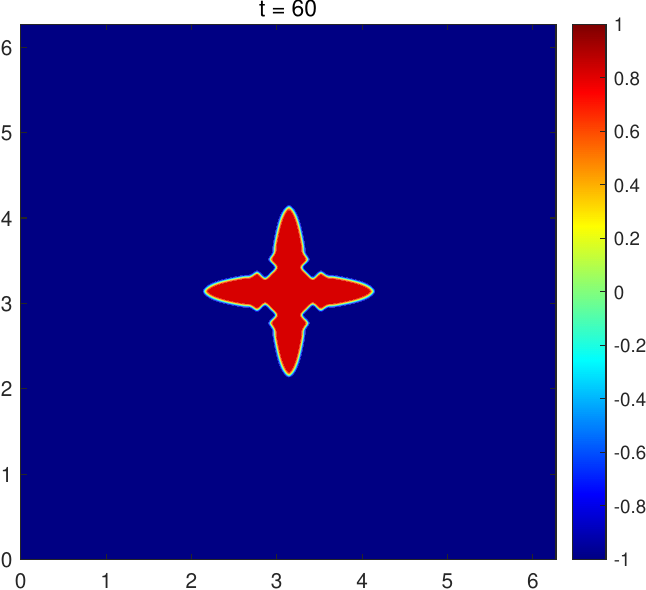}
		\end{minipage}
		\begin{minipage}[t]{0.185\linewidth}
			\centering
			\includegraphics[width=1\linewidth]{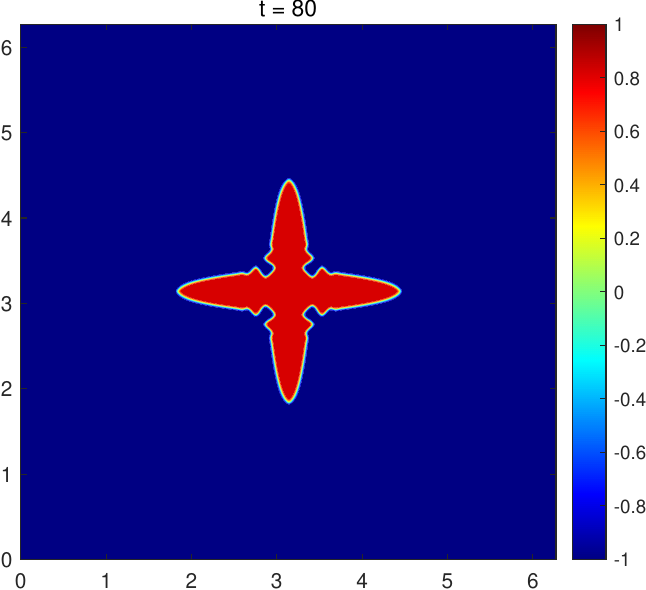}
		\end{minipage}
		\begin{minipage}[t]{0.185\linewidth}
			\centering
			\includegraphics[width=1\linewidth]{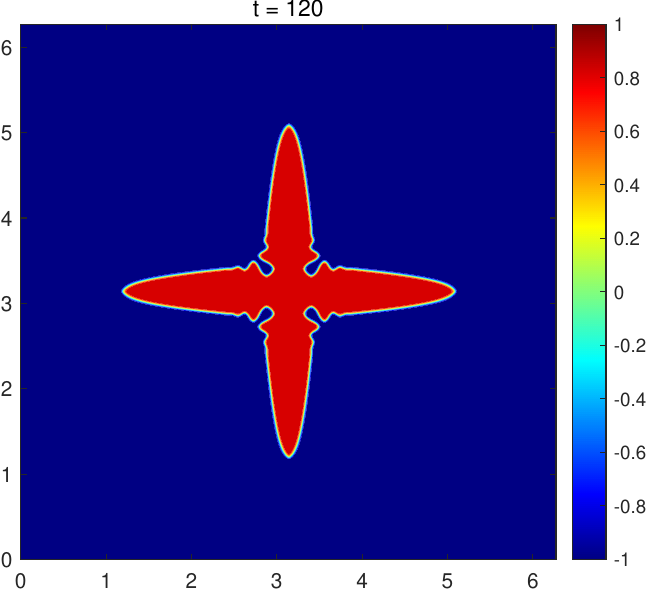}
		\end{minipage}
	}

    \subfigure[\scriptsize (b) Snapshots of the temperature field $U$.]{

		\begin{minipage}[t]{0.185\linewidth}
			\centering
			\includegraphics[width=1\linewidth]{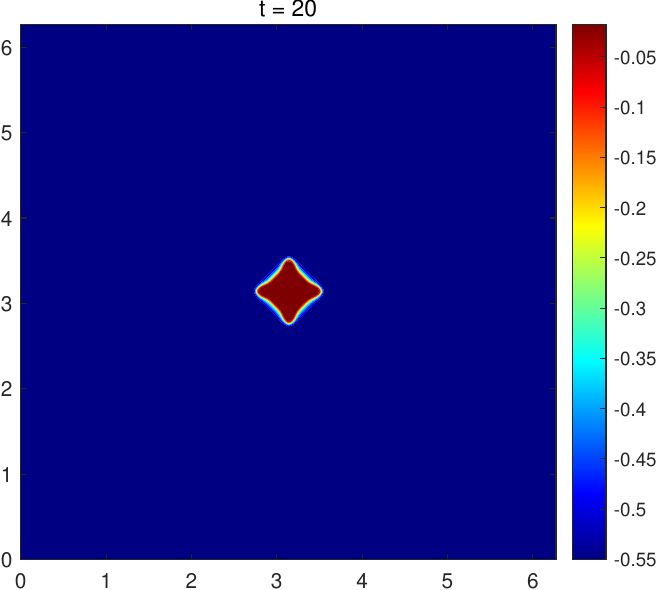}
		\end{minipage}
		\begin{minipage}[t]{0.185\linewidth}
			\centering
			\includegraphics[width=1\linewidth]{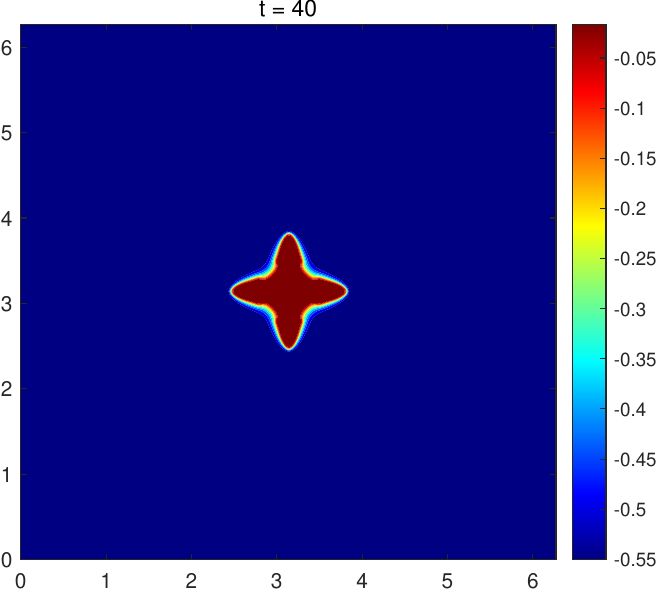}
		\end{minipage}
		\begin{minipage}[t]{0.185\linewidth}
			\centering
			\includegraphics[width=1\linewidth]{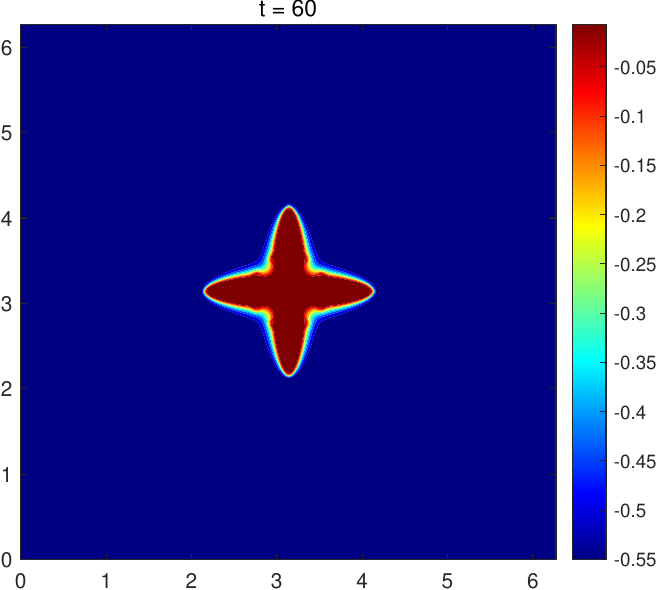}
		\end{minipage}
		\begin{minipage}[t]{0.185\linewidth}
			\centering
			\includegraphics[width=1\linewidth]{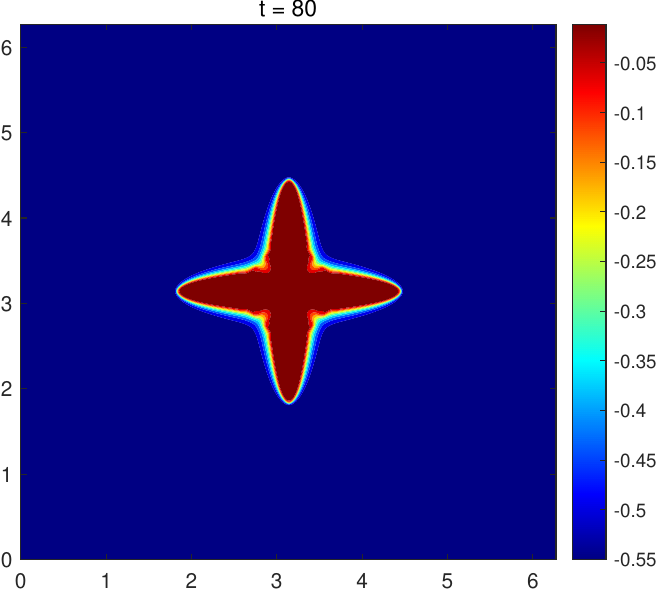}
		\end{minipage}
		\begin{minipage}[t]{0.185\linewidth}
			\centering
			\includegraphics[width=1\linewidth]{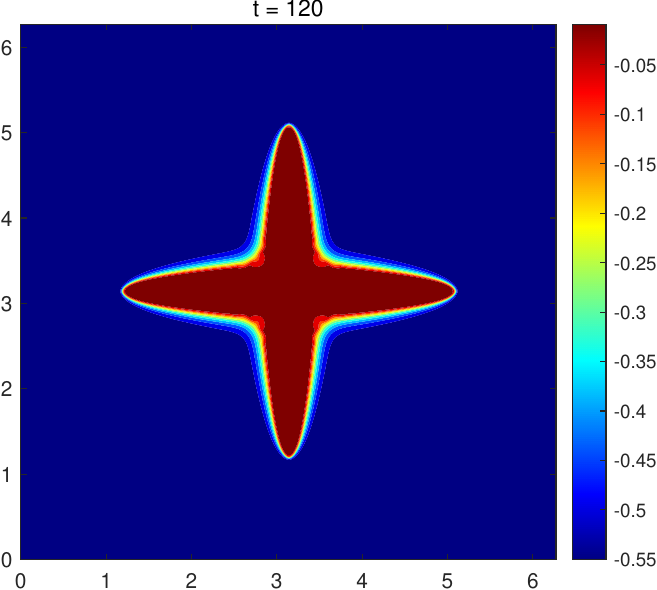}
		\end{minipage}
	}
	
	\caption{(Example \ref{exp.4}) Dynamics driven by the dendritic crystal growth model with fourfold anisotropy $\beta=4$ and the parameter $K=0.7$ using the adaptive variable time step R-GSAV BDF2 scheme.}
	\label{fig:fourfoldK07}
\end{figure}

\begin{figure}[htbp]

	\vspace{-0.35cm}
	\subfigtopskip=2pt
	\subfigbottomskip=2pt
	\subfigcapskip=-5pt
	\centering	

    \subfigure[\scriptsize (a) Snapshots of the phase field $\phi$.]{

		\begin{minipage}[t]{0.185\linewidth}
			\centering
			\includegraphics[width=1\linewidth]{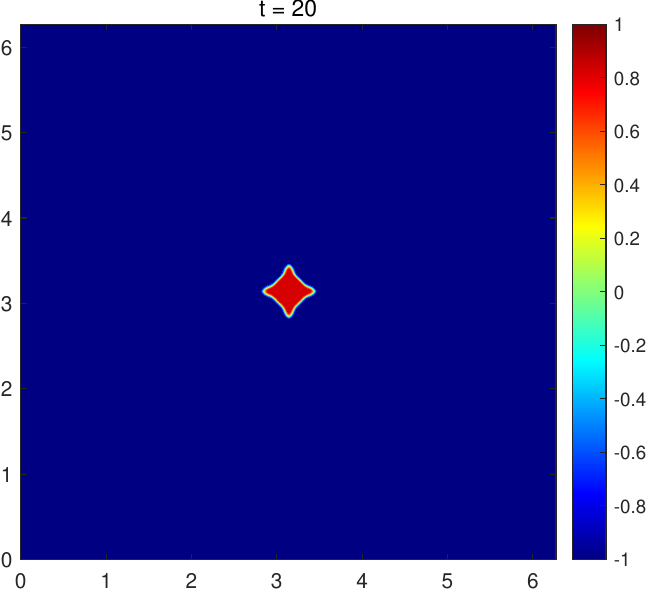}
		\end{minipage}
		\begin{minipage}[t]{0.185\linewidth}
			\centering
			\includegraphics[width=1\linewidth]{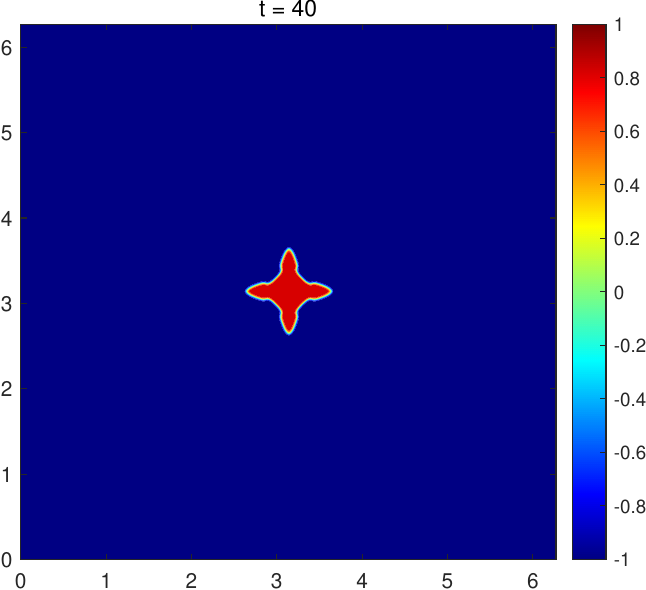}
		\end{minipage}
		\begin{minipage}[t]{0.185\linewidth}
			\centering
			\includegraphics[width=1\linewidth]{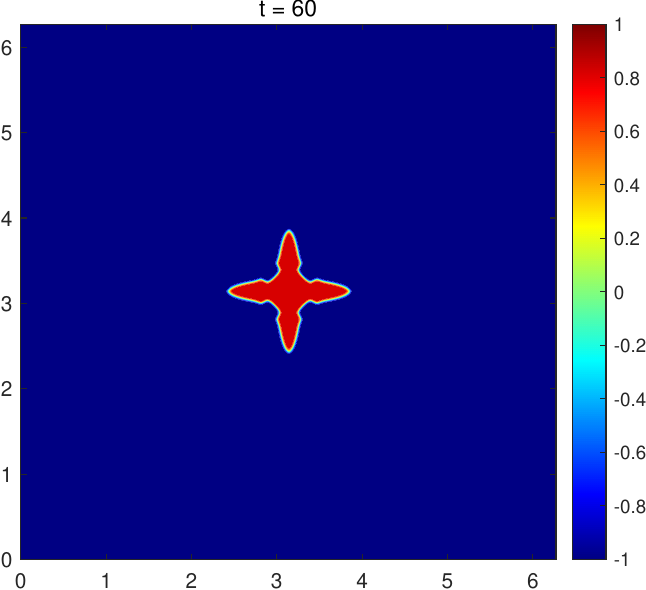}
		\end{minipage}
		\begin{minipage}[t]{0.185\linewidth}
			\centering
			\includegraphics[width=1\linewidth]{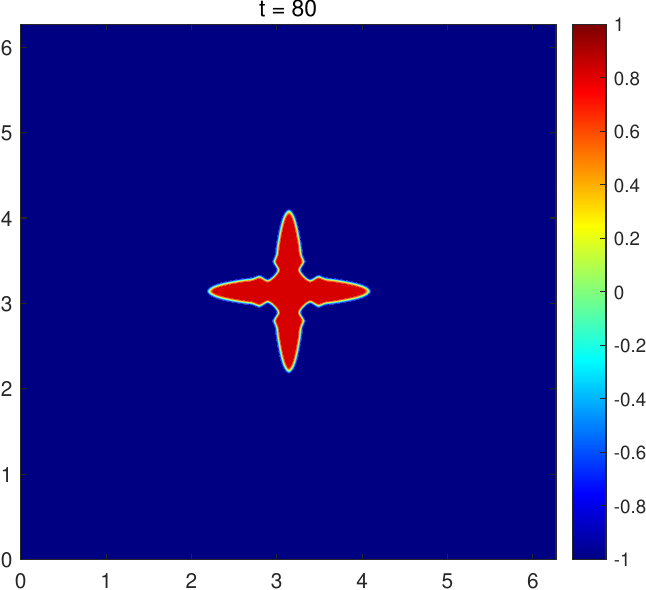}
		\end{minipage}
		\begin{minipage}[t]{0.185\linewidth}
			\centering
			\includegraphics[width=1\linewidth]{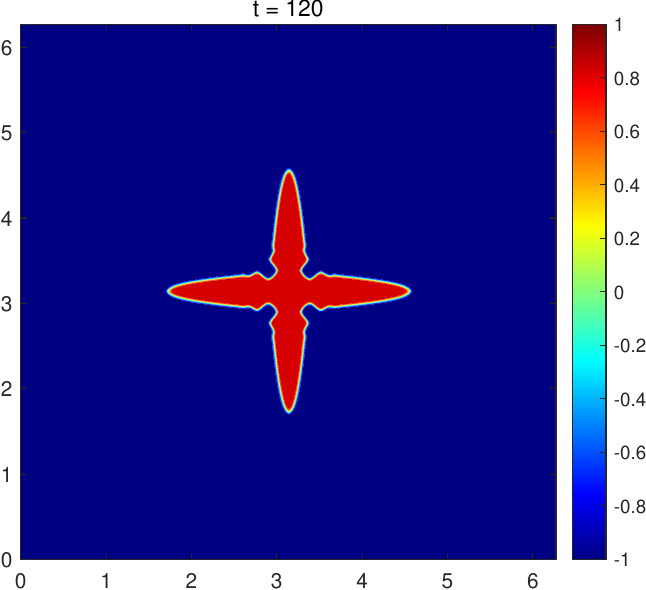}
		\end{minipage}
	}

    \subfigure[\scriptsize (b) Snapshots of the temperature field $U$.]{

		\begin{minipage}[t]{0.185\linewidth}
			\centering
			\includegraphics[width=1\linewidth]{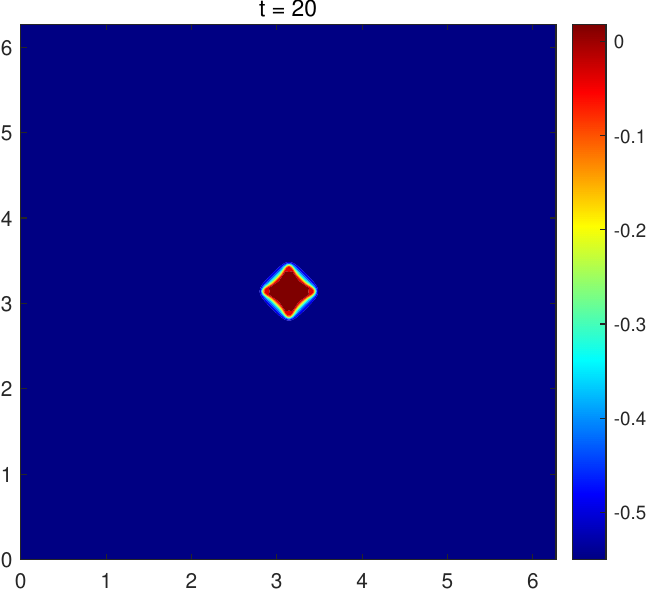}
		\end{minipage}
		\begin{minipage}[t]{0.185\linewidth}
			\centering
			\includegraphics[width=1\linewidth]{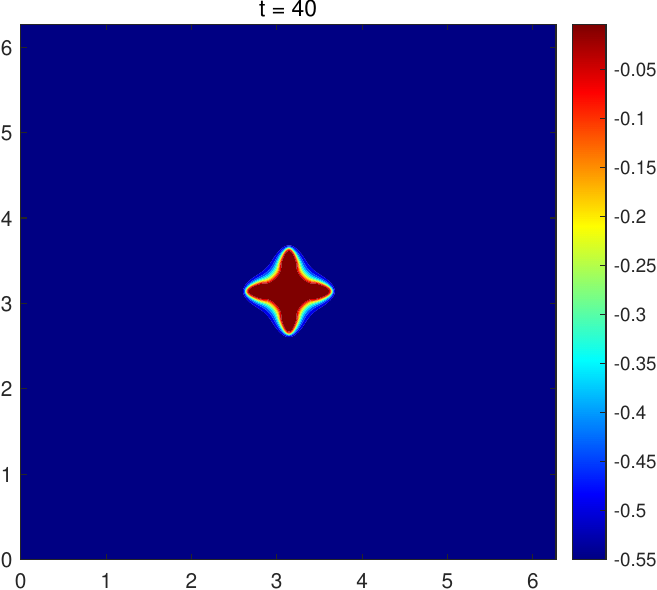}
		\end{minipage}
		\begin{minipage}[t]{0.185\linewidth}
			\centering
			\includegraphics[width=1\linewidth]{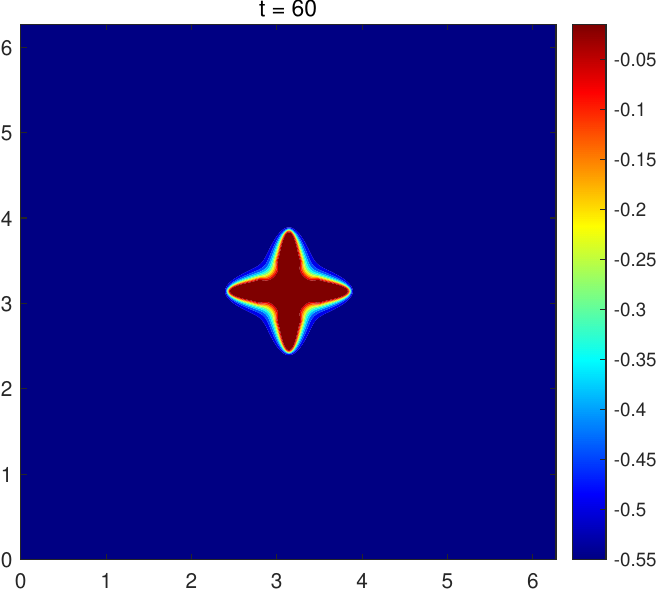}
		\end{minipage}
		\begin{minipage}[t]{0.185\linewidth}
			\centering
			\includegraphics[width=1\linewidth]{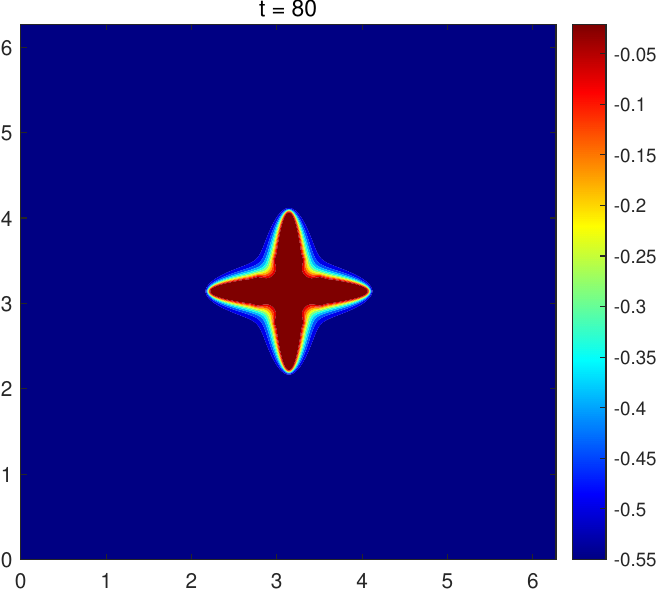}
		\end{minipage}
		\begin{minipage}[t]{0.185\linewidth}
			\centering
			\includegraphics[width=1\linewidth]{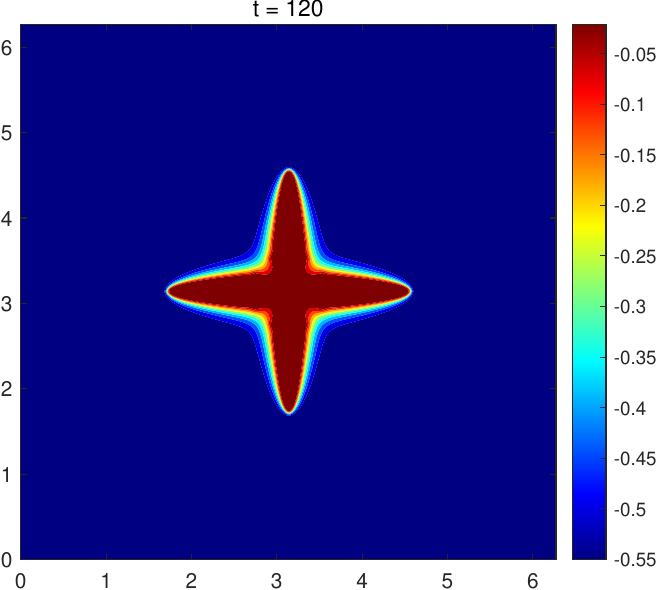}
		\end{minipage}
	}

	\caption{(Example \ref{exp.4}) Dynamics driven by the dendritic crystal growth model with fourfold anisotropy $\beta=4$ and the parameter $K=0.8$ using the adaptive variable time step R-GSAV BDF2 scheme.}
	\label{fig:fourfoldK08}
\end{figure}

\begin{figure}[htbp]
\subfigure[\scriptsize (a) Energy evolution]{
\begin{minipage}{0.45 \textwidth}
\centering
\includegraphics[width=\textwidth]{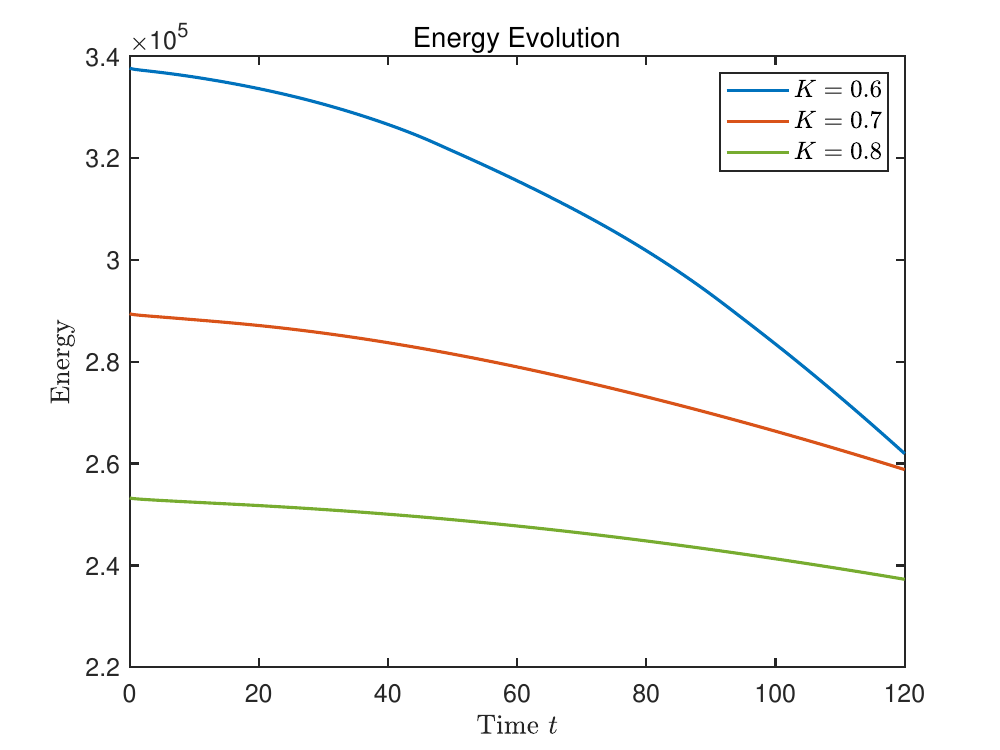}
\end{minipage}
}
\subfigure[\scriptsize (b) Crystal area]{
\begin{minipage}{0.45 \textwidth}
\centering
\includegraphics[width=\textwidth]{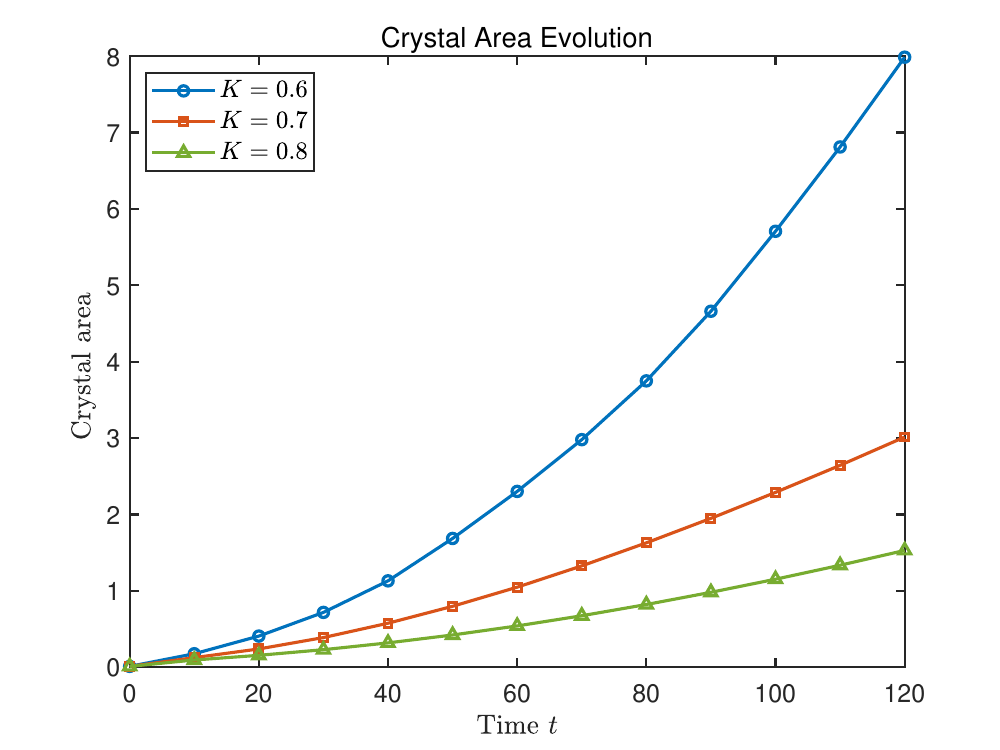}
\end{minipage}
}

\subfigure[\scriptsize (c) Adaptive time step]{
\begin{minipage}{0.45 \textwidth}
\centering
\includegraphics[width=\textwidth]{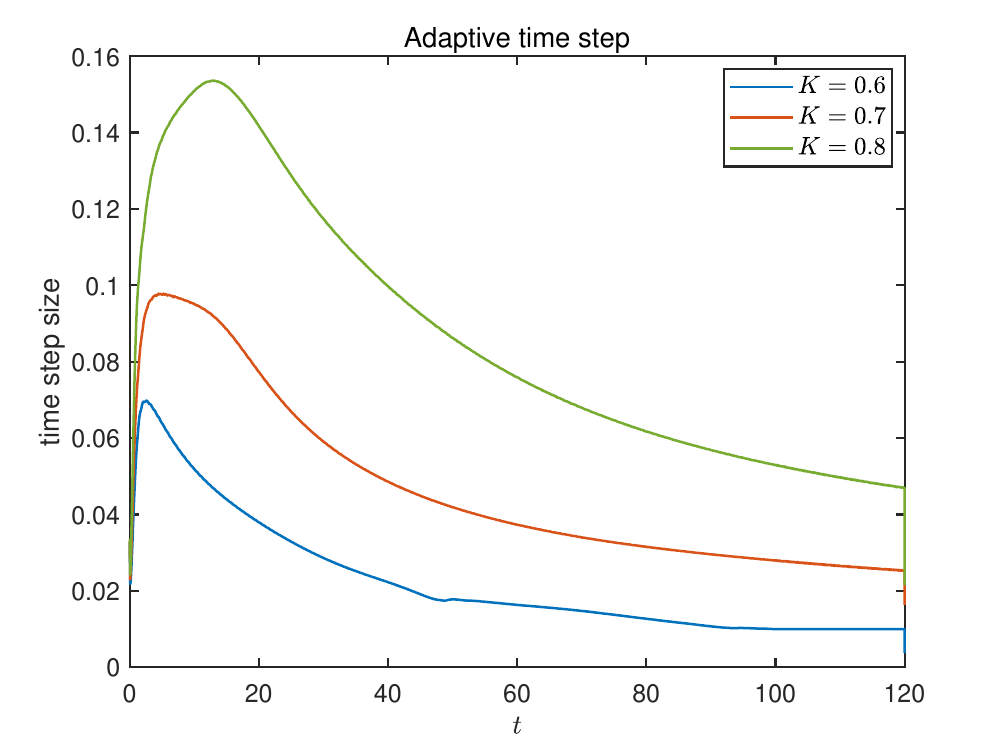}
\end{minipage}
}
\caption{(Example \ref{exp.4}) (a) Time evolution of the energy. (b) Crystal
area $\int_{\Omega}\frac{1+\phi}{2}dx$. (c) The evolution of the adaptive time step sizes for different values of the parameter $K$ with fourfold anisotropy $\beta=4$. }
\label{fourfold_energy}
\end{figure}

We plot in Figures \ref{fig:fourfoldK06snap}--\ref{fig:fourfoldK08} some snapshots of the phase field $\phi$ and the temperature field $U$ computed for the values $0.6$, $0.7$ and $0.8$ of the latent heat parameter $K$. We observe from these figures that the crystal surface grows in time, and a star shape with four branches is always formed in all cases starting with a tiny nucleus. However, the latent heat parameter has a significant impact on the width of the branches of the growing crystals: the larger $K$ is, the thinner the width of the branches is, and the sharper the tips are. This is consistent with the results reported in \cite{yang2020efficient,guo2024efficient}. The snapshots of the temperature field of each simulation are also shown in Figures \ref{fig:fourfoldK06snap}--\ref{fig:fourfoldK08}. We see that the contours of $U$ take the same dendrite crystal shape as the phase field. The physical explanation of this observation is that the heat propagates only at the interface. The dissipation behavior of the energy in time is shown in Figure \ref{fourfold_energy}(a). The decay feature of the energy for all tested $K$ reflects the good stability property of the proposed scheme employed in the calculation. Figure \ref{fourfold_energy}(b) displays the crystal area, defined by the quantity $\int_{\Omega}\frac{1+\phi}{2}dx$, indicating that the area of the crystal increases in time and with decreasing $K$. The corresponding evolution of the adaptive time step sizes is plotted in Figure \ref{fourfold_energy}(c). It is observed that, after a short initial stage, the time step size decreases monotonically and gradually approaches $\tau_{\min}$ as time evolves. This is because the step size selected by the energy based strategy \eqref{eq:strategy} is inversely proportional to the energy decay rate $|E'(t)|$: as the dendrite grows, the crystal area and, more importantly, the total length of the solid liquid interface keep increasing due to the continuous development of the four branches and the sidebranches, so that more latent heat is released per unit time and the total energy decays faster and faster (cf. Figure \ref{fourfold_energy}(a)). The strategy \eqref{eq:strategy} thus automatically refines the time step to accurately capture the accelerating growth dynamics, until the step size saturates at the lower bound $\tau_{\min}$. Moreover, the numbers of time steps of the uniform mesh with $\tau=\tau_{\min}=0.01$ and the adaptive strategy \eqref{eq:strategy} are compared in Table \ref{tab:steps43}: the numbers of time steps are reduced by $38.1\%$, $74.4\%$ and $86.6\%$ for $K=0.6$, $0.7$ and $0.8$, respectively.

\begin{table}[htbp]
\centering
\caption{(Example \ref{exp.4}) Comparison of the numbers of time steps up to $T=120$: the uniform mesh with $\tau=0.01$ and the adaptive strategy \eqref{eq:strategy}.}
\label{tab:steps43}
\begin{tabular}{c|c|c|c}
  \hline
 $K$ & uniform $\tau=0.01$ & adaptive strategy \eqref{eq:strategy} & step reduction \\
  \hline
  0.6 & 12000 & 7424 & 38.1\% \\
  0.7 & 12000 & 3068 & 74.4\% \\
  0.8 & 12000 & 1608 & 86.6\% \\
  \hline
\end{tabular}
\end{table}

\subsection{Sixfold dendrite crystal growth in 2D}
We next turn to the sixfold case ($\beta=6$) to examine the influence of the fold number of anisotropy on the growth pattern. All the discretization details, i.e., the proposed adaptive variable step R-GSAV BDF2 scheme \eqref{2.4} with the strategy \eqref{eq:strategy}, the corresponding adaptive parameters and the spatial resolution, are exactly the same as those in the previous subsection.

\begin{example}\label{exp.5}
\rm{We consider the sixfold dendrite crystal growth with the same initial conditions \eqref{e4} and the same parameters as in Example \ref{exp.4}, except that the number of folds of anisotropy is taken as $\beta=6$.}
\end{example}

\begin{figure}[htbp]
	\vspace{-0.35cm}
	\subfigtopskip=2pt
	\subfigbottomskip=2pt
	\subfigcapskip=-5pt
	\centering	
    \subfigure[\scriptsize (a) Snapshots of the phase field $\phi$.]{

		\begin{minipage}[t]{0.185\linewidth}
			\centering
			\includegraphics[width=1\linewidth]{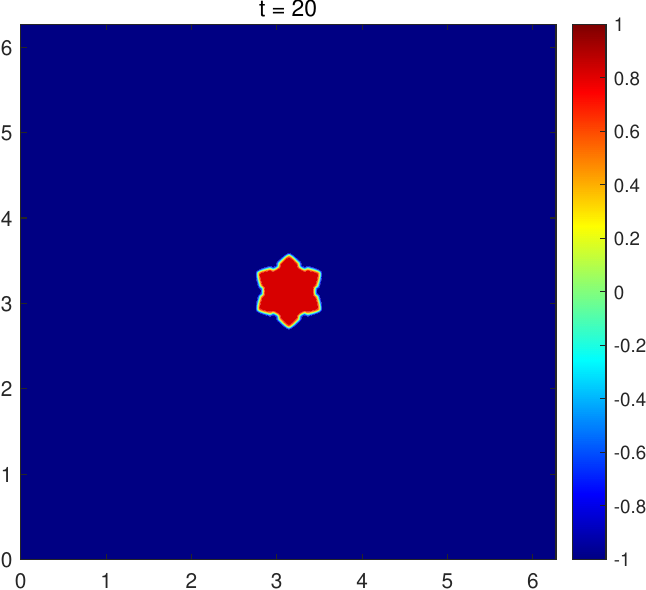}
		\end{minipage}
		\begin{minipage}[t]{0.185\linewidth}
			\centering
			\includegraphics[width=1\linewidth]{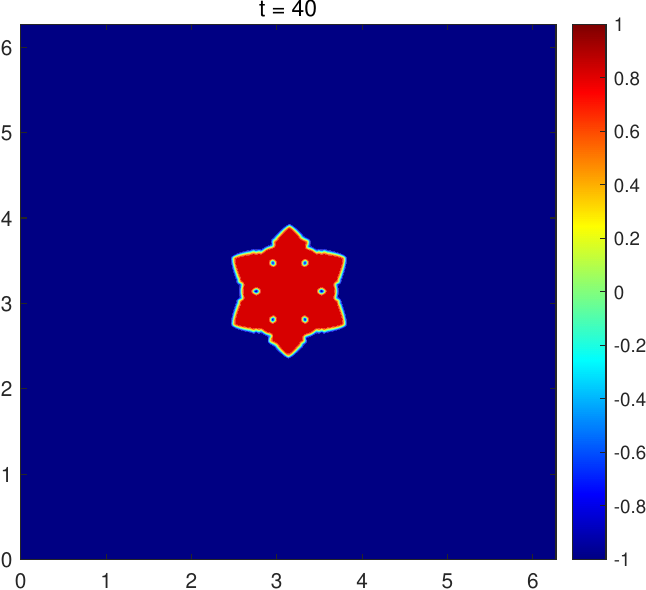}
		\end{minipage}
		\begin{minipage}[t]{0.185\linewidth}
			\centering
			\includegraphics[width=1\linewidth]{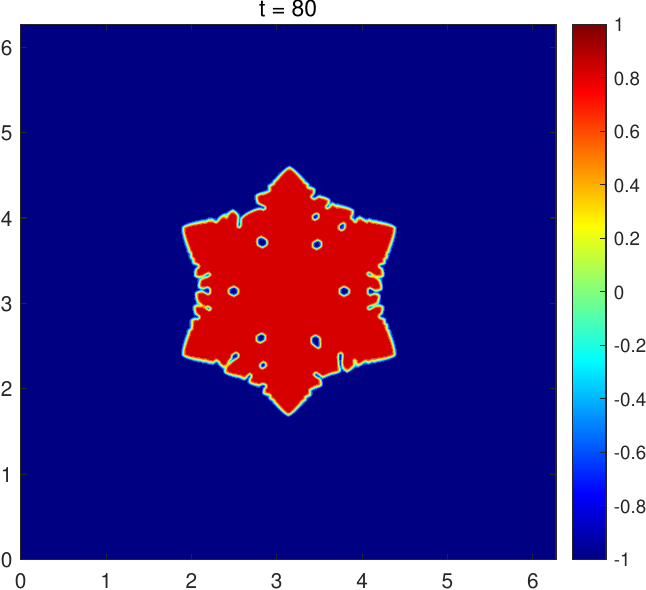}
		\end{minipage}
		\begin{minipage}[t]{0.185\linewidth}
			\centering
			\includegraphics[width=1\linewidth]{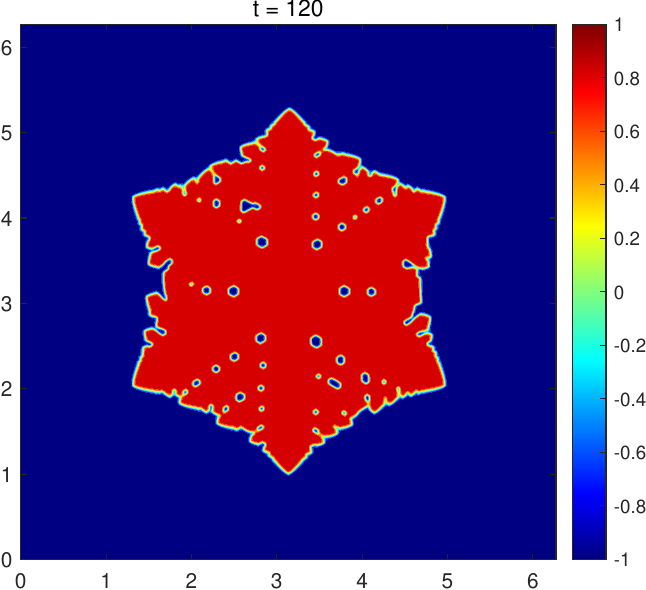}
		\end{minipage}
		\begin{minipage}[t]{0.185\linewidth}
			\centering
			\includegraphics[width=1\linewidth]{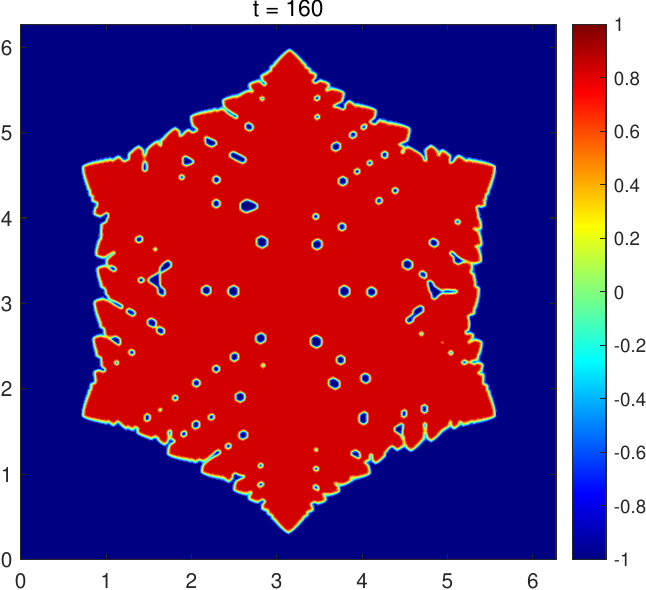}
		\end{minipage}
	}
    \subfigure[\scriptsize (b) Snapshots of the temperature field $U$.]{
		\begin{minipage}[t]{0.185\linewidth}
			\centering
			\includegraphics[width=1\linewidth]{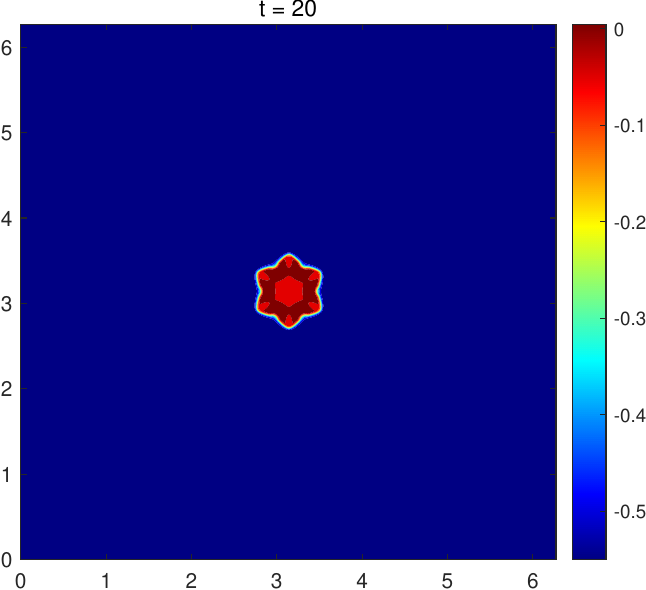}
		\end{minipage}
		\begin{minipage}[t]{0.185\linewidth}
			\centering
			\includegraphics[width=1\linewidth]{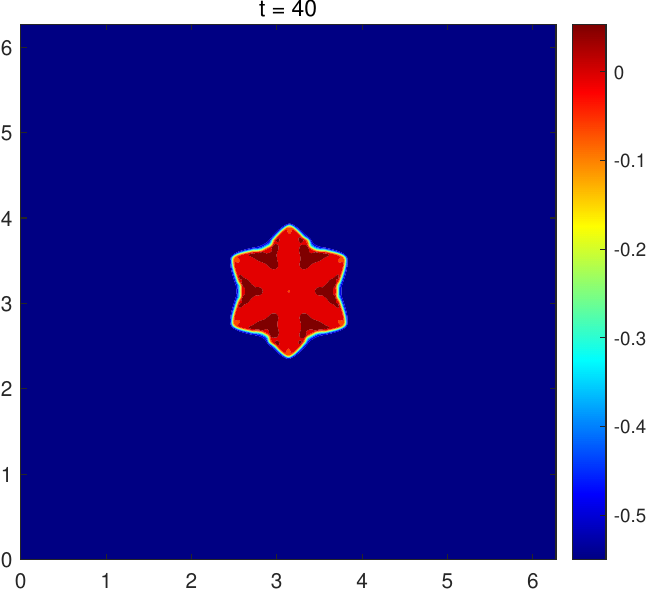}
		\end{minipage}
		\begin{minipage}[t]{0.185\linewidth}
			\centering
			\includegraphics[width=1\linewidth]{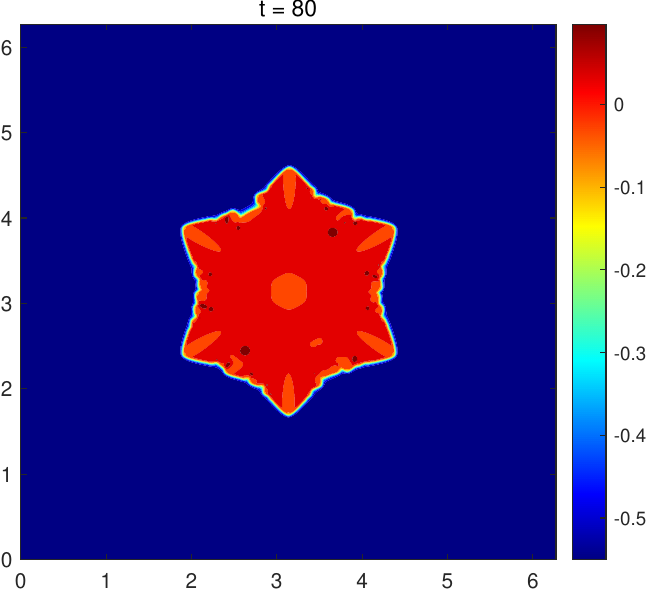}
		\end{minipage}
		\begin{minipage}[t]{0.185\linewidth}
			\centering
			\includegraphics[width=1\linewidth]{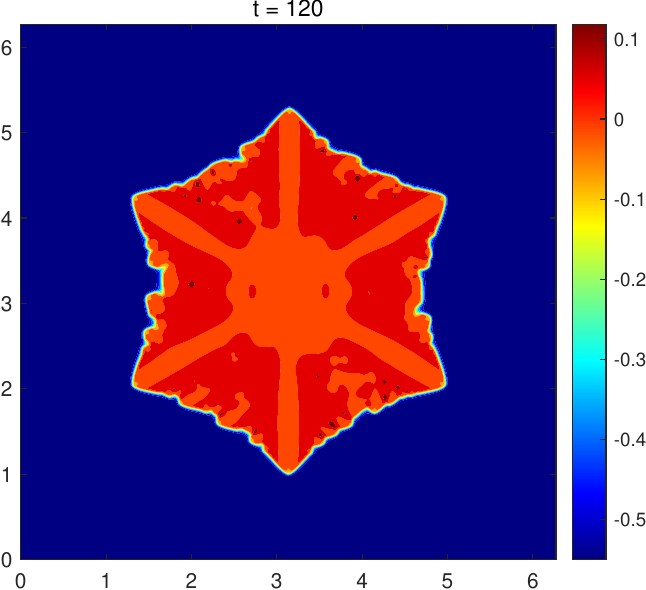}
		\end{minipage}
		\begin{minipage}[t]{0.185\linewidth}
			\centering
			\includegraphics[width=1\linewidth]{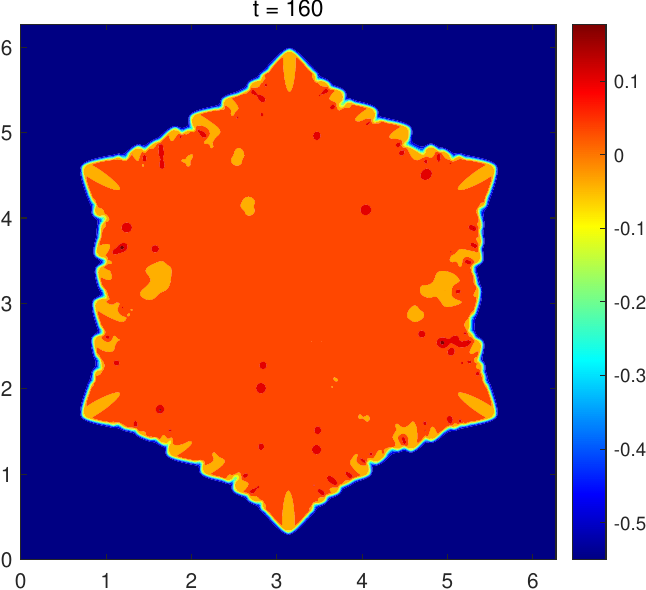}
		\end{minipage}
	}
	\caption{(Example \ref{exp.5}) Dynamics driven by the dendritic crystal growth model with sixfold anisotropy $\beta=6$ and the parameter $K=0.6$ using the adaptive variable time step R-GSAV BDF2 scheme.}
	\label{fig:sixfoldK06}
\end{figure}

\begin{figure}[htbp]
	\vspace{-0.35cm}
	\subfigtopskip=2pt
	\subfigbottomskip=2pt
	\subfigcapskip=-5pt
	\centering	
    \subfigure[\scriptsize (a) Snapshots of the phase field $\phi$.]{
		\begin{minipage}[t]{0.185\linewidth}
			\centering
			\includegraphics[width=1\linewidth]{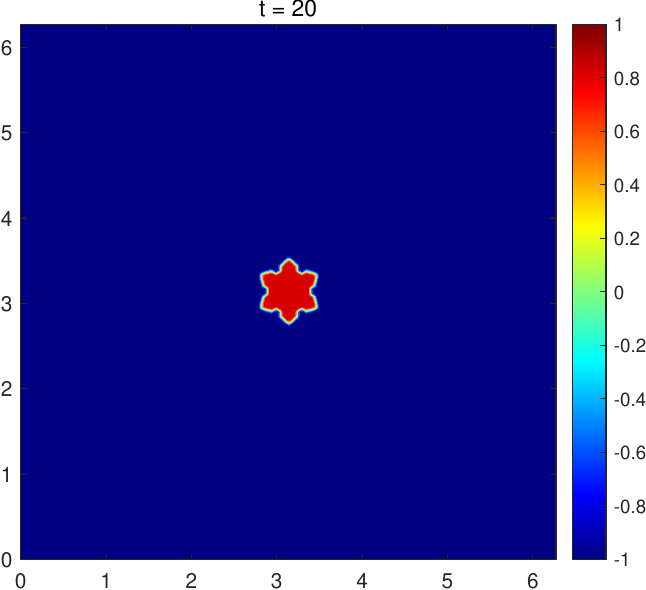}
		\end{minipage}
		\begin{minipage}[t]{0.185\linewidth}
			\centering
			\includegraphics[width=1\linewidth]{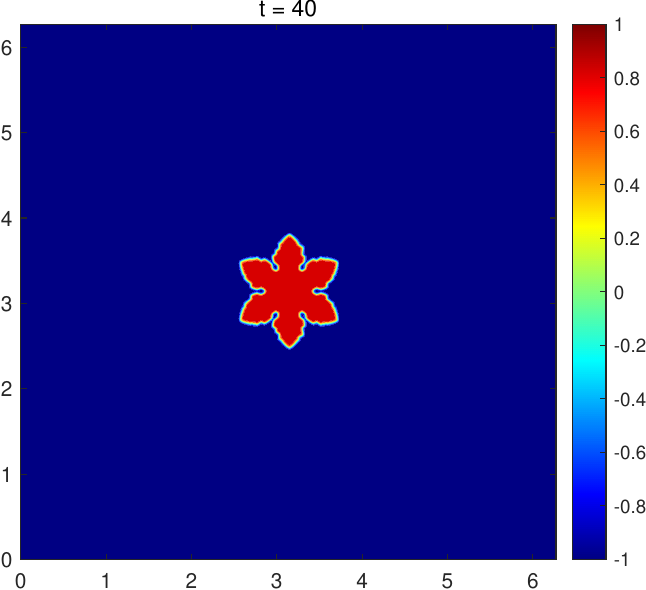}
		\end{minipage}
		\begin{minipage}[t]{0.185\linewidth}
			\centering
			\includegraphics[width=1\linewidth]{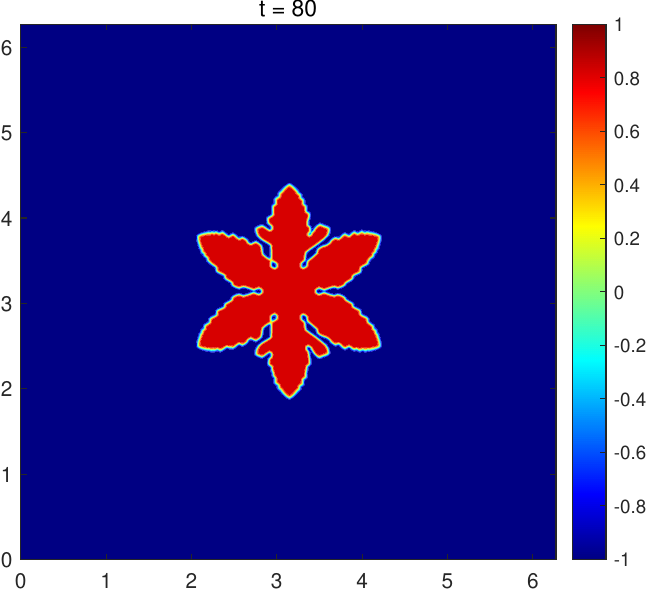}
		\end{minipage}
		\begin{minipage}[t]{0.185\linewidth}
			\centering
			\includegraphics[width=1\linewidth]{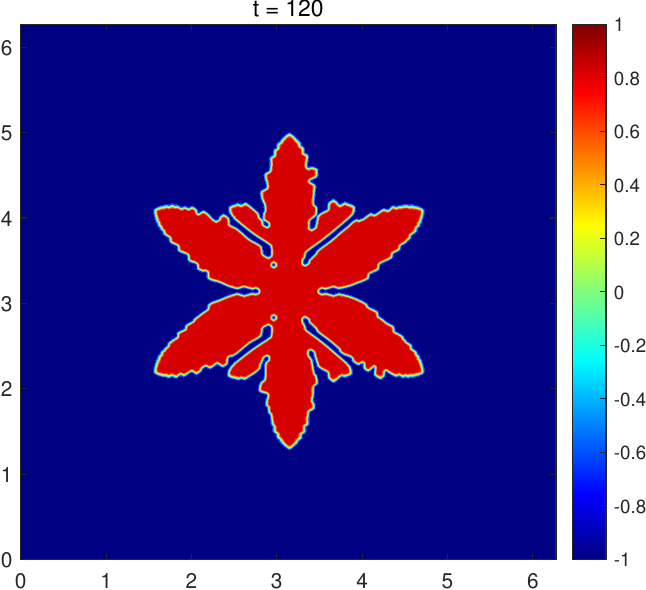}
		\end{minipage}
		\begin{minipage}[t]{0.185\linewidth}
			\centering
			\includegraphics[width=1\linewidth]{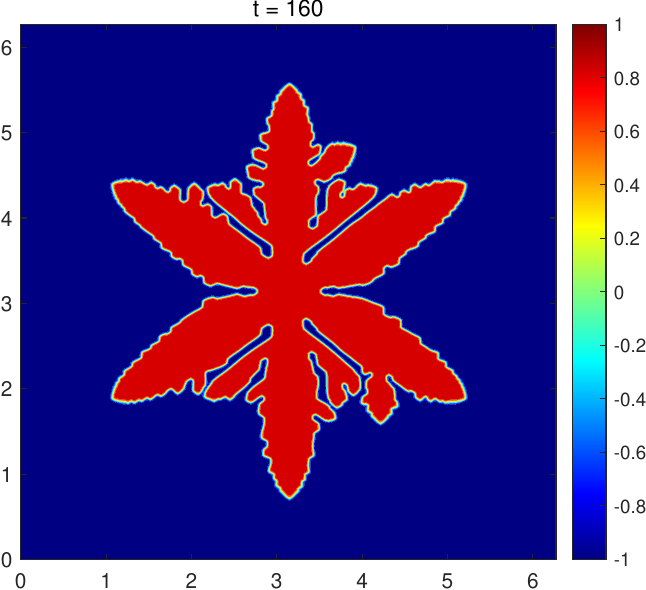}
		\end{minipage}
	}
    \subfigure[\scriptsize (b) Snapshots of the temperature field $U$.]{
		\begin{minipage}[t]{0.185\linewidth}
			\centering
			\includegraphics[width=1\linewidth]{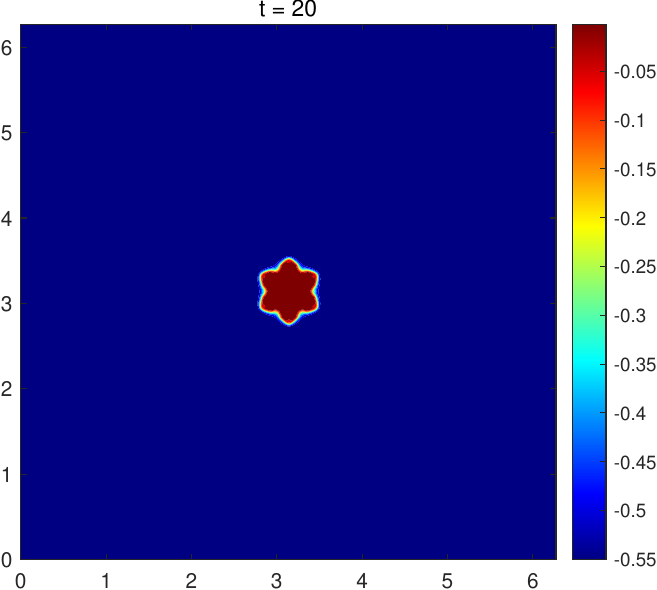}
		\end{minipage}
		\begin{minipage}[t]{0.185\linewidth}
			\centering
			\includegraphics[width=1\linewidth]{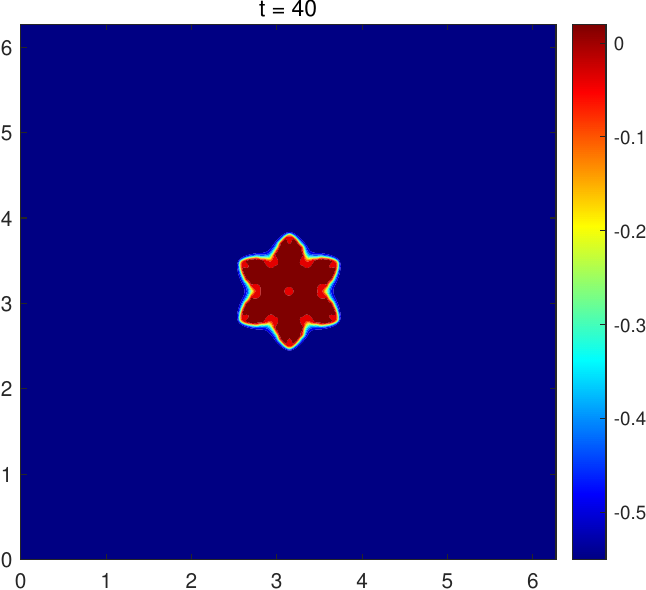}
		\end{minipage}
		\begin{minipage}[t]{0.185\linewidth}
			\centering
			\includegraphics[width=1\linewidth]{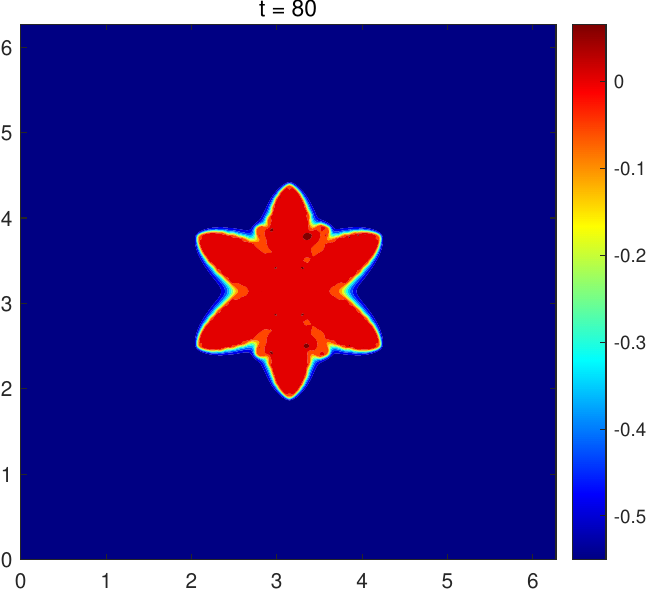}
		\end{minipage}
		\begin{minipage}[t]{0.185\linewidth}
			\centering
			\includegraphics[width=1\linewidth]{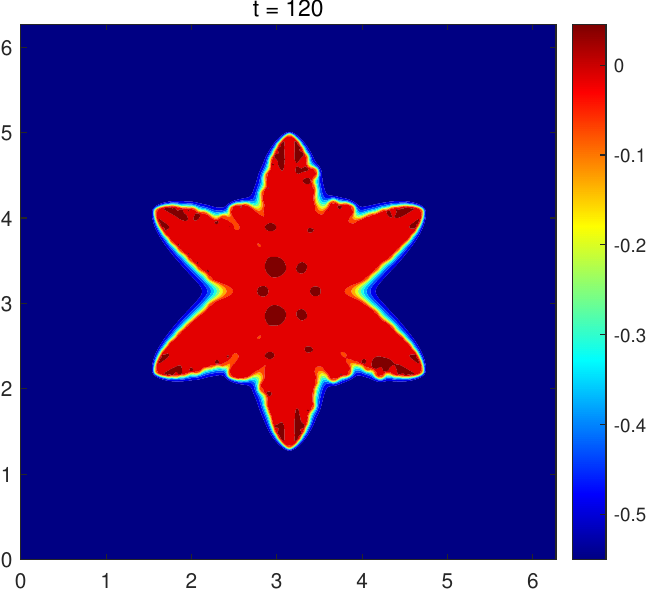}
		\end{minipage}
		\begin{minipage}[t]{0.185\linewidth}
			\centering
			\includegraphics[width=1\linewidth]{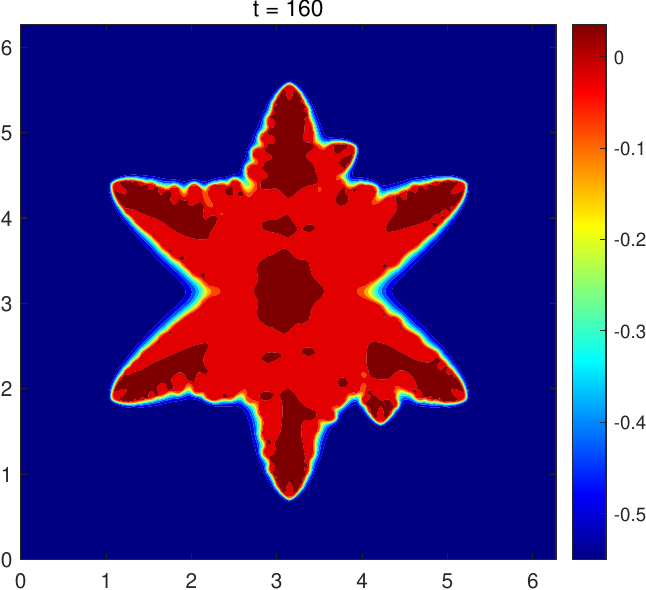}
		\end{minipage}
	}
	\caption{(Example \ref{exp.5}) Dynamics driven by the dendritic crystal growth model with sixfold anisotropy $\beta=6$ and the parameter $K=0.65$ using the adaptive variable time step R-GSAV BDF2 scheme.}
	\label{fig:sixfoldK065}
\end{figure}

\begin{figure}[htbp]
	\vspace{-0.35cm}
	\subfigtopskip=2pt
	\subfigbottomskip=2pt
	\subfigcapskip=-5pt
	\centering	
    \subfigure[\scriptsize (a) Snapshots of the phase field $\phi$.]{
		\begin{minipage}[t]{0.185\linewidth}
			\centering
			\includegraphics[width=1\linewidth]{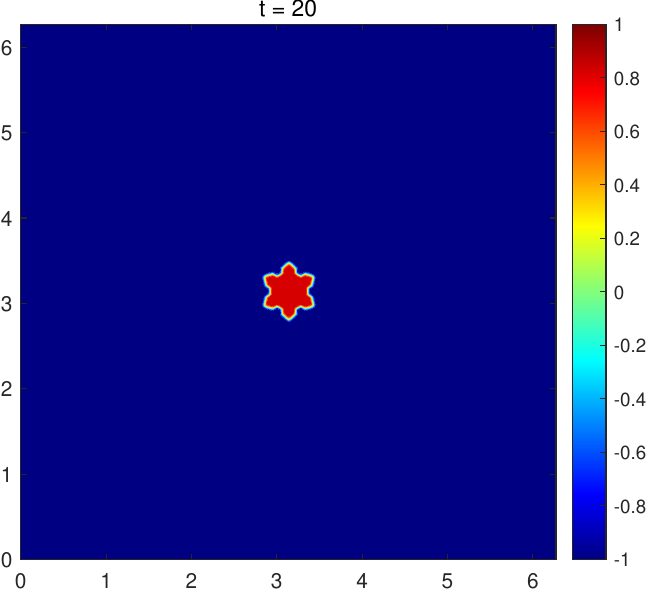}
		\end{minipage}
		\begin{minipage}[t]{0.185\linewidth}
			\centering
			\includegraphics[width=1\linewidth]{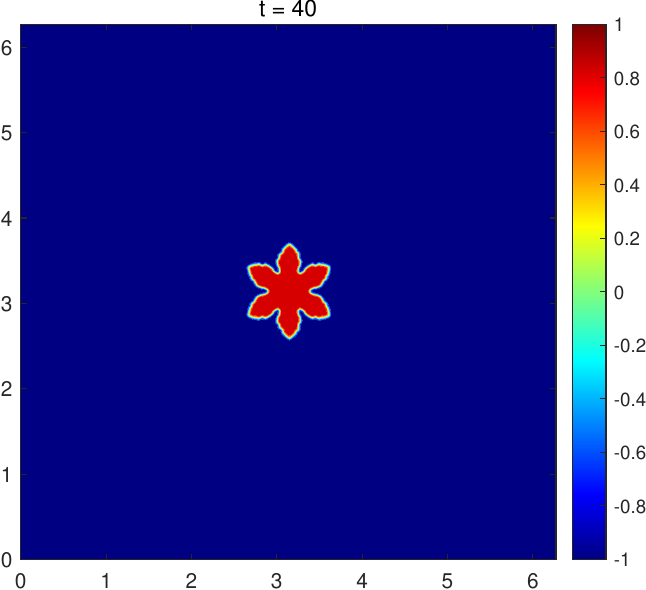}
		\end{minipage}
		\begin{minipage}[t]{0.185\linewidth}
			\centering
			\includegraphics[width=1\linewidth]{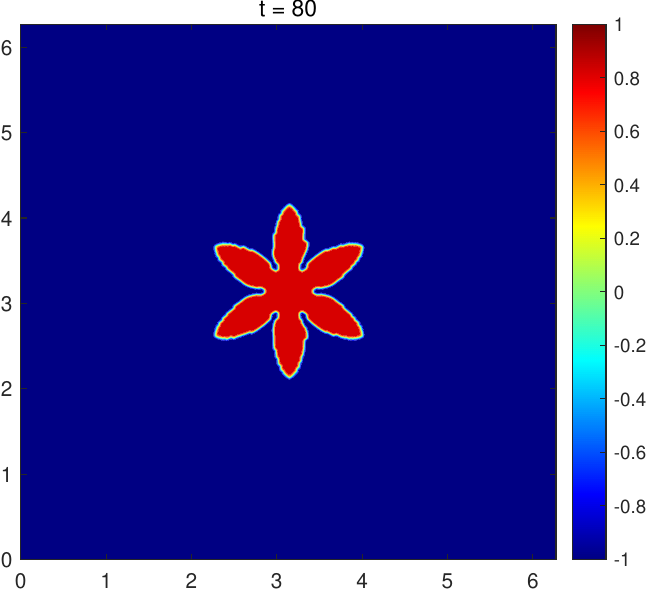}
		\end{minipage}
		\begin{minipage}[t]{0.185\linewidth}
			\centering
			\includegraphics[width=1\linewidth]{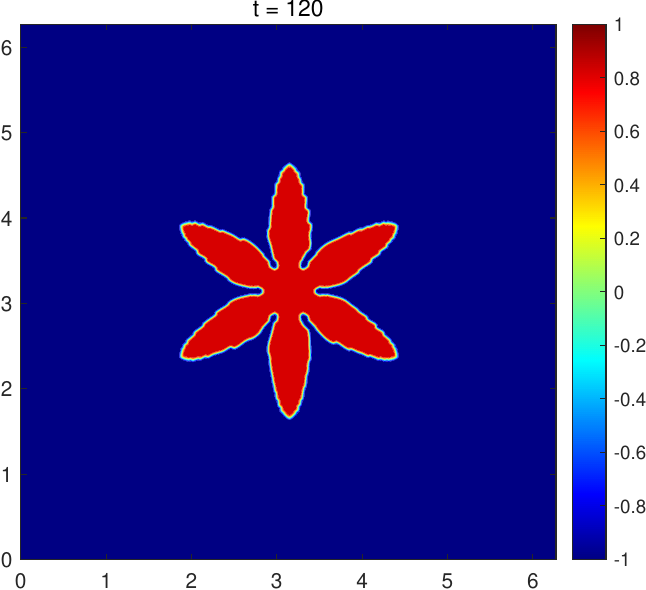}
		\end{minipage}
		\begin{minipage}[t]{0.185\linewidth}
			\centering
			\includegraphics[width=1\linewidth]{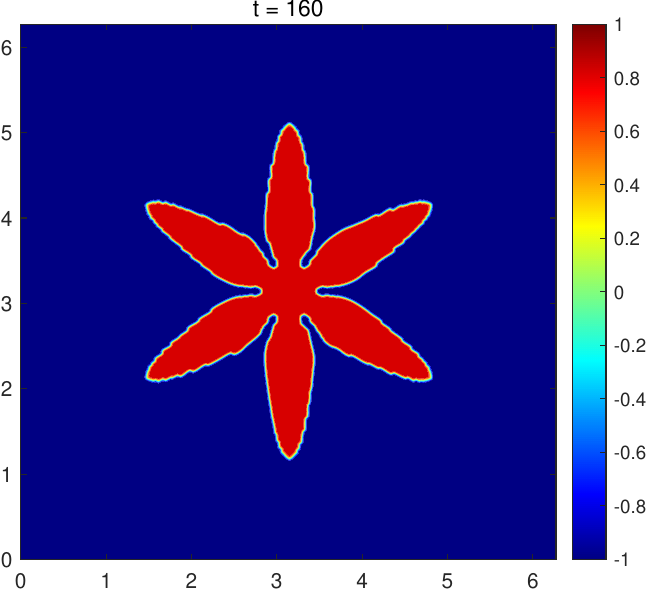}
		\end{minipage}
	}
    \subfigure[\scriptsize (b) Snapshots of the temperature field $U$.]{
		\begin{minipage}[t]{0.185\linewidth}
			\centering
			\includegraphics[width=1\linewidth]{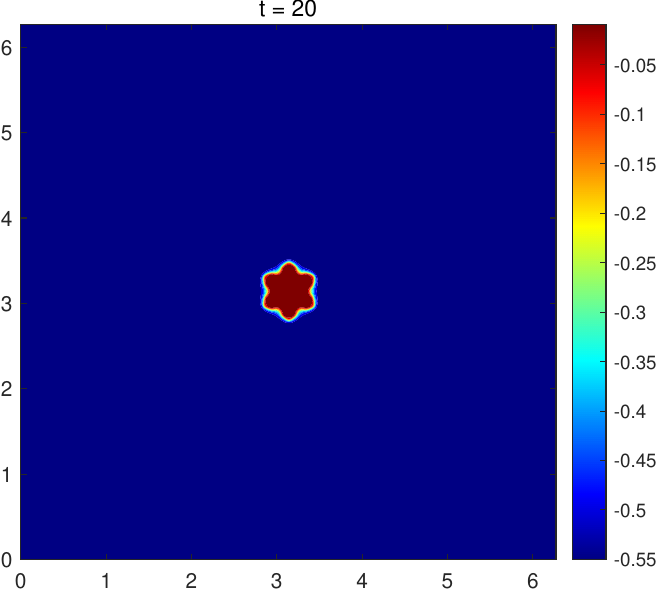}
		\end{minipage}
		\begin{minipage}[t]{0.185\linewidth}
			\centering
			\includegraphics[width=1\linewidth]{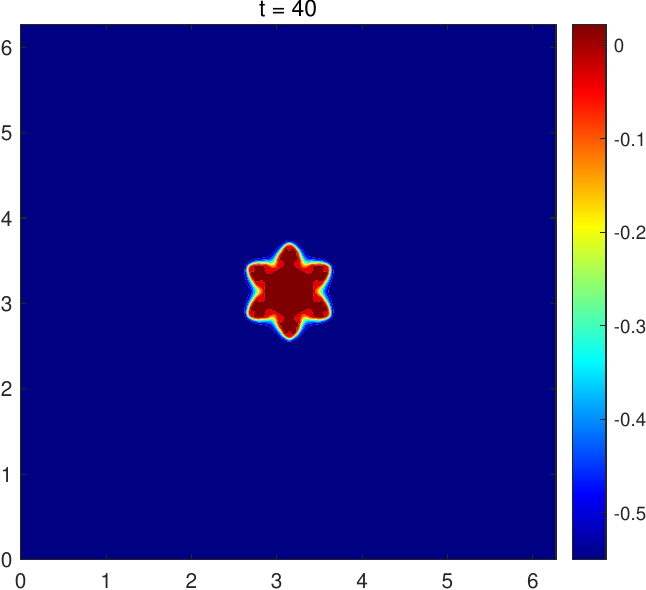}
		\end{minipage}
		\begin{minipage}[t]{0.185\linewidth}
			\centering
			\includegraphics[width=1\linewidth]{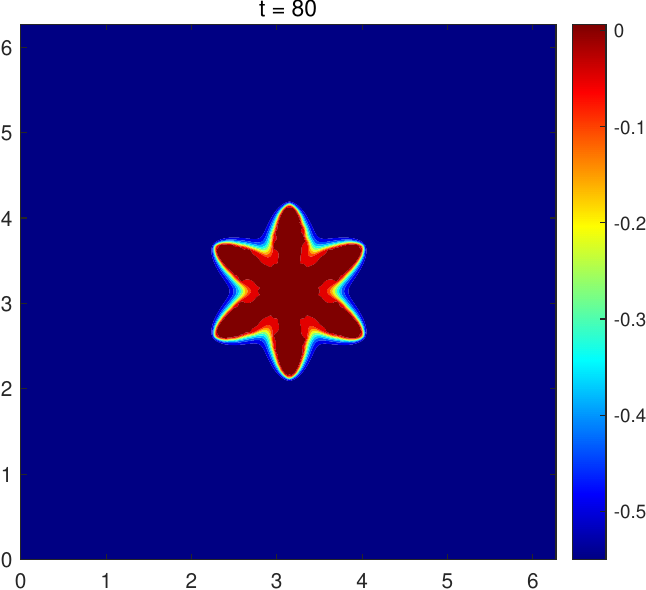}
		\end{minipage}
		\begin{minipage}[t]{0.185\linewidth}
			\centering
			\includegraphics[width=1\linewidth]{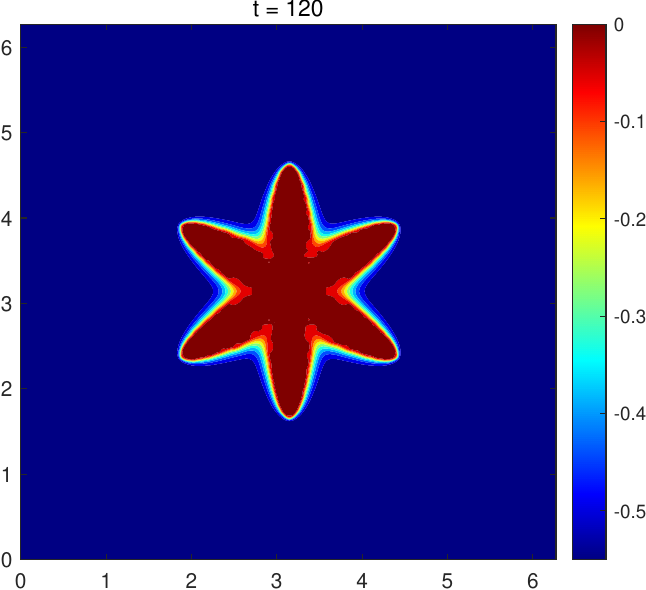}
		\end{minipage}
		\begin{minipage}[t]{0.185\linewidth}
			\centering
			\includegraphics[width=1\linewidth]{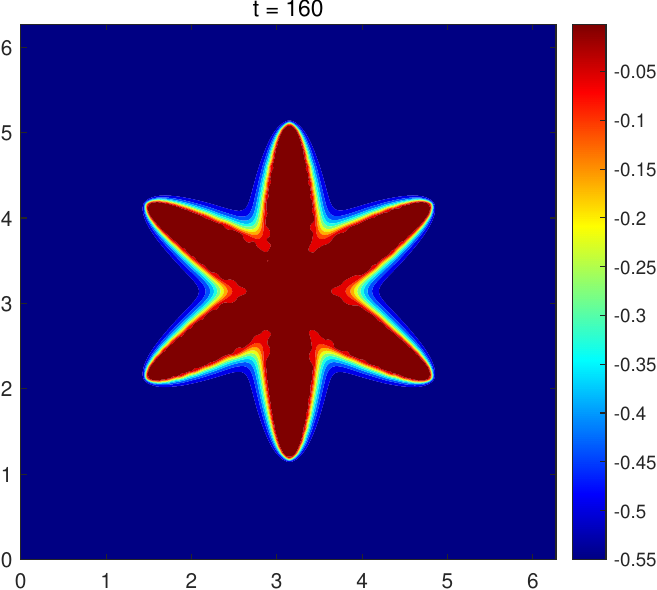}
		\end{minipage}
	}
	\caption{(Example \ref{exp.5}) Dynamics driven by the dendritic crystal growth model with sixfold anisotropy $\beta=6$ and the parameter $K=0.7$ using the adaptive variable time step R-GSAV BDF2 scheme.}
	\label{fig:sixfoldK07}
\end{figure}

\begin{figure}[htbp]
\subfigure[\scriptsize (a) Energy evolution]{
\begin{minipage}{0.45 \textwidth}
\centering
\includegraphics[width=\textwidth]{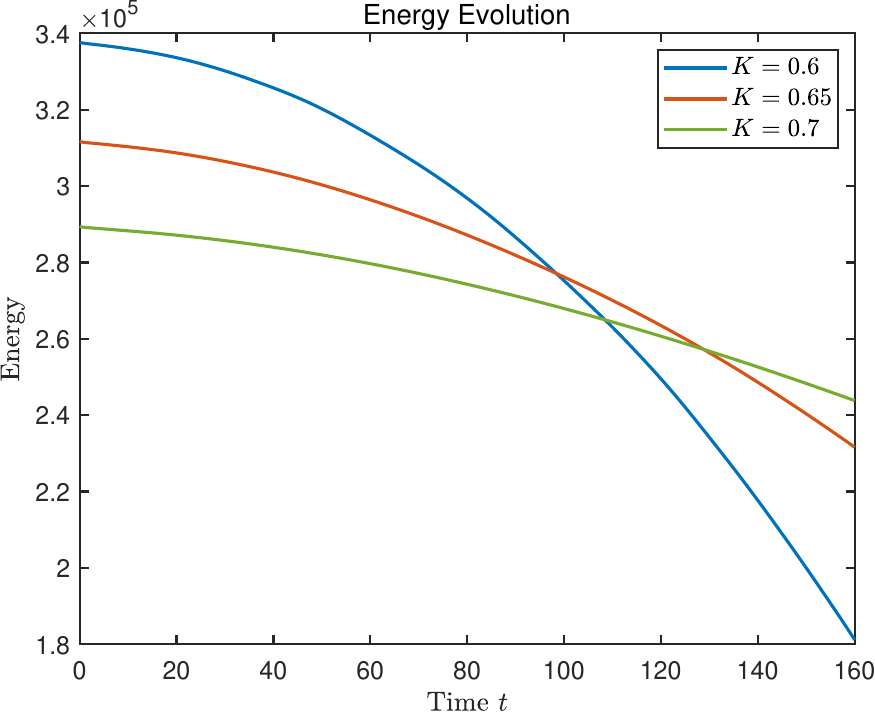}
\end{minipage}
}
\subfigure[\scriptsize (b) Crystal area]{
\begin{minipage}{0.45 \textwidth}
\centering
\includegraphics[width=\textwidth]{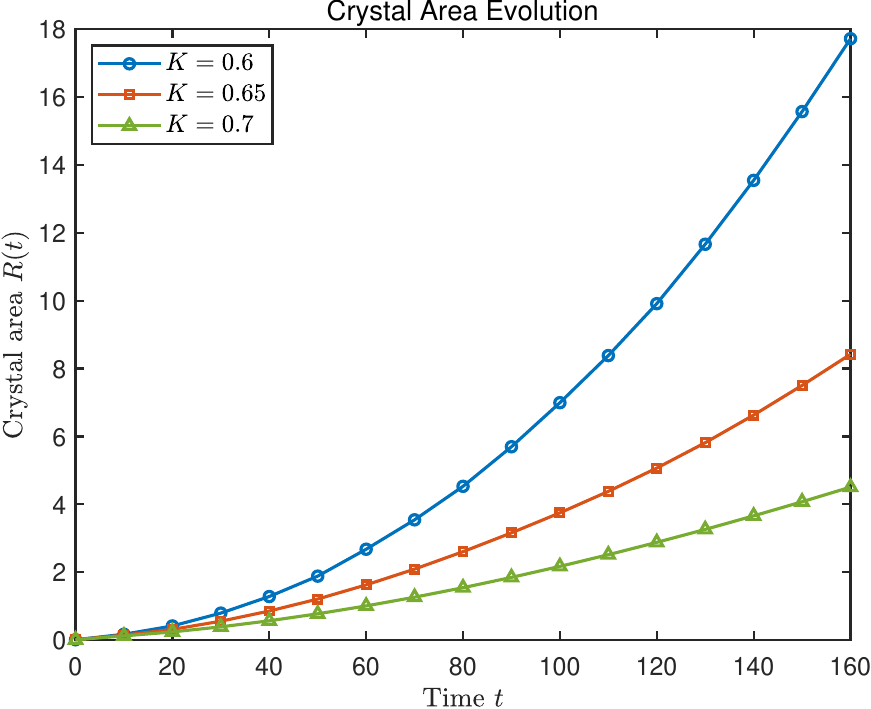}
\end{minipage}
}

\subfigure[\scriptsize (c) Adaptive time step]{
\begin{minipage}{0.45 \textwidth}
\centering
\includegraphics[width=\textwidth]{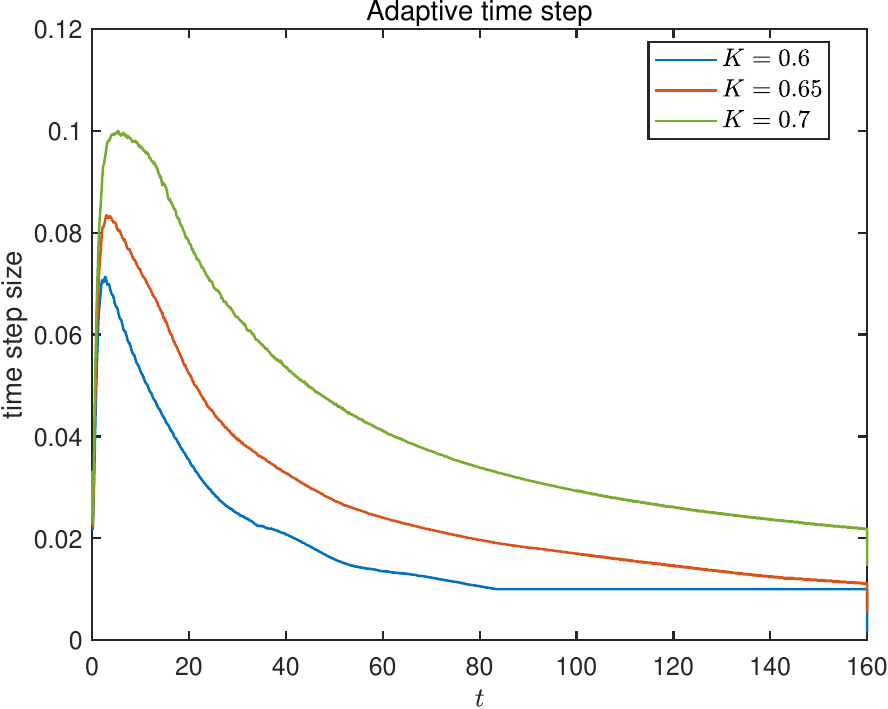}
\end{minipage}
}
\caption{(Example \ref{exp.5}) (a) Time evolution of the energy. (b) The crystal
area $\int_{\Omega}\frac{1+\phi}{2}dx$. (c) The adaptive time step sizes for different values of the parameter $K$ with sixfold anisotropy $\beta=6$. }
\label{sixfold_energy}
\end{figure}

We present in Figures \ref{fig:sixfoldK06}--\ref{fig:sixfoldK07} the snapshots of the phase field $\phi$ and the temperature field $U$ for $K=0.6$, $0.65$ and $0.7$. We observe that the small circle at the initial time first evolves into a hexagon shape, and then into a snowflake with six branches. There are microscopic structures that emerge in the growth process due to the anisotropy in the heat transfer. Comparing the crystal shapes in Figures \ref{fig:sixfoldK06}--\ref{fig:sixfoldK07} for different $K$, we see that the bigger $K$ is, the sharper the snow shaped tips are, the thinner the branches are, and the less pronounced the microstructures are. The evolution of the energy, the crystal area and the adaptive time step sizes  are plotted in Figure \ref{sixfold_energy}, which shows that the energy is dissipative in time, the crystal area increases in time and with decreasing $K$, and the adaptive time step size decreases in time towards $\tau_{\min}$ by the same mechanism as explained in the fourfold case, since the total length of the solid liquid interface and hence the energy decay rate $|E'(t)|$ keep increasing during the snowflake growth. Again, the simulation results are consistent with the reference \cite{yang2020efficient}. The numbers of time steps of the uniform mesh with $\tau=\tau_{\min}=0.01$ and the adaptive strategy \eqref{eq:strategy} are compared in Table \ref{tab:steps44}: the numbers of time steps are reduced by $24.6\%$, $49.4\%$ and $71.5\%$ for $K=0.6$, $0.65$ and $0.7$, respectively.

\begin{table}[htbp]
\centering
\caption{(Example \ref{exp.5}) Comparison of the numbers of time steps up to $T=160$: the uniform mesh with $\tau=0.01$ and the adaptive strategy \eqref{eq:strategy}.}
\label{tab:steps44}
\begin{tabular}{c|c|c|c}
  \hline
 $K$ & uniform $\tau=0.01$ & adaptive strategy \eqref{eq:strategy} & step reduction \\
  \hline
  0.6 & 16000 & 12070 & 24.6\% \\
  0.65 & 16000 & 8104 & 49.4\% \\
  0.7 & 16000 & 4566 & 71.5\% \\
  \hline
\end{tabular}
\end{table}
\subsection{Fourfold dendrite crystal growth in 3D}
Finally, we apply the proposed adaptive variable step R-GSAV BDF2 scheme \eqref{2.4} (with $\delta=1$) to the more challenging three dimensional dendritic crystal growth. The simulation is performed on the domain $\Omega=(0,2\pi)^3$ with $128\times128\times128$ Fourier modes for the spatial discretization, and the time adaptive strategy \eqref{eq:strategy} is adopted  with $\tau_{\min}=10^{-2}$, $\tau_{\max}=1$, $\mu=10^{-2}$ and the initial step size $\tau_0=10^{-2}$.
\begin{example}\label{exp.6}
\rm{The initial phase field is taken as }
\begin{equation}\label{e6}
   \phi(x,y,z,t=0)=\tanh\Big(\frac{0.2-\sqrt{(x-\pi)^2+(y-\pi)^2+(z-\pi)^2}}{0.072}\Big),
\end{equation}
\rm{and the initial temperature field is set}
\begin{equation*}
U(x,y,z,t=0)=\left\{
             \begin{array}{ll}
             0,  &\phi>0,\\
             -0.55, &otherwise.
                           \end{array}
           \right.
\end{equation*}
\rm{The parameters are taken as}
$$
\begin{aligned}
\sigma=2.5e4, \varepsilon=3e-2, \lambda=260, D=2.25e-4, \varepsilon_4=0.05,\beta=4, S_1=4, S_2=4,
\end{aligned}
$$
\rm{and the latent heat parameter is taken as $K=1$ and $K=1.5$.}
\end{example}

\begin{figure}[htbp]
	\vspace{-0.35cm}
	\subfigtopskip=2pt
	\subfigbottomskip=2pt
	\subfigcapskip=-5pt
	\centering	
    \subfigure[\scriptsize (a) Isosurfaces of the phase field $\phi=0$ for $K=1$.]{
	\begin{minipage}[t]{0.185\linewidth}
		\centering
		\includegraphics[width=1\linewidth]{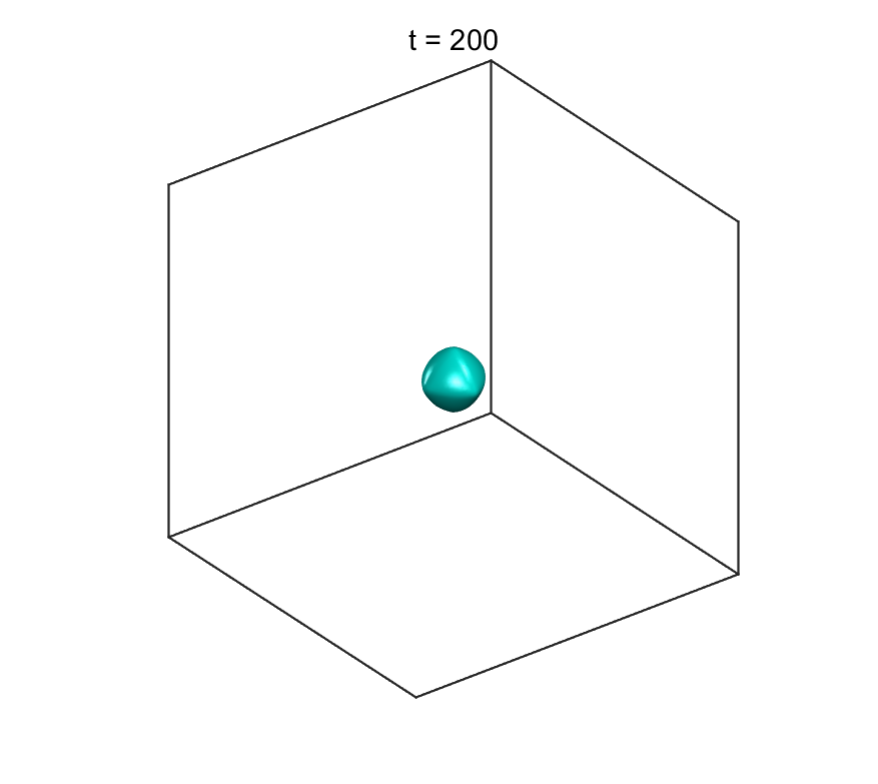}
		\end{minipage}
		\begin{minipage}[t]{0.185\linewidth}
			\centering
			\includegraphics[width=1\linewidth]{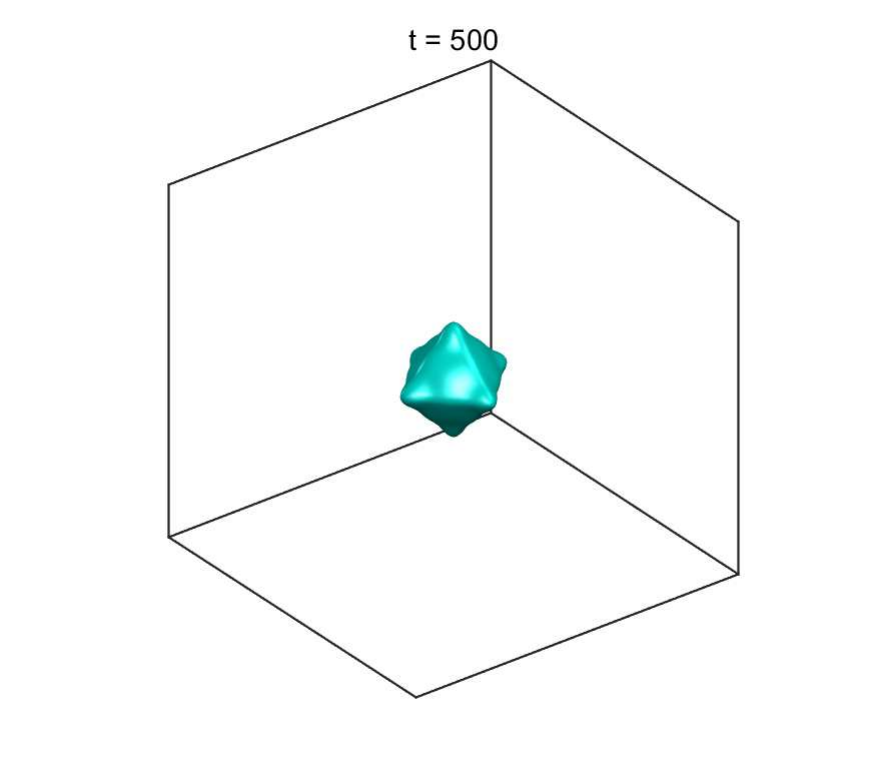}
		\end{minipage}
		\begin{minipage}[t]{0.185\linewidth}
			\centering
			\includegraphics[width=1\linewidth]{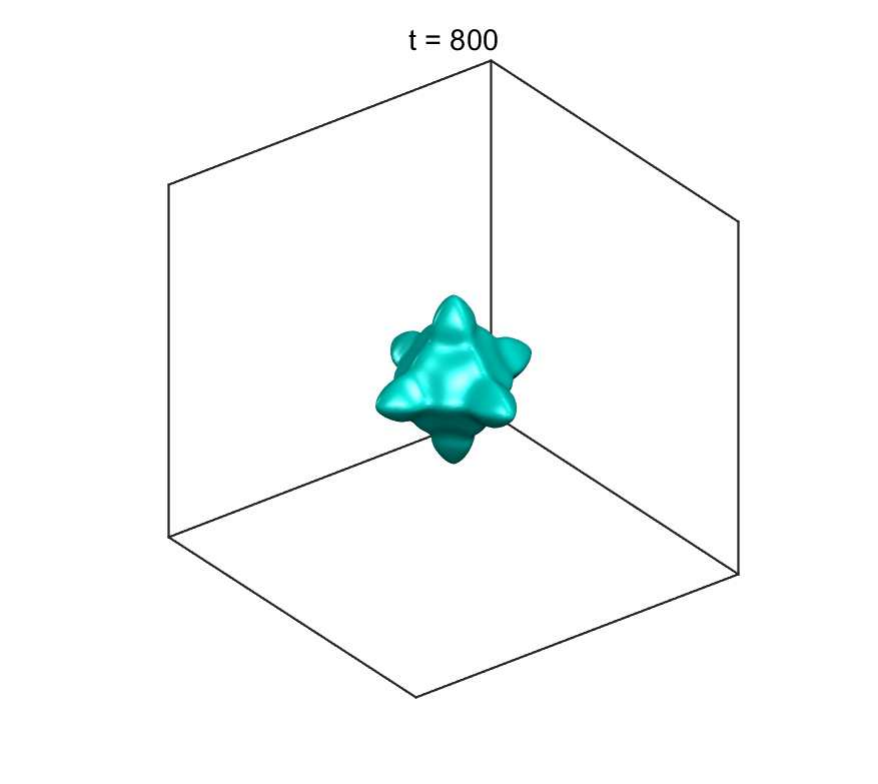}
		\end{minipage}
		\begin{minipage}[t]{0.185\linewidth}
			\centering
			\includegraphics[width=1\linewidth]{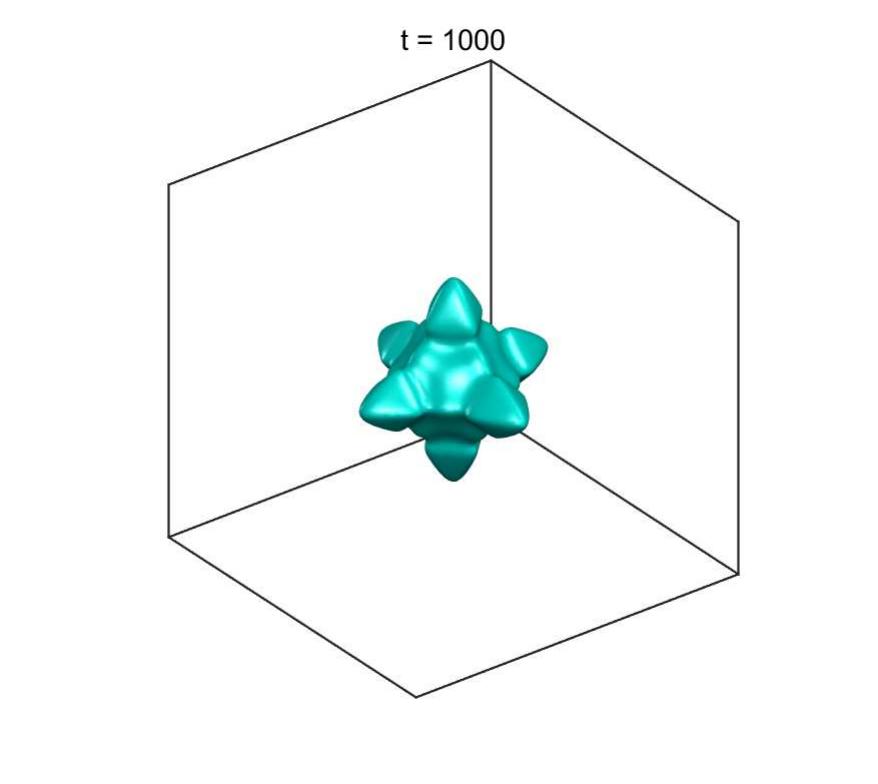}
		\end{minipage}
		\begin{minipage}[t]{0.185\linewidth}
			\centering
			\includegraphics[width=1\linewidth]{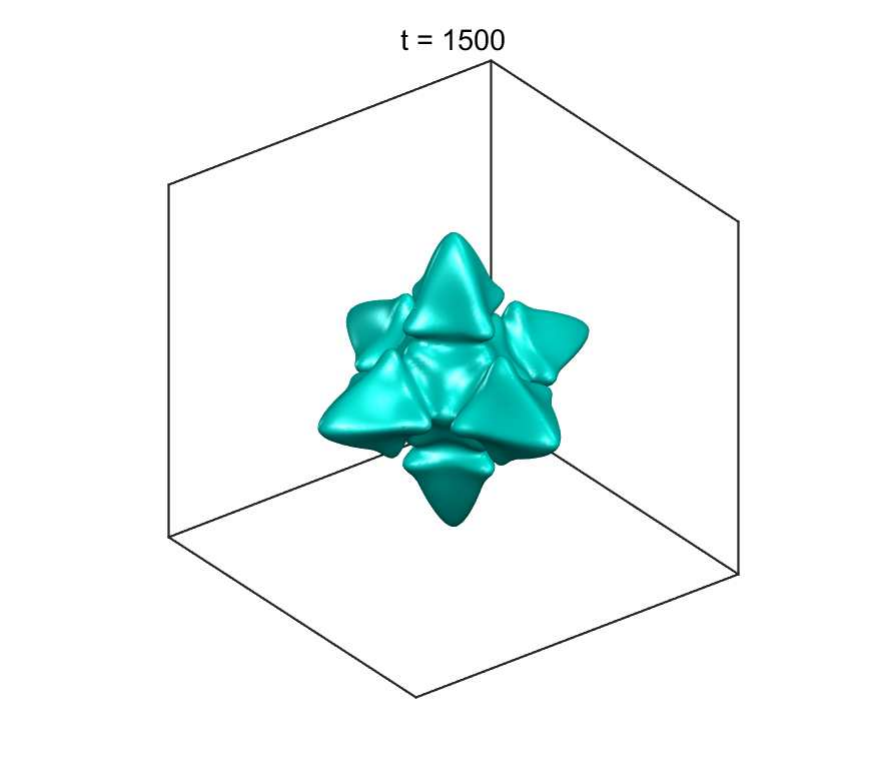}
		\end{minipage}
	}
    \subfigure[\scriptsize (b) Isosurfaces of the phase field $\phi=0$ for $K=1.5$.]{
	\begin{minipage}[t]{0.185\linewidth}
		\centering
		\includegraphics[width=1\linewidth]{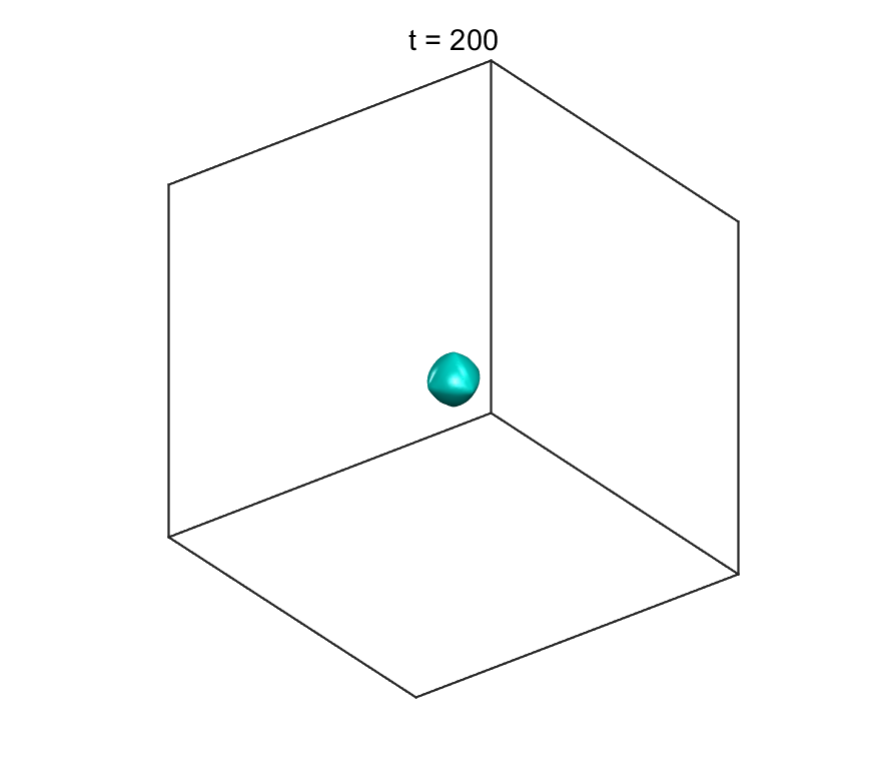}
		\end{minipage}
		\begin{minipage}[t]{0.185\linewidth}
			\centering
			\includegraphics[width=1\linewidth]{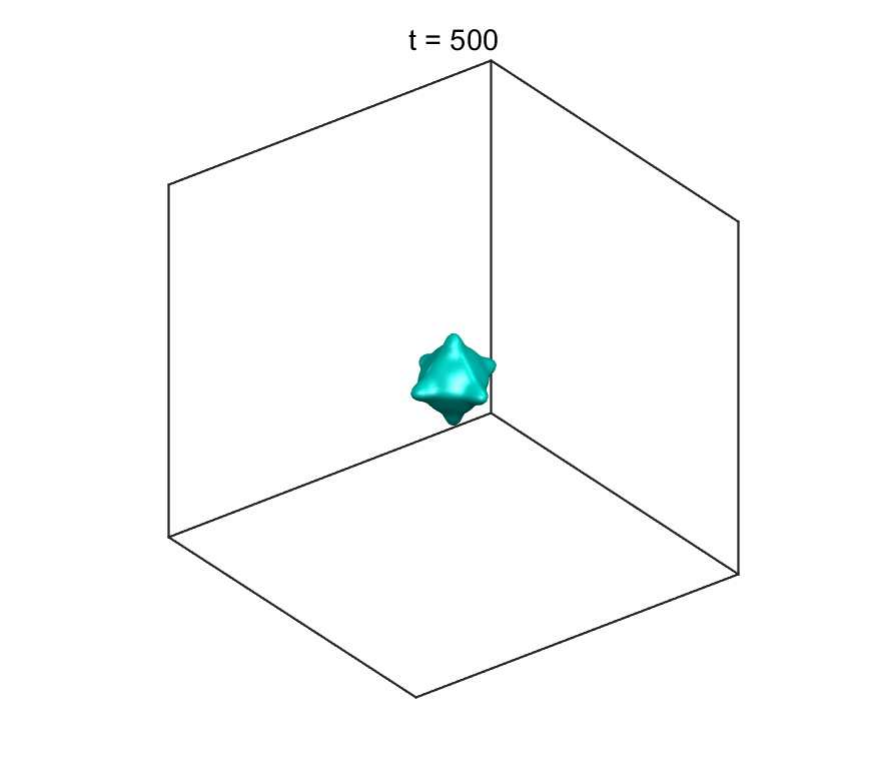}
		\end{minipage}
		\begin{minipage}[t]{0.185\linewidth}
			\centering
			\includegraphics[width=1\linewidth]{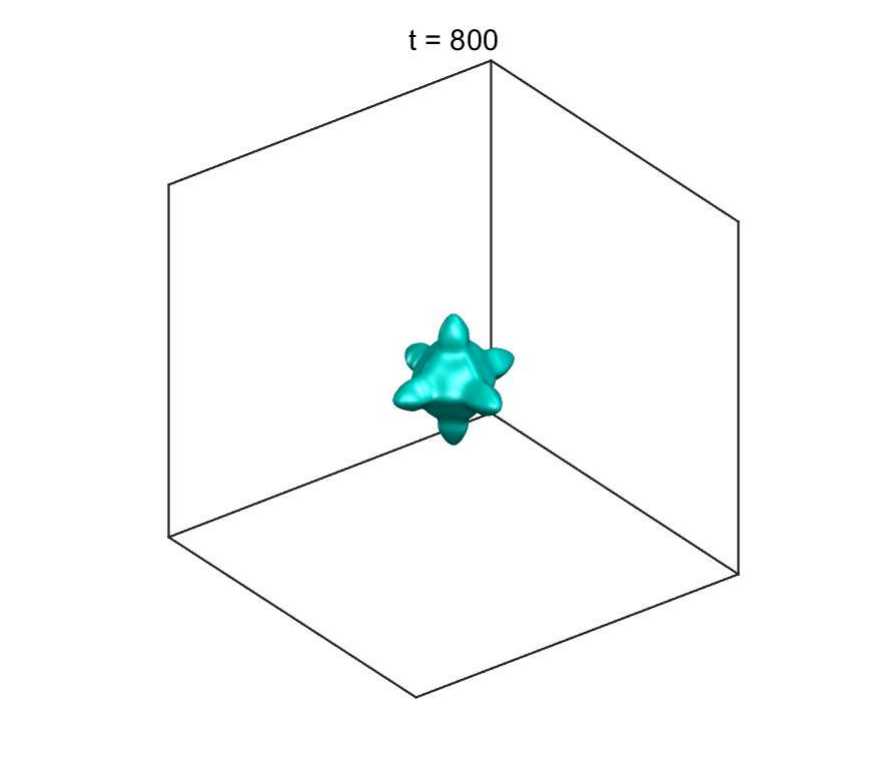}
		\end{minipage}
		\begin{minipage}[t]{0.185\linewidth}
			\centering
			\includegraphics[width=1\linewidth]{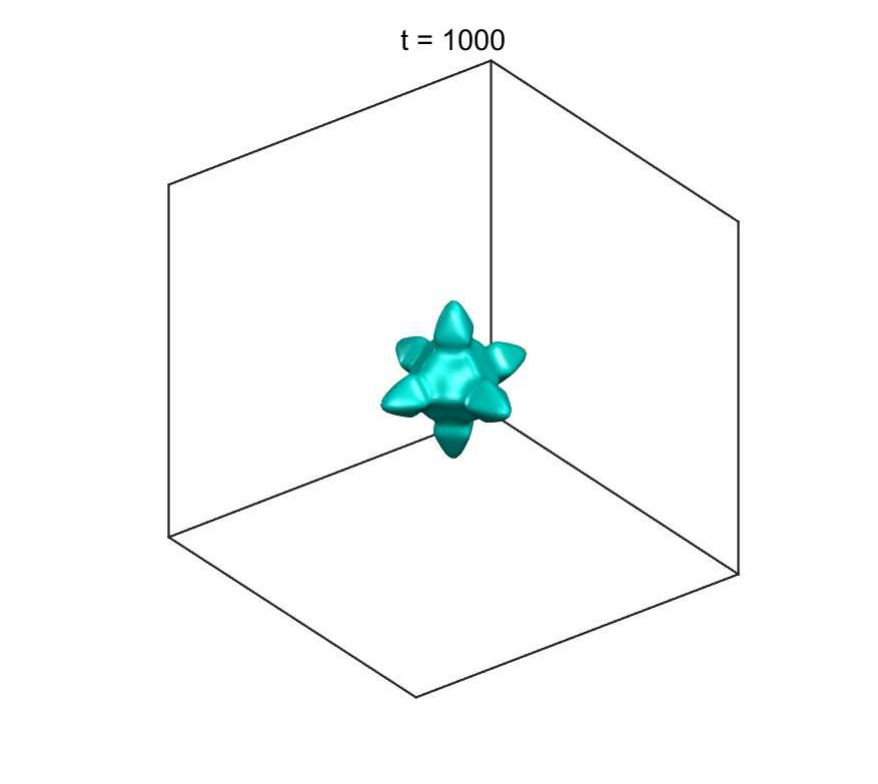}
		\end{minipage}
		\begin{minipage}[t]{0.185\linewidth}
			\centering
			\includegraphics[width=1\linewidth]{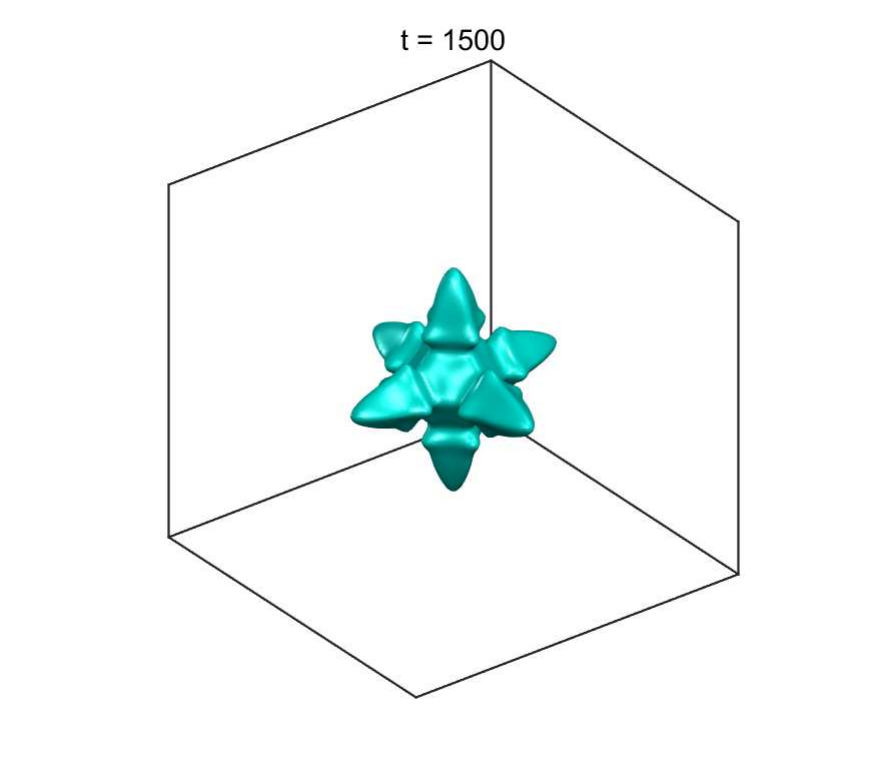}
		\end{minipage}
	}
	\caption{(Example \ref{exp.6}) Three dimensional dendritic crystal growth with fourfold anisotropy $\beta=4$ using the adaptive variable step R-GSAV BDF2 scheme: isosurfaces of the phase field $\phi=0$ for $K=1$ (top row) and $K=1.5$ (bottom row) at $t=200, 500, 800, 1000$ and $1500$.}
	\label{fig:3Dsnap}
\end{figure}

\begin{figure}[htbp]
\subfigure[\scriptsize (a) Energy evolution]{
\begin{minipage}{0.45 \textwidth}
\centering
\includegraphics[width=\textwidth]{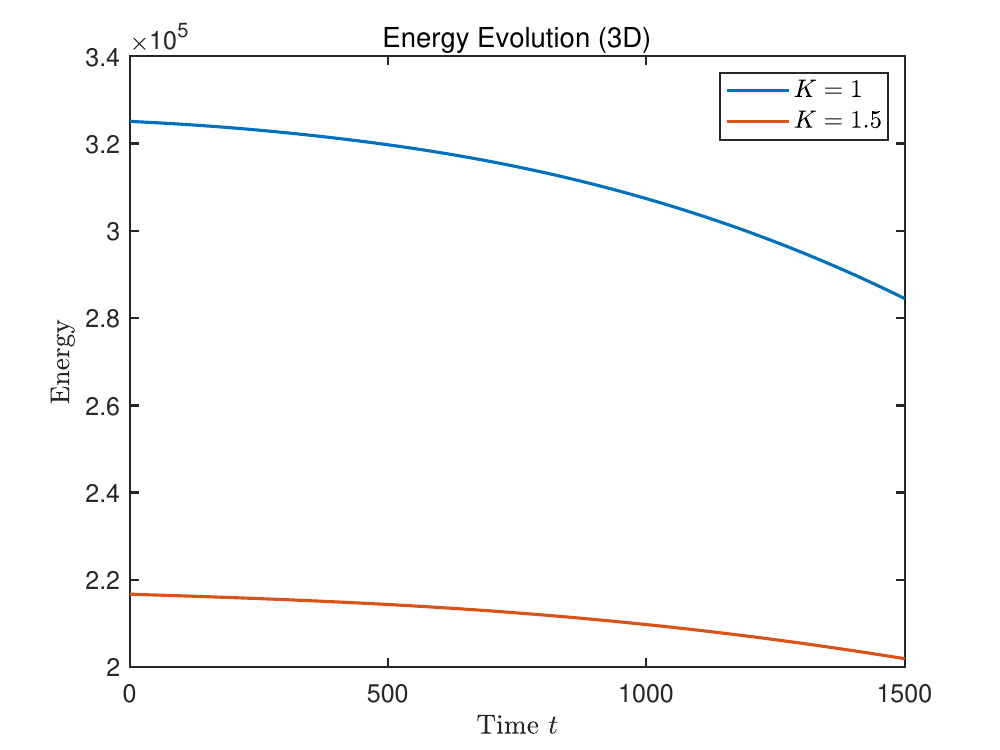}
\end{minipage}
}
\subfigure[\scriptsize (b) Crystal volume]{
\begin{minipage}{0.45 \textwidth}
\centering
\includegraphics[width=\textwidth]{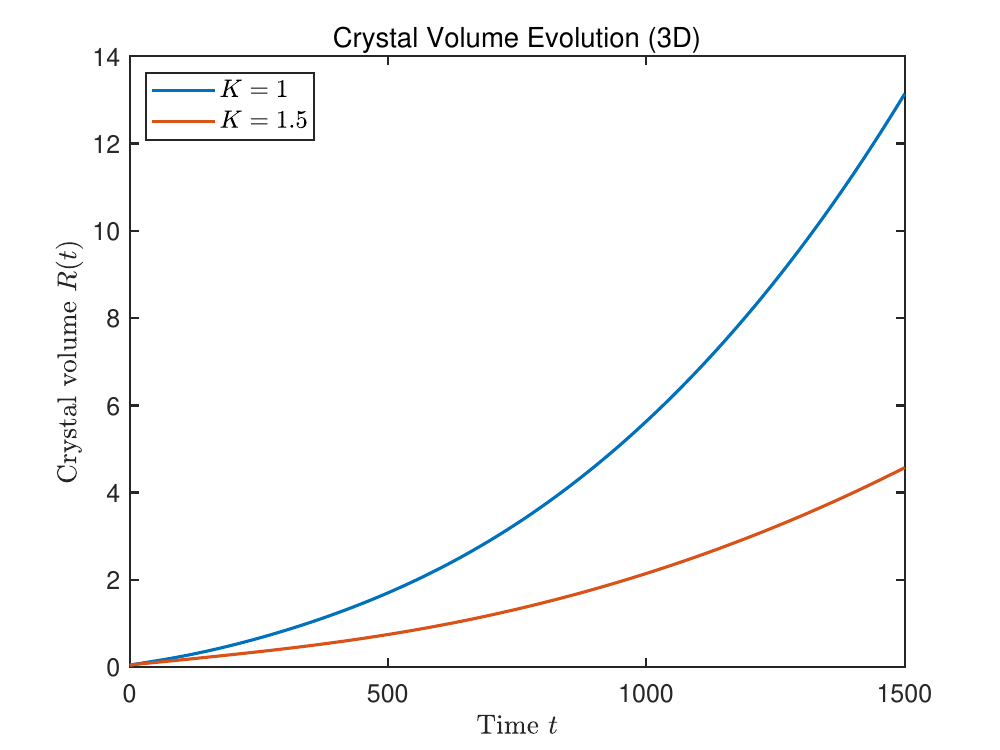}
\end{minipage}
}

\subfigure[\scriptsize (c) Adaptive time step]{
\begin{minipage}{0.45 \textwidth}
\centering
\includegraphics[width=\textwidth]{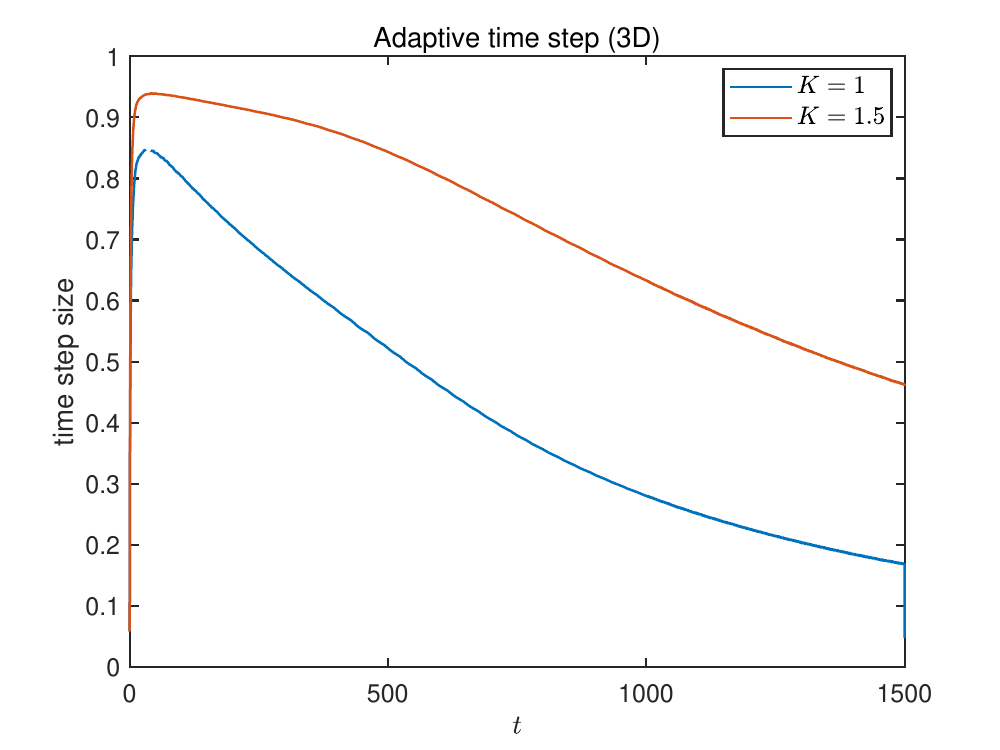}
\end{minipage}
}
\caption{(Example \ref{exp.6})  (a) Time evolution of the energy. (b) The crystal
volume $\int_{\Omega}\frac{1+\phi}{2}d\mathbf{x}$. (c) The adaptive time step sizes for $K=1$ and $K=1.5$ in the three dimensional fourfold case. }
\label{3D_energy}
\end{figure}

Figure \ref{fig:3Dsnap} displays the isosurfaces of the phase field $\phi=0$ at $t=200,500,800,1000$ and $1500$ for $K=1$ and $K=1.5$. It is observed that the tiny spherical seed grows rapidly into a three dimensional dendrite with arms along the six coordinate directions, which is consistent with the fourfold anisotropy of the model. Similar to the two dimensional case, the latent heat parameter $K$ has a significant influence on the growth: the smaller $K$ is, the faster the crystal grows and the thicker the branches are. The time evolution of the energy, the crystal volume and the adaptive time step sizes  are shown in Figure \ref{3D_energy}. The energy is dissipative in time for both values of $K$ (Figure \ref{3D_energy}(a)), and the relaxation keeps the auxiliary variable an accurate approximation of $1$, with $|\xi^{n}-1|$ of order $10^{-9}$--$10^{-8}$ throughout the simulations. The crystal volume increases in time and with decreasing $K$ (Figure \ref{3D_energy}(b)); in particular, the crystal volume at the final time for $K=1$ ($\approx13.1$) is about twice that for $K=1.5$ ($\approx6.6$). As shown in Figure \ref{3D_energy}(c), after a short initial stage, the adaptive time step size decreases gradually in time for both $K$, which, as explained in the two dimensional fourfold case, is because the total area of the solid liquid interface and hence the energy decay rate $|E'(t)|$ keep increasing during the dendrite growth, so that the energy based strategy \eqref{eq:strategy} automatically refines the time step to capture the accelerating growth dynamics. The advantage of the adaptive strategy is also reflected in this 3D simulation: as listed in Table \ref{tab:steps45}, compared with the uniform mesh with $\tau=0.1$, the number of time steps of the adaptive strategy \eqref{eq:strategy} is reduced by $70.5\%$ and $85.5\%$ for $K=1$ and $K=1.5$, respectively.

\begin{table}[htbp]
\centering
\caption{(Example \ref{exp.6}) Comparison of the numbers of time steps: the uniform mesh with $\tau=0.1$ and the adaptive strategy \eqref{eq:strategy}.}
\label{tab:steps45}
\begin{tabular}{c|c|c|c}
  \hline
 $K$ & uniform $\tau=0.1$ & adaptive strategy \eqref{eq:strategy} & step reduction \\
  \hline
  1  & 15000 & 4427 & 70.5\% \\
  1.5  & 15000 & 2171 & 85.5\% \\
  \hline
\end{tabular}
\end{table}
\section{Conclusion}
In this work,  a class of second-order variable time-stepping schemes has been developed for the anisotropic phase field dendritic crystal growth model. The key idea is to combine the R-GSAV approach with a second order BDF operator on
nonuniform temporal meshes, supplemented with two linear stabilizers, so that
the phase field and the temperature field are decoupled and only two linear
equations together with a few algebraic equations need to be solved at each
time step.

 The energy stability of the proposed schemes has been rigorously established for arbitrary positive adjacent step ratios,
while the mild ratio bound retained in the adaptive strategy serves only as a safeguard for temporal accuracy.
Benefiting from the nonuniform mesh formulation, the schemes work seamlessly with energy based time adaptive strategies, whose effectiveness has been verified by the long time simulations.
Moreover, the relaxation technique keeps the auxiliary variable an accurate approximation of $1$ even on irregular temporal meshes, and the computed modified energy coincides with
the original energy in our experiments. The numerical results of the
fourfold and sixfold dendrites in 2D and the fourfold dendrite in 3D agree
well with the existing results in the literature.

In our future work, we will consider extending the present approach to
higher order time stepping schemes on nonuniform temporal meshes and to more
general phase field models.

\section*{CRediT authorship contribution statement}
\textbf{Yayu Guo:} Methodology, Software, Validation, Formal analysis,
Investigation, Writing -- original draft, Visualization.
\textbf{Ziqi Ma:} Conceptualization, Supervision, Project administration,
Writing -- review \& editing.
\section*{Declaration of competing interest}
The authors have no relevant financial or non-financial interests to disclose.
\section*{Data availability statement}
Data will be made available on request.

\clearpage
\bibliographystyle{unsrt}
\bibliography{reference}

\end{document}